\RequirePackage{ifpdf}
\ifpdf 
\documentclass[pdftex]{sigma}
\else
\documentclass{sigma}
\fi

\numberwithin{equation}{section}
\numberwithin{theorem}{section}
\numberwithin{proposition}{section}
\numberwithin{lemma}{section}
\numberwithin{corollary}{section}
\numberwithin{conjecture}{section}
\numberwithin{definition}{section}
\numberwithin{example}{section}
\numberwithin{remark}{section}

\begin{document}

\allowdisplaybreaks

\renewcommand{\thefootnote}{$\star$}

\renewcommand{\PaperNumber}{013}

\FirstPageHeading

\ShortArticleName{Soliton Cellular Automaton}

\ArticleName{Bethe Ansatz, Inverse Scattering Transform\\
and Tropical Riemann Theta Function\\ in a Periodic
Soliton Cellular Automaton for $\boldsymbol{A^{(1)}_n}$\footnote{This paper is a
contribution to the Proceedings of the Workshop ``Geometric Aspects of Discrete and Ultra-Discrete Integrable Systems'' (March 30 -- April 3, 2009, University of Glasgow, UK). The full collection is
available at
\href{http://www.emis.de/journals/SIGMA/GADUDIS2009.html}{http://www.emis.de/journals/SIGMA/GADUDIS2009.html}}}

\Author{Atsuo KUNIBA~$^\dag$ and Taichiro TAKAGI~$^\ddag$}

\AuthorNameForHeading{A.~Kuniba and T.~Takagi}

\Address{$^\dag$~Institute of Physics, University of Tokyo, Komaba,
Tokyo 153-8902, Japan}
\EmailD{\href{mailto:atsuo@gokutan.c.u-tokyo.ac.jp}{atsuo@gokutan.c.u-tokyo.ac.jp}}

\Address{$^\ddag$~Department of Applied Physics, National Defense Academy,
Kanagawa 239-8686, Japan}
\EmailD{\href{mailto:takagi@nda.ac.jp}{takagi@nda.ac.jp}}

\ArticleDates{Received September 21, 2009;  Published online January 31, 2010}

\Abstract{We study an integrable
vertex model with a periodic boundary condition
associa\-ted with $U_q\big(A^{(1)}_n\big)$ at the crystallizing point $q=0$.
It is an $(n+1)$-state cellular automaton
describing the factorized scattering of solitons.
The dynamics originates in the
commuting family of fusion transfer matrices and generalizes
the ultradiscrete Toda/KP f\/low corresponding to
the periodic box-ball system.
Combining Bethe ansatz and crystal theory in quantum group,
we develop an inverse scattering/spectral formalism
and solve the initial value problem
based on several conjectures.
The action-angle variables are constructed representing
the amplitudes and phases of solitons.
By the direct and inverse scattering maps,
separation of variables into solitons is achieved and
nonlinear dynamics is transformed into a~straight motion
on a tropical analogue of the Jacobi variety.
We decompose the level set into connected components
under the commuting family of time evolutions
and identify each of them with the set of integer points on a torus.
The weight multiplicity formula derived from the $q=0$
Bethe equation acquires an elegant interpretation as the volume
of the phase space expressed by the size and multiplicity of these tori.
The dynamical period is determined as an explicit
arithmetical function of the $n$-tuple of Young diagrams
specifying the level set.
The inverse map, i.e.,
tropical Jacobi inversion is expressed in terms of a tropical
Riemann theta function associated with the Bethe ansatz data.
As an application, time average of some local variable is
calculated.}


\Classification{82B23; 37K15; 68R15; 37B15}


\tableofcontents

\renewcommand{\thefootnote}{\arabic{footnote}}
\setcounter{footnote}{0}

\section{Introduction}\label{k:sec:intro}

\subsection{Background}
In \cite{KT1}, a class of soliton cellular automata (SCA)
on a periodic one dimensional
lattice was introduced associated with
the non-exceptional quantum af\/f\/ine algebras
$U_q({\mathfrak g}_n)$ at $q=0$.
In this paper we focus on the type ${\mathfrak g}_n = A^{(1)}_n$.
The case $n=1$ was introduced earlier as
the periodic box-ball system \cite{YT}.
It is a periodic version of Takahashi--Satsuma's
soliton cellular automaton~\cite{TS} originally def\/ined on the inf\/inite lattice.
In the pioneering work~\cite{YYT},
a closed formula for the dynamical period
of the periodic box-ball system was found
for the states with a trivial internal symmetry.

After some while, the f\/irst complete solution of the initial value problem
and its explicit formula by a tropical analogue of Riemann theta
function were obtained in \cite{KTT,KS1,KS2}.
This was done by developing the
inverse scattering/spectral method\footnote{Although it is more customary to use `spectral' in periodic systems,
we shall simply call it inverse scattering hereafter for short.}~\cite{GGKM}
in a tropical (ultradiscrete) setting
primarily based on the {\em quantum integrability} of the system.
It provides a linearization
scheme of the nonlinear dynamics into straight motions.
The outcome is practical as well as conceptual.
For instance practically,
a closed formula of the exact dynamical period
for {\em arbitrary} states \cite[Theorem 4.9]{KTT}
follows straightforwardly without a cumbersome
combinatorial argument, see~(\ref{nhf}).
Perhaps what is more important conceptually is,
it manifested a decent tropical analogue
of the theory of quasi-periodic soliton solutions
\cite{DT,DMN} and its quantum aspects in the realm of SCA.
In fact, several essential tools in classical
integrable systems have found their quantum counterparts
at a combinatorial level.
Here is a table of rough correspondence.
\begin{center}
\begin{tabular}{c|c}
classical  &  quantum \\
\hline
solitons  & Bethe strings \\
discrete KP/Toda f\/low &
fusion transfer matrix \\
action-angle variable & extended rigged conf\/iguration mod
($q=0$ Bethe eq.)\\
Abel--Jacobi map,
Jacobi inversion& modif\/ied Kerov--Kirillov--Reshetikhin bijection \\
Riemann theta function & charge of extended rigged conf\/iguration
\end{tabular}
\end{center}
\noindent
Here the Kerov--Kirillov--Reshetikhin (KKR) bijection \cite{KKR,KR}
is the celebrated algorithm on rigged conf\/igurations
in the combinatorial Bethe ansatz
explained in Appendix \ref{k:app:kkr}.
As for the last line, see (\ref{k:star}).

In a more recent development \cite{T},
a complete decomposition
of the phase space of the periodic box-ball system has been
accomplished into connected components under the time evolutions,
and each of them is identif\/ied with the set
of integer points on a certain torus ${\mathbb Z}^g/F{\mathbb Z}^g$ explicitly.

\subsection{A quick exposition}
The aim of this paper is to extend all
these results \cite{KTT,KS1,KS2,T}
to a higher rank case based on several conjectures.
The results will be demonstrated with a number of examples.
We call the system the
{\em periodic $A^{(1)}_n$ soliton cellular automaton $($SCA$)$}.
It is a dynamical system on the~$L$~numbers
$\{1,2,\ldots, n+1\}^L$ which we call {\em paths}.
There is the $n$-tuple of commuting series of
time evolutions $T^{(1)}_l, \ldots, T^{(n)}_l$ with
$l \ge 1$.
The system is periodic in that
the cyclic shift~(\ref{k:shf}) is contained in the
member as $T^{(1)}_1$.
The casual exposition in the following example may be helpful
to grasp the idea quickly.

\begin{example}\label{k:ex:int}
Take the path $p=321113211222111223331111$
of length $L=24$ with $n=2$,
which will further be treated in
Examples \ref{k:ex:arc}, \ref{k:ex:aa2}, \ref{k:ex:ga2}
\ref{k:ex:deco} and \ref{ex:owari}.
\begin{equation*}
\begin{split}
&\qquad\qquad\; \big(T^{(1)}_3\big)^t(p)\qquad \qquad\qquad\qquad\qquad\qquad
\big(T^{(2)}_4\big)^t(p)\\
t=0:&\quad 321113211222111223331111 \qquad\qquad 321113211222111223331111 \\
t=1:&\quad 113211132111222112213331 \qquad\qquad 132111321333111221112221 \\
t=2:&\quad 331132111321111221122213 \qquad\qquad 112211132111222331113321 \\
t=3:&\quad 213311332113211112211122 \qquad\qquad 213311113211322112211132 \\
t=4:&\quad 122233111332132111122111 \qquad\qquad 221122111321133213311112 \\
t=5:&\quad 111122332111321332111221 \qquad\qquad 331132111122111321122213 \\
t=6:&\quad 211111221332113211332112 \qquad\qquad 112213221133111132132211 \\
t=7:&\quad 122211112211332132111331 \qquad\qquad 113311222111221113213321 \\
t=8:&\quad 311122211122111321332113 \qquad\qquad 211122333111321111221132 \\
t=9:&\quad 133111122211221113211332 \qquad\qquad 221132111222132111331113
\end{split}
\end{equation*}
By regarding $1$ as an empty background, these patterns
show the repeated collisions of
solitons like 22, 32, 33, 222, 322, 332 and 333.
At each time step there are always 5 solitons with
amplitudes 3, 3, 2, 2, 2.
The dynamics is highly nonlinear.
Ultimately in our approach, it will be
transformed into the straight motion
(Hint: right hand side of (\ref{k:rc35}) and Remark \ref{k:re:t11}):
\begin{equation*}
\begin{picture}(200,75)(-170,13)
\setlength{\unitlength}{0.5mm}

\multiput(-30,0)(0,0){1}{

\put(-16,42){$\Phi$}
\put(-84,35){$\big(T^{(1)}_l\big)^t\big(T^{(2)}_m\big)^s(p)\quad\longmapsto$}

\multiput(0,0)(0,-10){3}{
\put(10,60){\line(1,0){30}}}
\multiput(0,0)(0,10){4}{
\put(10,10){\line(1,0){20}}}
\multiput(0,0)(10,0){3}{
\put(10,60){\line(0,-1){50}}}
\put(40,60){\line(0,-1){20}}

\put(42,52){$1+\min(l,3)t$}
\put(42,42){$0+\min(l,3)t$}
\put(32,32){$7+\min(l,2)t$}
\put(32,22){$4+\min(l,2)t$}
\put(32,12){$2+\min(l,2)t$}
}

\multiput(-30,0)(0,0){1}{
\put(143,52){$11+\min(m,4)s$}
\put(113,42){$1+\min(m,1)s$}
\put(100,60){\line(1,0){40}}
\put(100,50){\line(1,0){40}}
\put(100,40){\line(1,0){10}}
\put(100,60){\line(0,-1){20}}
\put(110,60){\line(0,-1){20}}
\multiput(0,0)(10,0){3}{
\put(120,60){\line(0,-1){10}}}}
\end{picture}
\end{equation*}
Here $\Phi$ is the modif\/ied Kerov--Kirillov--Reshetikhin (KKR) bijection
mentioned in the table.
One notices that
$T^{(1)}_3 = T^{(1)}_\infty$
and $T^{(2)}_4 = T^{(4)}_\infty$ in this example.
(Of course $\min(m,1)=1$.)
Note also that the f\/irst Young diagram $(33222)$ gives the
list of amplitudes of solitons.
The two Young diagrams are conserved quantities (action variable).
A part of a Young diagram having equal width is called a {\em block}.
There are 4 blocks containing 2, 3, 1, 1 rows in this example.
The 7 dimensional vector consisting of the listed numbers
is the angle variable f\/lowing linearly.
It should be understood as an element of (Hint:~(\ref{k:waru}))
\begin{equation*}
{\mathbb Z}^7/B{\mathbb Z}^7/({\frak S}_2 \times {\frak S}_3 \times
{\frak S}_1 \times {\frak S}_1),
\end{equation*}
where ${\frak S}_m$ is the symmetric group of degree $m$.
We have exhibited the trivial ${\frak S}_1$'s as another hint.
The matrix $B$ is given by (\ref{k:Bex}).
The description is simplif\/ied further if one realizes that the
relative (i.e., dif\/ferences of)
components within each block remain unchanged throughout.
Thus picking the bottom component only from each block, we obtain
\begin{equation*}
\big(T^{(1)}_l\big)^t\big(T^{(2)}_m\big)^s(p)
\overset{\Phi_\chi}{\longmapsto}
\begin{pmatrix}
0+\min(l,3)t\\
2+\min(l,2)t\\
11+\min(m,4)s\\
1+\min(m,1)s
\end{pmatrix}
\mod F{\mathbb Z}^4,\qquad
F = \begin{pmatrix}
16 & 12 & -3 & -1 \\
8 & 19 & -2 & -1 \\
-6 & -6 & 10 & 2 \\
-2 & -3 & 2 & 3
\end{pmatrix}.
\end{equation*}
The 4 dimensional vector here is the reduced angle variable
(Section \ref{k:sec:dcc}).
After all, the dynamics is just a straight motion
in ${\mathbb Z}^4/F{\mathbb Z}^4$.
The abelian group of the time evolutions acts transitively on this torus.
Thus the size of the connected component of $p$ is $\det F = 4656$.
The dynamical period ${\mathcal N}^{(r)}_l$ of $T^{(r)}_l$, i.e.,
the smallest positive integer satisfying
$(T^{(r)}_l)^{{\mathcal N}^{(r)}_l}(p)=p$ is just
the smallest positive integer satisfying
\begin{equation*}
{\mathcal N}^{(1)}_l
\begin{pmatrix}
\min(l,3)\\
\min(l,2)\\
0\\
0
\end{pmatrix}\in F{\mathbb Z}^4,\qquad
{\mathcal N}^{(2)}_m
\begin{pmatrix}
0\\
0\\
\min(m,4)\\
\min(m,1)
\end{pmatrix}\in F{\mathbb Z}^4,
\end{equation*}
where the coef\/f\/icient vectors are velocities.
Solving the elementary exercise, one f\/inds
\[
\mbox{\begin{tabular}{c|c|c|c|c|c|c}
${\mathcal N}^{(1)}_1$ & ${\mathcal N}^{(1)}_2$
& ${\mathcal N}^{(1)}_\infty$
& ${\mathcal N}^{(2)}_1$ & ${\mathcal N}^{(2)}_2$
& ${\mathcal N}^{(2)}_3$ & ${\mathcal N}^{(2)}_\infty$  \\
\hline
24 & 12 & 194 & 1164 & 776 & 582 & 2328
\end{tabular}}
\]
The level set characterized by the above pair of Young diagrams
$((33222),(41))$ consists of 139680 paths.
It is decomposed into 36 connected components as
(see Example \ref{k:ex:deco})
\begin{equation*}
24\big(\mathbb{Z}^4 /F \mathbb{Z}^4\big)
\sqcup 12\big(\mathbb{Z}^4 /F' \mathbb{Z}^4\big),\qquad
139680
=24(\det F)+12(\det F'),
\end{equation*}
where $F' = F\cdot\text{diag}(\frac{1}{2},1,1,1)$.
Here the two kinds of tori arise ref\/lecting that the two kinds of
internal symmetries are allowed for the
pair $((33222),(41))$ under consideration.
\end{example}

\subsection{Related works}
A few remarks are in order on related works.

1.  The solution of initial value problem of
the $n=1$ periodic box-ball system \cite{KTT,KS1,KS2}
has been reproduced partially by
the procedure called 10-elimination \cite{MIT2}.
By now the precise relation between
the two approaches has been shown \cite{KiS}.
It is not known whether the 10-elimination admits a decent generalization.
The KKR bijection on the other hand is a canonical algorithm
allowing generalizations
not only to type $A^{(1)}_n$ \cite{KSS}
but also beyond \cite{OSS}.
Conjecturally the approach based on the KKR type bijection
will work universally for the periodic SCA~\cite{KT1}
associated with Kirillov--Reshetikhin crystals
as exemplif\/ied for $A^{(1)}_1$
higher spin case \cite{KS3} and~$A^{(1)}_n$ \cite[this paper]{KT1, KT2}.
In fact, the essential results like torus, dynamical period and
phase space volume formula are all presented in a~universal form by the data from Bethe ansatz and the
root system\footnote{It is a Bethe ansatz folklore that
``the solution exists before the model is constructed''
according to V.V.~Bazhanov.}.
From a practical point of view,
it should also be recognized that
the KKR map $\phi^{\pm1}$ is
a delightfully elementary algorithm described in less than half a page
in our setting in Appendix \ref{k:app:kkr}.
It is a useful exercise to program it to follow examples in this paper.

2. When there is no duplication of the amplitudes of solitons,
the tropical analogue of the Jacobian ${\mathcal J}(\mu)$
obtained in \cite{KTT} has an interpretation from
the tropical geometry point of view~\cite{IT0,MZ}.
In this paper we are naturally led to a higher rank version of
${\mathcal J}(\mu)$ and the relevant tropical analogue
of the Riemann theta function.
We will call them `tropical $\ldots$' rather casually without
identifying an underlying tropical geometric objects hoping not to cause
a too much embarrassment.

3. In \cite{MIT1}, the dynamical system on
$B^{1,l_1}\otimes \cdots \otimes B^{1,l_L}$
equipped with the unique time evolu\-tion~$T^{(1)}_\infty$ of the form (\ref{k:tkk}) was studied under
the name of the generalized periodic box-ball system.
See Sections \ref{k:sec:rect} and
\ref{k:sec:DD} for the notations $B^{r,l}$ and $T^{(r)}_l$.
Approaches by ultradiscretization of the periodic Toda lattice
also capture the $T^{(1)}_\infty$ but conceivably
$\{T^{(1)}_l\mid  l\ge 1\}$ at most.
Our periodic $A^{(1)}_n$ SCA is the generalization of
the generalized periodic box-ball system
with $\forall\, l_i = 1$ that is furnished with
the commuting family of time evolutions
$\{T^{(r)}_l\,|\,1 \le r \le n, l\ge 1\}$.
As we will see in the rest of the paper, it is crucial to consider
this wider variety of dynamics and
their whole joint spectrum $\{E^{(r)}_l\}$.
They turn out to be the necessary and suf\/f\/icient conserved quantities
to formulate the inverse scattering method.
Such a usage of the full family $\{T^{(r)}_l\}$
was f\/irstly proposed in the periodic setting in \cite{KT1,KT2}.
In particular in the latter reference, the most general
periodic $A^{(1)}_n$
SCA on $B^{r_1,l_1}\otimes \cdots \otimes B^{r_1,l_L}$
endowed with the dynamics $\{T^{(r)}_l\}$
was investigated, and the dynamical period and
a phase space volume formula were conjectured
using some heuristic connection with combinatorial Bethe ansatz.
This paper concerns  the basic case $\forall \, r_k = \forall \, l_k=1$ only
but goes deeper to explore the linearization scheme under which
the earlier conjectures \cite{KT1,KT2}
get ref\/ined and become simple corollaries.
We expect that essential features in the inverse scattering formalism
are not too much inf\/luenced by the choice of $\{r_k, l_k\}$.

\subsection{Contents of paper}
Let us digest the contents of the paper
along Sections \ref{k:sec:model}--\ref{k:sec:trt}.

In Section \ref{k:sec:model}, we formulate the periodic $A^{(1)}_n$ SCA.
It is a solvable vertex model \cite{Ba}
associated with the quantum af\/f\/ine algebra
$U_q(A^{(1)}_n)$ at $q=0$ in the sense that notions concerning
$U_q(A^{(1)}_n)$ are replaced by those from
the crystal theory \cite{K, KMN, KMN2},
a theory of quantum group at $q=0$.
The correspondence is shown between the
f\/irst and the second columns of the following table,
where we hope the third column may be more friendly to the readers
not necessarily familiar with the crystal theory.

\begin{center}
\begin{tabular}{@{}c|c|c@{}}
$U_q(A^{(1)}_n)$ vertex model & SCA $(q=0)$ & combinatorial description\\
\hline
Kirillov--Reshetikhin module $W^{(r)}_l$ &
crystal $B^{r,l}$ & shape $r\times l$ semistandard tableaux \\
quantum R & combinatorial R
&rule (\ref{k:bccb}) on Schensted insertions \\
fusion transfer matrix & time evolution $T^{(r)}_l$
& diagram (\ref{k:hone}) with $v=v' \in B^{r,l}$\\
conserved quantity & energy $E^{(r)}_l$
& $n$-tuple of Young diagrams $\mu$\\
\end{tabular}
\end{center}

The def\/inition of the Kirillov--Reshetikhin module
was f\/irstly given in \cite[Def\/inition 1.1]{KN}, although we do not use
it in this paper.
As mentioned previously,
there is the $n$-tuple of commuting series of time evolutions
$T^{(1)}_l, \ldots, T^{(n)}_l$ with $l \ge 1$.
The f\/irst series $T^{(1)}_l$ is
the ultradiscrete Toda/KP f\/low \cite{HHIKTT, Yy}, and
especially its top $T^{(1)}_\infty$
admits a simple description by a ball-moving algorithm
like the periodic box-ball system.
See Theorem \ref{k:th:fac} and comments following it.
However, as again said previously,
what is more essential in our approach is
to make use of the entire family  $\big\{T^{(r)}_l\big\}$ and
the associated conserved quantities called energy~$\big\{E^{(r)}_l\big\}$.
It is the totality of the joint spectrum $\big\{E^{(r)}_l\big\}$
that makes it possible to characterize
the level set~${\mathcal P}(\mu)$~(\ref{k:pmu})
by the $n$-tuple of Young diagrams
and the whole development thereafter.

As an illustration,
$T^{(2)}_2(241123431) = 423141213$
is calculated as
\begin{equation*}
\mbox{\footnotesize
\begin{picture}(200,30)(-60,-7)

\put(18,14){2}\put(48,14){4}\put(78,14){1}
\put(108,14){1}\put(138,14){2}\put(168,14){3}
\put(198,14){4}\put(228,14){3}\put(258,14){1}

\multiput(20,3)(30,0){9}{
\put(-6,0){\line(1,0){12}}\put(0,-8){\line(0,1){16}}}

\put(0,3.5){12}\put(0,-3.5){34}
\put(30,3.5){12}\put(30,-3.5){23}
\put(60,3.5){13}\put(60,-3.5){24}
\put(90,3.5){11}\put(90,-3.5){24}
\put(120,3.5){11}\put(120,-3.5){24}
\put(150,3.5){11}\put(150,-3.5){22}
\put(180,3.5){12}\put(180,-3.5){23}
\put(210,3.5){13}\put(210,-3.5){24}
\put(240,3.5){23}\put(240,-3.5){34}
\put(270,3.5){12}\put(270,-3.5){34}

\put(18,-14){4}\put(48,-14){2}\put(78,-14){3}
\put(108,-14){1}\put(138,-14){4}\put(168,-14){1}
\put(198,-14){2}\put(228,-14){1}\put(258,-14){3}

\end{picture}}
\end{equation*}
by using the `hidden variable' called carrier on horizontal edges
belonging to $B^{2,2}$.
This is a~conventional diagram representing
the row transfer matrix of a vertex model~\cite{Ba} whose
auxiliary (horizontal) space has the `fusion type' $B^{2,2}$.
The general case~$T^{(r)}_l$ is similarly def\/ined by using~$B^{r,l}$.
The peculiarity as a vertex model is that there is no thermal
f\/luctuation due to the crystallizing choice $q=0$
resulting in a deterministic dynamics.
Note that the carrier has been chosen
specially so that the leftmost and the rightmost
ones coincide to match the periodic boundary condition.
This non-local
postulate makes the well-def\/inedness of the dynamics highly nontrivial
and is in fact a source of the most intriguing features of our SCA.
It is an important problem to characterize the situation in which
all the time evolutions act stably.
With regard to this, we propose a neat suf\/f\/icient condition in~(\ref{k:pp}) under which we will mainly work.
See Conjecture~\ref{t:conj:aug10_3}.
Our most general claim is an elaborate one in
Conjecture~\ref{t:conj:aug11_1}.
In Section~\ref{k:sec:rc},
we explain the rigged conf\/igurations and
Kerov--Kirillov--Reshetikhin bijection
which are essential tools from
the combinatorial Bethe ansatz~\cite{KKR,KR} at $q=1$.
Based on Theorem~\ref{k:th:rc}, we identify the
solitons and strings in~(\ref{k:sost}).

In Section \ref{k:sec:ist},  we present the inverse scattering
formalism,  the solution algorithm of the initial value problem
together with applications to the volume formula of the
phase space and the dynamical period.
The content is based on a couple of conjectures whose status is
summarized in Section \ref{sec:suma}.
Our strategy is a synthesis of the combinatorial Bethe ans\"atze
at $q=1$ \cite{KKR,KR} and $q=0$ \cite{KN}, or in other words,
a modif\/ication of the KKR theory to match the periodic boundary condition.
The action-angle variables are constructed from
the rigged conf\/igurations invented at $q=1$
by a quasi-periodic extension (\ref{k:qpr}) followed by an
identif\/ication compatible with the $q=0$ Bethe equation (\ref{eq:sce}).
More concretely, the action variable is the $n$-tuple of Young diagrams
$\mu=(\mu^{(1)},\ldots, \mu^{(n)})$ preserved under time evolutions.
The angle variables live in a tropical analogue of the
Jacobi-variety ${\mathcal J}(\mu)$ and undergo a straight motion
with the velocity corresponding to a given time evolution.
Roughly speaking, the action and angle variables
represent the amplitudes and the phases of solitons, respectively
as demonstrated in Example \ref{k:ex:int}.
The modif\/ied KKR bijection $\Phi$, $\Phi^{-1}$
yield the direct and inverse scattering maps.
Schematically these aspects are summarized in the commutative
diagram (Conjecture~\ref{t:conj:aug21_5}):
\begin{equation*}
\begin{CD}
{\mathcal P}(\mu) @>{\Phi}>> {\mathcal J}(\mu) \\
@V{\mathcal T}VV @VV{\mathcal T}V\\
{\mathcal P}(\mu) @>{\Phi}>> {\mathcal J}(\mu)
\end{CD}
\end{equation*}
where ${\mathcal T}$ stands for the commuting family of time evolutions
$\big\{T^{(r)}_l\big\}$.
Since its action on ${\mathcal J}(\mu)$ is linear,
this diagram achieves the solution of the initial value problem conceptually.
Practical calculations can be found in
Examples \ref{k:ex:aa2} and \ref{k:ex:ivp}.

In Section \ref{k:sec:dcc} we decompose the level set ${\mathcal P}(\mu)$
further into connected components,
i.e., ${\mathcal T}$-orbits, and identify each of them with
${\mathbb Z}^g / F_{\boldsymbol{\gamma}}\, {\mathbb Z}^g$,
the set of integer points on a torus
with $g$ and~$F_{\boldsymbol{\gamma}}$ explicitly specif\/ied
in (\ref{k:gdef}) and (\ref{k:ftdef}).
Here $\boldsymbol{\gamma}$ denotes the order of symmetry
in the angle variables.
As a result we obtain (Theorem~\ref{k:th:bmain})
\begin{gather*}
|{\mathcal P}(\mu)|
= \sum_{\boldsymbol{\gamma}}\!\!
\underbrace{\det F_{\boldsymbol{\gamma}}}_{
\text{size of a ${\mathcal T}$-orbit}}
\underbrace{
\prod_{(ai) \in \overline{H}}\!\!
\frac{|\Lambda_{\gamma^{(a)}_i}(m_i^{(a)},p_i^{(a)})|}
{m_i^{(a)}/\gamma_i^{(a)}}}_{\text{number of ${\mathcal T}$-orbits}}
= (\det F)\prod_{(a i) \in \overline{H}} \frac{1}{m_i^{(a)}}
\binom{p_i^{(a)} + m_i^{(a)} - 1}{m_i^{(a)} - 1}.
\end{gather*}
The f\/irst equality follows from the decomposition into tori
and the second one is a slight calculation.
The last expression was known as the number of
Bethe roots at $q=0$ \cite{KN}.
Thus this identity of\/fers the Bethe ansatz formula a most elegant
interpretation by the structure of the phase space of the periodic SCA.

Once the linearization scheme is formulated,
it is straightforward to determine the dynamical period,
the smallest positive integer ${\mathcal N}$ satisfying
$T^{\mathcal N}(p) =p$ for any time evolution $T \in {\mathcal T}$
and path $p \in {\mathcal P}(\mu)$.
The result is given in Theorem \ref{k:th:dp} and Remark \ref{k:re:com}.
We emphasize that (\ref{t:eq:aug18_3}) is a
closed formula that gives the {\em exact} (not a multiple of)
dynamical period even when there are more than one solitons
with equal amplitudes and their order of symmetry
$\boldsymbol{\gamma}$ is nontrivial.
See Example \ref{k:ex:dp} in this paper for $n=2$
and also \cite[Example 4.10]{KTT} for $n=1$.
In Section \ref{k:sec:rba},
we explain the precise relation of our linearization scheme to the
Bethe ansatz at $q=0$ \cite{KN}.
The angle variables are actually in one to one correspondence
with what we call the Bethe root at $q=0$.
These results are natural generalizations of the $n=1$ case
proved in \cite{KTT,T}.
In Section \ref{k:sec:gen}, we treat the general case (\ref{k:pg})
relaxing the condition (\ref{k:pp}).
It turns out that the linearization scheme
remains the same provided one discards some time evolutions
and restricts the dynamics to a subgroup
${\mathcal T}'$ of ${\mathcal T}$.
This uncovers a new feature at $n>1$.

In Section \ref{k:sec:trt}, we derive an explicit formula
for the path $p \in {\mathcal P}(\mu)$ that
corresponds to a given action-angle variable
(Theorem \ref{k:con:ivp}).
This is a tropical analogue of the Jacobi inversion problem,
and the result is indeed expressed
by a tropical analogue of the Riemann theta function
(\ref{k:trt}):
\begin{equation*}
\Theta({\bf z}) = -\min_{{\bf n} \in {\mathbb Z}^G}
\left\{\tfrac{1}{2}\, {}^t{\bf n}B{\bf n}+ {}^t{\bf z}{\bf n}\right\}
\end{equation*}
having the quasi-periodicity
$\Theta({\bf z}+{\bf v}) = \Theta({\bf z})
+ {}^t{\bf v}B^{-1}\bigl({\bf z}+\frac{\bf v}{2}\bigr)$
for ${\bf v} \in B{\mathbb Z}^G$.
Here the $G \times G$ period matrix $B$ is specif\/ied by (\ref{k:B})
from the $n$-tuple of Young diagrams $\mu$
as in (\ref{k:rca0})--(\ref{k:gdef}).
Theorems  \ref{k:th:xt} and \ref{k:con:ivp} are derived from
the explicit formula of the KKR map
by the tropical tau function
$\tau_{\text{trop}}$ for $A^{(1)}_n $ \cite{KSY}.
It is a piecewise linear function on a rigged conf\/iguration
related to its charge.
The key to the derivation is the identity
\begin{equation}\label{k:star}
\Theta =
\lim_{\text{RC}\rightarrow \infty}\left(
\tau_{\text{trop}}(\text{RC})
-\text{divergent part}\right)
\end{equation}
f\/irst discovered in \cite{KS1} for $n=1$.
Here the limit sends the rigged conf\/iguration $\text{RC}$
into the inf\/initely large rigged conf\/iguration obtained by
the quasi-periodic extension.
See (\ref{k:taut}) for the precise form.
It is known \cite{KSY} that  $\tau_{\text{trop}}$ is indeed the
tropical analogue (ultradiscrete limit)
of a tau function in the KP hierarchy \cite{JM}.
Thus our result can be viewed as a
{\em fermionization} of the Bethe ansatz
and quasi-periodic solitons at a combinatorial level.
As applications, joint eigenvectors of $\big\{T^{(r)}_l\big\}$ are constructed
that possess every aspect as the Bethe eigenvectors at $q=0$
(Section \ref{k:sec:bv}),
the dynamical period is linked with the
$q=0$ Bethe eigenvalue (Section~\ref{k:sec:bedp}),
and miscellaneous calculations of some time average
are presented (Section~\ref{k:sec:mis}).

The main text is followed by 4 appendices.
Appendix~\ref{k:app:ins} recalls the row and column
insertions following~\cite{F} which is necessary to
understand the rule (\ref{k:bccb}) that governs
the local dynamics.
Appendix~\ref{k:app:fac} contains a proof of Theorem~\ref{k:th:fac}
based on crystal theory.
Appendix~\ref{k:app:kkr} is a quickest and
self-contained exposition of the algorithm
for the KKR bijection.
Appendix~\ref{k:app:proof1} is a proof of Theorem~\ref{k:th:rc}
using a result by Sakamoto~\cite{Sa} on energy of paths.

\section[Periodic $A^{(1)}_n$ soliton cellular automaton]{Periodic $\boldsymbol{A^{(1)}_n}$ soliton cellular automaton}
\label{k:sec:model}

\subsection[Crystal $B^{r,l}$ and combinatorial $R$]{Crystal $\boldsymbol{B^{r,l}}$ and combinatorial $\boldsymbol{R}$}
\label{k:sec:rect}

Let $n$ be a positive integer.
For any pair of integers $r$, $l$ satisfying $1 \leq r \leq n$ and $l \geq 1$,
we denote by $B^{r,l}$ the set of all $r \times l$
rectangular semistandard tableaux with entries chosen from the set
$\{1,2,\ldots,n+1\}$.
For example when $n=2$, one has
$B^{1,1} = \{ 1,2,3 \}$, $B^{1,2} = \{ 11,12,13,22,23,33 \}$
and
$B^{2,2} = \left\{ {1 1 \atop 2 2}, {1 1 \atop 2 3},
{1 1 \atop 3 3},{1 2 \atop 2 3},{1 2 \atop 3 3},{2 2 \atop 3 3} \right\}$.
The set $B^{r,l}$ is equipped with the structure of
crystal \cite{KMN2,Sh} for the Kirillov--Reshetikhin module of
the quantum af\/f\/ine algebra $U_q\big(A^{(1)}_n\big)$.

There is a canonical bijection, the isomorphism of crystals, between
$B^{r,l} \otimes B^{1,1}$ and $B^{1,1} \otimes B^{r,l}$
called {\em combinatorial R}.
It is denoted either by $R(b\otimes c) =  c' \otimes b' $
or $R(c' \otimes b' )=b \otimes c$,
or simply $b\otimes c \simeq c' \otimes b'$,
where  $b \otimes c \in B^{r,l} \otimes B^{1,1}$ and
$c' \otimes b' \in B^{1,1} \otimes B^{r,l}$.
$R$ is uniquely determined by the condition that the product
tableaux $c\cdot b$ and $b' \cdot c'$ coincide \cite{Sh}.
Here, the product $c\cdot b$ for example signif\/ies
the column insertion of~$c$ into~$b$, which is also obtained
by the row insertion of~$b$ into~$c$~\cite{F}.
Since~$c$ and~$c'$ are single numbers in our case,
it is simplest to demand the equality
\begin{equation}\label{k:bccb}
(c \rightarrow b) = (b' \leftarrow c')
\end{equation}
to f\/ind the image $b\otimes c \simeq c' \otimes b'$.
See Appendix \ref{k:app:ins} for
the def\/initions of the row and column insertions.
The insertion procedure also determines an integer
$H(b\otimes c) = H(c'\otimes b')$ called {\em local energy}.
We specify it as $H=0$ or $1$
according as the shape of the common product tableau is
$((l+1), l^{r-1})$ or $(l^r,1)$, respectively\footnote{Up to a constant shift, this is
the original def\/inition of local energy \cite{KMN} times $(-1)$.}.
Note that $R$ and $H$ refer to the pair $(r,l)$
although we suppress the dependence on it in the notation.

\begin{example}
Let $r=2$, $l=3$.
\begin{gather*}
R \left(
\begin{array}{l}
1 1 2 \\
2 2 3 \\
\end{array}
\otimes 1
\right) = 2 \otimes
\begin{array}{l}
1 1 1 \\
2 2 3 \\
\end{array}
, \qquad
\mbox{product tableau}=
\begin{array}{l}
1 1 1 2 \\
2 2 3 \\
\end{array}
, \qquad
H=0,
\\
R \left(
\begin{array}{l}
1 1 2 \\
2 2 3 \\
\end{array}
\otimes 3
\right) = 1 \otimes
\begin{array}{l}
1 2 2 \\
2 3 3 \\
\end{array}
,\qquad
\mbox{product tableau}=
\begin{array}{l}
1 1 2 \\
2 2 3 \\
3
\end{array}
, \qquad
H=1.
\end{gather*}
\end{example}

We depict the relation
$R(b \otimes c)= c' \otimes b'$ (\ref{k:bccb})
and $H(b\otimes c) = e$ as

\begin{picture}(50,60)(-170,-17)
\put(-10,7){$b$}\put(12,30){$c$}
\put(15,-5){\line(0,1){30}}\put(-2,10){\line(1,0){34}}
\put(34,7){$b'$}
\put(12,-16){$c'$}

\put(7,13){$e$}

\end{picture}

We will often suppress $e$.
The horizontal and vertical lines here carry
an element from $B^{r,l}$ and $B^{1,1}$, respectively.
We remark that $R$ is trivial,  namely
$R(b \otimes c ) =b \otimes c$,
for $(r,l)=(1,1)$ by the def\/inition.

\subsection{Def\/inition of dynamics}\label{k:sec:DD}

Fix a positive integer $L$ and set $B = (B^{1,1})^{\otimes L}$.
An element of $B$ is called a {\em path}.
A path $b_1\otimes \cdots \otimes b_L$ will often
be written simply as a word $b_1b_2\ldots b_L$.
Our periodic $A^{(1)}_n$ soliton cellular automaton (SCA)
is a dynamical system on a subset of $B$.
To def\/ine the time evolution $T^{(r)}_l$ associated with $B^{r,l}$, we consider
a bijective map $B^{r,l} \otimes B \rightarrow B \otimes B^{r,l}$ and the
local energy obtained by repeated use of $R$.
Schematically, the map
\begin{equation*}
\begin{split}
B^{r,l} \otimes B \qquad &\rightarrow \qquad B \otimes B^{r,l}\\
v\otimes (b_1\otimes \cdots \otimes b_L)
&\mapsto
(b'_1\otimes \cdots \otimes b'_L) \otimes v'
\end{split}
\end{equation*}
and the local energy $e_1, \ldots, e_L$ are def\/ined
by the composition of the previous diagram as follows.

\begin{equation}\label{k:hone}
\begin{picture}(150,40)(-140,-13)
\put(-10,7){$v$}
\put(11,30){$b_1$}\put(15,-5){\line(0,1){30}}\put(12,-16){$b'_1$}
\put(4,14){$e_1$}

\put(36,30){$b_2$}\put(40,-5){\line(0,1){30}}\put(37,-16){$b'_2$}
\put(29,14){$e_2$}

\put(-2,10){\line(1,0){57}}

\put(67,6){$\cdots$}

\put(90,10){\line(1,0){30}}
\put(107,-5){\line(0,1){32}}
\put(103,30){$b_L$}\put(104,-16){$b'_L$}
\put(94,14){$e_L$}\put(123,7){$v'$}
\end{picture}
\end{equation}

The element $v$ or $v'$ is called a {\em carrier}.
Set $p'=b'_1 \otimes \cdots \otimes b'_L$.
Given a path $p=b_1\otimes \cdots \otimes b_L \in B$,
we regard $v'=v'(v; p), e_k = e_k(v; p)$ and $p'=p'(v; p)$
as the functions of $v$ containing $p$ as a~`parameter'.
Naively, we wish to def\/ine the time evolution $T^{(r)}_l$ of the path $p$
as $T^{(r)}_l(p) = p'$ by using a carrier $v$ that satisf\/ies the
periodic boundary condition $v=v'(v; p)$.
This idea indeed works without a dif\/f\/iculty for $n=1$ \cite{KTT}\footnote{Even for $n=1$, it becomes nontrivial for higher spin case \cite{KS3}.}.
In addition, $T^{(1)}_1$ is always well def\/ined for general $n$,
yielding the cyclic shift:
\begin{equation}\label{k:shf}
T^{(1)}_1(b_1\otimes \cdots \otimes b_{L-1}\otimes b_L )
= b_L\otimes b_1 \otimes \cdots \otimes b_{L-1}.
\end{equation}
This is due to the triviality of $R$ on $B^{1,1} \otimes B^{1,1}$
mentioned above.
In fact the unique carrier is specif\/ied as $v=v'=b_L$.
Apart from this however,
one encounters the three problems $(i)$--$(iii)$
to overcome in general for $n\ge 2$.
For some $p\in B$, one may suf\/fer from
\begin{enumerate}\itemsep=0pt
\item[$(i)$] Non-existence.  There may be no carrier $v$ satisfying
$v=v'(v;p)$.

\item[$(ii)$]
Non-uniqueness. There may be more than one carriers,
say  $v_1, v_2, \ldots, v_m$
such that $v_j = v'(v_j; p)$ but $p'(v_i; p) \neq p'(v_j; p)$ or
$e_k(v_i; p) \neq e_k(v_j; p)$ for some $i$, $j$, $k$.

\end{enumerate}
To cope with these problems, we introduce the notion of evolvability of a path
according to \cite{KT2}.
A path $p \in B$ is said $T_l^{(r)}$-{\em evolvable} if
there exist $v \in B^{r,l}$ such that $v=v'(v;p)$, and moreover
$p'=p'(v;p), e_k=e_k(v;p)$ are unique for possibly non-unique
choice of $v$.\footnote{We expect that the uniqueness of the local energy $e_k(v;p)$
follows from the uniqueness of $p'(v;p)$ alone.}
In this case we def\/ine
\begin{gather}\label{k:te}
T_l^{(r)} (p)=p', \qquad  E^{(r)}_l(p) = e_1+\cdots + e_L,
\end{gather}
indicating that the time evolution operator $T_l^{(r)}$ is acting on $p$.
The quantity $E^{(r)}_l(p)\in {\mathbb Z}_{\ge 0}$
is called an {\em energy} of the path $p$.

By the standard argument using the Yang--Baxter equation
of the combinatorial $R$, one can show
\begin{proposition}[\cite{KT2}]\label{k:pr:tte}
The commutativity
$T^{(a)}_iT^{(b)}_j(p) = T^{(b)}_jT^{(a)}_i(p)$ is valid if all the
time evolutions $T^{(r)}_l$ here
are acting on $T^{(r)}_l$-evolvable paths.
Moreover,
$E^{(a)}_i\big(T^{(b)}_j(p)\big)=E^{(a)}_i(p)$ and
$E^{(b)}_j\big(T^{(a)}_i(p)\big)=E^{(b)}_j(p)$ hold.
\end{proposition}

If $p$ is $T_l^{(r)}$-evolvable for all $r$, $l$, then it is simply
called {\em evolvable}.
The third problem in def\/ining the dynamics is
\begin{enumerate}\itemsep=0pt
\item[$(iii)$] Even if $p$ is evolvable, $T^{(r)}_l(p)$ is not necessarily so in general.
\end{enumerate}

For instance,
$p = 112233 \in \big(B^{1,1}\big)^{\otimes 6}$ is evolvable
but $T^{(2)}_1(p) = 213213$ is not
$T^{(2)}_1$-evolvable.
In fact, the non-uniqueness problem $(ii)$ takes place as follows:
\begin{equation*}
\begin{picture}(400,40)(0,-20)

\multiput(0,0)(32,0){6}{
\put(-8,0){\line(1,0){16}}\put(0,-10){\line(0,1){20}}
}

\put(-3,14){2}\put(-18,2){1}\put(-18,-9){3}
\put(29,14){1}\put(14,2){1}\put(14,-9){2}
\put(61,14){3}\put(46,2){1}\put(46,-9){2}
\put(93,14){2}\put(78,2){2}\put(78,-9){3}
\put(125,14){1}\put(110,2){2}\put(110,-9){3}
\put(157,14){3}\put(142,2){1}\put(142,-9){3}
                          \put(174,2){1}\put(174,-9){3}

\put(-3,-20){3}
\put(29,-20){1}
\put(61,-20){1}
\put(93,-20){2}
\put(125,-20){2}
\put(157,-20){3}

\multiput(225,0)(0,0){1}{
\multiput(0,0)(32,0){6}{
\put(-8,0){\line(1,0){16}}\put(0,-10){\line(0,1){20}}
}

\put(-3,14){2}\put(-18,2){2}\put(-18,-9){3}
\put(29,14){1}\put(14,2){2}\put(14,-9){3}
\put(61,14){3}\put(46,2){1}\put(46,-9){3}
\put(93,14){2}\put(78,2){1}\put(78,-9){3}
\put(125,14){1}\put(110,2){1}\put(110,-9){2}
\put(157,14){3}\put(142,2){1}\put(142,-9){2}
                        \put(174,2){2}\put(174,-9){3}

\put(-3,-20){2}
\put(29,-20){2}
\put(61,-20){3}
\put(93,-20){3}
\put(125,-20){1}
\put(157,-20){1}
}

\end{picture}
\end{equation*}
To discuss the issue $(iii)$,
we need to prepare several def\/initions describing the
spectrum $\big\{E^{(r)}_l\big\}$.
Let $\mu=(\mu^{(1)}, \ldots, \mu^{(n)})$ be
an $n$-tuple of Young diagrams.
From $\mu^{(a)}$,
we specify the data $m^{(a)}_i$, $l^{(a)}_i$, $g_a$  as in the
left diagram in (\ref{k:rca0}).
\begin{equation}\label{k:rca0}
\includegraphics{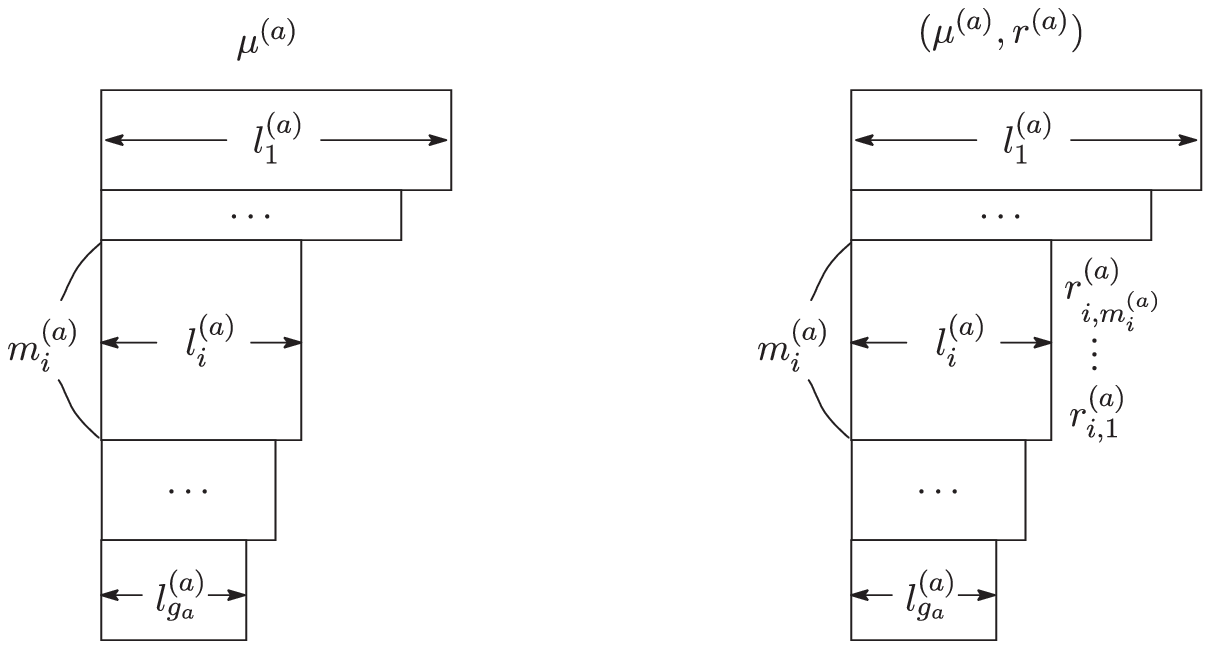}
\end{equation}
Here, $m^{(a)}_i\times l^{(a)}_i$ rectangle part
is called the $(a,i)$ block.
Furthermore we introduce\footnote{The meaning of $i$
in $m^{(a)}_i$ and  $p^{(a)}_i$ here has been altered(!) from
the literatures \cite{KN,KS1,KS2,KS3,KSY,KTT,KT1,KT2}.
When the space is tight, we will often write $(ai\alpha)$
for $(a,i,\alpha)\in H$ and $(ai)$ for
$(a,i) \in \overline{H}$, etc. }
\begin{gather}
\begin{split}
&H  = \big\{(a, i, \alpha)\mid
1 \le a \le n, \,1 \le i \le g_a, \,
1 \le \alpha \le m^{(a)}_i\big\},
\\
&\overline{H} = \{(a, i) \mid
1 \le a \le n, \,1 \le i \le g_a\},
\end{split}\label{k:hdef}\\
p^{(a)}_i  =
 L\delta_{a 1}-\sum_{(b j \beta) \in H}
C_{ab}\min\big(l^{(a)}_i, l^{(b)}_j\big)
=L\delta_{a 1}-\sum_{(bj) \in \overline{H}}
C_{ab}\min\big(l^{(a)}_i, l^{(b)}_j\big)m^{(b)}_j, \label{k:pdef}\\
F   = (F_{ai, bj})_{(ai),(bj) \in \overline{H}},\qquad
F_{ai,bj} = \delta_{ab}\delta_{ij}p^{(a)}_i
+ C_{ab}\min \big(l^{(a)}_i, l^{(b)}_j\big)m^{(b)}_j, \label{k:f}
\\
G =|H| =\sum_{a=1}^n\sum_{i=1}^{g_a}m^{(a)}_i,\qquad
g = |\overline{H}| = \sum_{a=1}^ng_a,\label{k:gdef}
\end{gather}
where $(C_{ab})_{1 \le a,b \le n}$ is the Cartan
matrix of $A_n$, i.e.,
$C_{ab}=2\delta_{ab}-\delta_{a,b+1}-\delta_{a,b-1}$.
The quantity $p^{(a)}_i$ is called a {\em vacancy number}.

Given $\mu=(\mu^{(1)}, \ldots, \mu^{(n)})$, we def\/ine
the sets of paths
$B \supset {\mathcal P} \supset {\mathcal P}(\mu)
\supset {\mathcal P}_+(\mu)$ by
\begin{gather}
{\mathcal P}  = \{p \in B\mid
\#(1) \geq \#(2) \geq \cdots \geq \#(n+1)\},\label{k:pop}\\
{\mathcal P}(\mu)  = \Bigg\{p \in {\mathcal P}\mid
p: \text{evolvable},\;
E^{(a)}_l(p) = \sum_{i=1}^{g_a}\min\big(l,l^{(a)}_i\big)m^{(a)}_i\Bigg\},
\label{k:pmu}\\
{\mathcal P}_+(\mu)  =
\{p \in {\mathcal P}(\mu)\mid
p: \text{highest} \},
\label{k:hip}
\end{gather}
where $\#(a)$ is the number of $a \in B^{1,1}$
in $p = b_1 \otimes \cdots \otimes b_L$.
In view of the Weyl group symmetry \cite[Theorem 2.2]{KT2},
and the obvious property that $T^{(r)}_l$ is weight preserving,
we restrict ourselves to ${\mathcal P}$
which is the set of paths with nonnegative weight.
A path $b_1 \otimes \cdots \otimes b_L$ is {\em highest} if
the pref\/ix $b_1\otimes \cdots \otimes b_j$
satisf\/ies the condition
$\#(1) \geq \#(2) \geq \cdots \geq \#(n+1)$ for all $1 \le j \le L$.
We call ${\mathcal P}(\mu)$ the {\em level set}
associated with $\mu$.
The $n$-tuple $\mu$ of Young diagrams or equivalently the data
$(m^{(a)}_i, l^{(a)}_i)_{(ai) \in \overline{H}}$ will be referred to
as the {\em soliton content} of ${\mathcal P}(\mu)$
or the paths contained in~it.
This terminology comes from the fact that
when $L$ is large and $\#(1)\gg \#(2),\ldots, \#(n+1)$,
it is known that the paths in ${\mathcal P}(\mu)$ consist of
$m^{(1)}_i$ solitons with amplitude $l^{(1)}_i$
separated from each other or in the course of collisions.
Here a soliton with amplitude $l$ is a part of a path of the form
$j_1 \otimes j_2\otimes \cdots \otimes j_l \in (B^{1,1})^{\otimes l}$
with $j_1\ge \cdots \ge j_l\ge 2$
by regarding $1 \in B^{1,1}$ as an empty background.
Note that a soliton is endowed not only with the amplitude $l$ but also
the internal degrees of freedom $\{j_i\}$.
The data $m^{(a)}_i$, $l^{(a)}_i$ with $a>1$ are relevant to
the internal degrees of freedom of solitons.

Note that the relation in (\ref{k:pmu}) can be inverted as
($E^{(a)}_0 :=0$)
\begin{equation}\label{k:eee}
-E^{(a)}_{l-1}+2E^{(a)}_l-E^{(a)}_{l+1}
=\begin{cases}
m^{(a)}_i & \text{ if } l= l^{(a)}_i \text{ for some }
1 \le i \le g_a,\\
0 & \text{ otherwise}.
\end{cases}
\end{equation}
Therefore the correspondence between the soliton content
$\mu=\big(\mu^{(1)}, \ldots, \mu^{(n)}\big)$ and
the energy $\{E^{(a)}_l\}$ is one to one.
Either $\{E^{(a)}_l\}$ or $\mu$
are conserved quantities in the following sense;
if $p \in {\mathcal P}(\mu)$ and
$T^{(r)}_l(p)$ is evolvable, then
Proposition \ref{k:pr:tte} tells that
$T^{(r)}_l(p) \in {\mathcal P}(\mu)$.
Comparing~(\ref{k:pdef}) and~(\ref{k:pmu}),
one can express the vacancy number also as
\begin{equation*}
p^{(a)}_i = L\delta_{a1}-\sum_{b=1}^n C_{ab}E^{(b)}_{l^{(a)}_i}.
\end{equation*}

Now we are ready to propose a suf\/f\/icient condition under which
the annoying feature in problem $(iii)$ is absent:
\begin{equation}\label{k:pp}
\mu=\big(\mu^{(1)}, \ldots, \mu^{(n)}\big) \ \
\text{satisf\/ies}\ \  p^{(a)}_i \ge 1 \ \ \text{for all} \ \
(a,i) \in \overline{H}.
\end{equation}
This condition was f\/irst introduced in \cite{KS3}.

\begin{conjecture}\label{t:conj:aug10_3}
Under the condition \eqref{k:pp},
$T_l^{(r)}({\mathcal P}(\mu)) = {\mathcal P}(\mu)$ holds
for any $r$, $l$.
\end{conjecture}

For a most general claim, see Conjecture~\ref{t:conj:aug11_1}.

The equality $T_l^{(r)}({\mathcal P}(\mu)) = {\mathcal P}(\mu)$
implies that one can apply any time evolution
in any order for arbitrary times on
any path $p \in {\mathcal P}(\mu)$.
In particular the inverse $\big(T_l^{(r)}\big)^{-1}$ exists.
We denote by $\Sigma (p)$ the set of all paths generated from $p$
in such a manner.
By the def\/inition it is a subset of the level set ${\mathcal P}(\mu)$.
Call $\Sigma (p)$ the {\em connected component} of the level set
containing the path $p$.

One can give another characterization of $\Sigma (p)$ in terms of group actions.
Let ${\mathcal T}$ be the abelian group generated by all $T_l^{(r)}$'s.
Then ${\mathcal T}$ acts on the level set ${\mathcal P}(\mu)$.
This group action is not transitive in general,
i.e., ${\mathcal P}(\mu)$ is generically
decomposed into several ${\mathcal T}$-orbits.
In what follows, connected components and
${\mathcal T}$-orbits mean the same thing.

\begin{example}\label{t:ex:sept1_2}
Let $n=2, L=8$ and $\mu = ((211),(1))$.
The vacancy numbers are given as
$p^{(1)}_1 =1$, $p^{(1)}_2 =3$, $p^{(2)}_1 =1$, hence
the condition (\ref{k:pp}) is satisf\/ied.
\begin{equation*}
\begin{picture}(200,65)(-120,40)
\setlength{\unitlength}{0.5mm}

\put(15,68){$\mu^{(1)}$}
\put(65,68){$\mu^{(2)}$}
\multiput(0,0)(0,-10){2}{
\put(10,60){\line(1,0){20}}}

\multiput(0,0)(0,10){2}{
\put(10,30){\line(1,0){10}}}

\put(10,30){\line(0,1){30}}
\put(20,30){\line(0,1){30}}
\put(30,50){\line(0,1){10}}

\multiput(-35,0)(0,0){1}{
\put(100,60){\line(1,0){10}}
\put(100,50){\line(1,0){10}}
\put(100,50){\line(0,1){10}}
\put(110,50){\line(0,1){10}}}

\end{picture}
\end{equation*}
Then we have ${\mathcal P}(\mu)
= \Sigma (p_X) \sqcup \Sigma (p_Y)$ with
$p_X = 11221123$ and $p_Y=11221213$.
The connected components are written as $\Sigma (p_X) = \bigsqcup_{i=1}^9
X_i$ and $\Sigma (p_Y) = \bigsqcup_{i=1}^9 Y_i$
where $X_i$ and $Y_i$ are given in the following table.
\begin{center}
\begin{tabular}{c|cc}
\hline
$i$ & $X_i$ & $Y_i$ \\
\hline
1 & [11221123] & [11221213] \\
2 & [11122123] & [11121223] \\
3 & [21112123] & [21121123] \\
4 & [21121213] & [21211213] \\
5 & [21212113] & [12212113] \\
6 & [12122113] & [11212213] \\
7 & [12112213] & [11211223] \\
8 & [12111223] & [21211123] \\
9 & [12211123] & [12211213] \\
\hline
\end{tabular}
\end{center}
Here the symbol [   ] denotes the set of all cyclic shifts of the entry.
For instance, $[123] = \{ 123, 231, 312 \}$.
For any $i \in {\mathbb Z} / 9 {\mathbb Z}$
we have $T_1^{(1)} (X_i) = X_i$, $T_{\geq 2}^{(1)} (X_i) = X_{i+1}$, $T_{\geq 1}
^{(2)} (X_i) = X_{i+3}$ and the
same relations for $Y_i$.
Hence the claim of Conjecture \ref{t:conj:aug10_3} is valid for this
${\mathcal P}(\mu)$.
\end{example}

The time evolution $T^{(1)}_\infty$ has an especially simple description,
which we shall now explain.
Let~$B_1$ be the set of paths in which
$1 \in B^{1,1}$ is contained most. Namely,
\begin{equation}\label{k:B1}
B_1 =\{p \in B\mid
\sharp(1)\ge \sharp(a), \;2\le a\le n+1\}.
\end{equation}
Thus we see ${\mathcal P} \subset B_1 \subset B$.
We def\/ine the weight preserving map $K_a: B_1 \rightarrow B_1$
for $a=2,3,\ldots, n+1$.
Given a path $p \in B_1$ regarded as a word
$p=b_1b_2\ldots b_L \in \{1,\ldots, n+1\}^L$,
the image $K_a(p)$ is determined  by the following procedure.
\begin{enumerate}\itemsep=0pt

\item[$(i)$] Ignore all the numbers except $1$ and $a$.

\item[$(ii)$] Connect every adjacent pair $a1$ (not $1a$) by an arc, where
$1$ is on the right of $a$ {\em cyclically}.

\item[$(iii)$] Repeat $(ii)$ ignoring the already connected pairs until
all $a$'s are connected to some $1$.

\item[$(iv)$] Exchange $a$ and $1$ within each connected pair.

\end{enumerate}

\begin{theorem}\label{k:th:fac}
Suppose $p \in B_1$ and $p$ is $T^{(1)}_\infty$-evolvable.
Then, $T^{(1)}_\infty$ is factorized as follows:
\begin{equation}\label{k:tkk}
T^{(1)}_\infty(p) = K_2K_3\cdots K_{n+1}(p).
\end{equation}
\end{theorem}

A proof is available in Appendix \ref{k:app:fac}, where
each $K_a$ is obtained as
a gauge transformed simple ref\/lection for type $A^{(1)}_n$
af\/f\/ine Weyl group.
The dynamics of the form (\ref{k:tkk})
has an interpretation as the ultradiscrete
Toda or KP f\/low \cite{HHIKTT,Yy}.
See also \cite[Theorem VI.2]{MIT1}.
As an additional remark,
a factorized time evolution similar to (\ref{k:tkk}) or
equivalently (\ref{k:tss}) has been formulated
for a periodic SCA associated to
any non-exceptional af\/f\/ine Lie algebras \cite{KT1}.
In the inf\/inite (non-periodic) lattice case,
it was f\/irst invented by Takahashi \cite{Th} for type $A^{(1)}_n$.
The origin of such a factorization is the
factorization of the combinatorial R itself in a certain asymptotic domain,
which has been proved uniformly for all non-exceptional types
in \cite{HKT}.

\begin{example}\label{k:ex:arc}
In the left case of Example \ref{k:ex:int},
the time evolution $T^{(1)}_3$ is actually equal to $T^{(1)}_\infty$.
The time evolution from  $t=2$ to $t=3$ is reproduced
by $T^{(1)}_3=T^{(1)}_\infty = K_2K_3$ as follows:
\begin{equation*}
\begin{picture}(250,100)

\put(190,70){\vector(0,-1){20}}\put(195,58){$K_3$}

\put(190,37){\vector(0,-1){20}}\put(195,25){$K_2$}

\multiput(0,70)(0,0){1}{

\put(47,10){\line(0,1){8}}\put(47,18){\line(1,0){16.4}}
\put(52.5,10){\line(0,1){5}}\put(52.5,15){\line(1,0){5.78}}
\put(58,10){\line(0,1){5}}
\put(63.4,10){\line(0,1){8}}
\put(68.9,10){\line(0,1){5}}\put(68.9,15){\line(1,0){10.9}}
\put(79.9,10){\line(0,1){5}}
\put(85.3,10){\line(0,1){11}}\put(85.3,21){\line(-1,0){49}}
\put(96.3,10){\line(0,1){5}}\put(96.3,15){\line(1,0){10.9}}
\put(107.3,10){\line(0,1){5}}
\put(173,10){\line(0,1){11}}\put(173,21){\line(1,0){5}}

\put(0,0){$t=2:\quad 331132111321111221122213$}
}

\multiput(0,40)(0,0){1}{

\put(47,10){\line(0,1){11}}\put(47,21){\line(-1,0){5}}
\put(74.4,10){\line(0,1){5}}\put(74.4,15){\line(1,0){16.4}}
\put(90.8,10){\line(0,1){5}}
\put(101.8,10){\line(0,1){5}}\put(101.8,15){\line(1,0){10.96}}
\put(112.7,10){\line(0,1){5}}
\put(129.2,10){\line(0,1){8}}\put(129.2,18){\line(1,0){16.4}}
\put(134.6,10){\line(0,1){5}}\put(134.6,15){\line(1,0){5.48}}
\put(140.1,10){\line(0,1){5}}
\put(145.6,10){\line(0,1){8}}
\put(151.1,10){\line(0,1){11}}\put(151.1,21){\line(1,0){25}}
\put(156.5,10){\line(0,1){8}}\put(156.5,18){\line(1,0){16.4}}
\put(162.0,10){\line(0,1){5}}\put(162.0,15){\line(1,0){5.48}}
\put(167.5,10){\line(0,1){5}}
\put(173,10){\line(0,1){8}}

\put(0,0){\phantom{$t=2:\quad$ }$113312331123111221122211$}
}

\multiput(0,10)(0,0){1}{
\put(0,0){$t=3:\quad 213311332113211112211122$}
}

\end{picture}
\end{equation*}
\end{example}

The procedure $(i)$--$(iv)$ is
delightfully simple but inevitably non-local as the
result of expelling the carrier out from the description.
It is an interesting question whether the other typical time evolutions
$T^{(2)}_\infty, \ldots, T^{(n)}_\infty$ admit a similar description
without a carrier.

\subsection{Rigged conf\/iguration}\label{k:sec:rc}

An $n$-tuple of Young diagrams
$\mu=(\mu^{(1)}, \ldots, \mu^{(n)})$ satisfying
$p^{(a)}_i \ge 0$ for all $(a,i) \in \overline{H}$ is
called a {\em configuration}.
Thus, those $\mu$ satisfying (\ref{k:pp}) form a subset of conf\/igurations.
Consider the conf\/iguration
$\mu$ attached with the integer arrays called {\em rigging}
${\bf r}= (r^{(a)})_{1 \le a \le n}
= (r^{(a)}_{i,\alpha})_{(ai\alpha)\in H}$
as in the right diagram in (\ref{k:rca0}).
The combined data
$(\mu, {\bf r})
= ((\mu^{(1)}, r^{(1)}), \ldots, (\mu^{(n)}, r^{(n)}))$
is called a {\em rigged configuration} if the condition
\begin{equation}\label{k:rcon}
0 \le r^{(a)}_{i,1}\le  r^{(a)}_{i,2}\le \cdots
\le r^{(a)}_{i, m^{(a)}_i} \le p^{(a)}_i\quad\text{for all }\ \
(a,i) \in \overline{H}
\end{equation}
is satisf\/ied\footnote{We do not include $L$ in the def\/inition of rigged conf\/iguration
understanding that it is f\/ixed.
It is actually necessary to determine the vacancy number
$p^{(a)}_i$ (\ref{k:pdef}) appearing as
the upper bound of the rigging (\ref{k:rcon}).}.
We let $\text{RC}(\mu)$ denote the set of
rigged conf\/igurations whose conf\/iguration is $\mu$.
It is well known \cite{KKR, KR} that the highest paths in $B$ are in one to one
correspondence with the rigged conf\/igurations
by the Kerov--Kirillov--Reshetikhin (KKR) map
\begin{equation}\label{k:kkr}
\phi : \{p \in B \mid
p: \text{highest}\}
\ \rightarrow\
\sqcup_\mu \text{RC}(\mu),
\end{equation}
where the union is taken over all the conf\/igurations.
The both $\phi$ and $\phi^{-1}$ can be described by an
explicit algorithm as described in Appendix~\ref{k:app:kkr}.
Our main claim in this subsection is the following,
which is a ref\/inement of~(\ref{k:kkr})
with respect to $\mu$ and an adaptation to the periodic setting.

\begin{theorem}\label{k:th:rc}
The restriction of the KKR map $\phi$ to
${\mathcal P}_+(\mu)$ separates the image
according to~$\mu$:
\begin{equation*}
\phi: \ {\mathcal P}_+(\mu) \rightarrow \text{\rm RC}(\mu).
\end{equation*}
\end{theorem}

The proof is available in Appendix~\ref{k:app:proof1}.
We expect that this restricted injection is still a~bijection
but we do not need this fact in this paper.

Let us explain some background and signif\/icance of Theorem \ref{k:th:rc}.
Rigged conf\/igurations were invented \cite{KKR} as combinatorial
substitutes of solutions to the Bethe equation under
{\em string hypothesis} \cite{Be}.
The KKR map $\phi^{-1}$ is the combinatorial analogue of
producing the Bethe vector from the solutions to the Bethe equation.
The relevant integrable system is $sl_{n+1}$ Heisenberg chain,
or more generally the rational vertex model
associated with $U_q\big(A^{(1)}_n\big)$ at $q=1$.
In this context, the conf\/iguration $\mu$ (\ref{k:rca0})
is the {\em string content} specifying that there are $m^{(a)}_i$ strings
with color $a$ and length $l^{(a)}_i$.
Theorem \ref{k:th:rc} identif\/ies the two meanings of $\mu$.
Namely,  the soliton content for ${\mathcal P}(\mu)$
measured by energy is equal to the string content for $\text{RC}(\mu)$
determined by the KKR bijection.
Symbolically, we have the identity
\begin{equation}\label{k:sost}
\text{soliton} = \text{string},
\end{equation}
which lies at the heart of the whole
combinatorial Bethe ansatz approaches to the
soliton cellular automata \cite{KS1,KS2,KS3,KSY,KTT,KT1,KT2}.
It connects the energy and the conf\/iguration by
(\ref{k:eee}),  and in a broader sense, crystal theory and Bethe ansatz.
Further arguments will be given around~(\ref{k:ss}).
For $n=1$, Theorem~\ref{k:th:rc}
has been obtained in \cite[Proposition~3.4]{KTT}.
In the sequel, we shall call the $n$-tuple of Young diagrams
$\mu=(\mu^{(1)},\ldots, \mu^{(n)})$
either as soliton content, string content or conf\/iguration
when $p^{(a)}_i \ge 0$ for all $(a,i) \in \overline{H}$.

A consequence of Theorem \ref{k:th:rc} is that for any path
$b_1\otimes \cdots \otimes b_L \in {\mathcal P}(\mu)$
with soliton content $\mu$, the number $\#(a)$ of $a \in B^{1,1}$
in $b_1,\ldots, b_L$ is given by
\begin{equation}\label{k:an}
\#(a) = \big|\mu^{(a-1)}\big|-\big|\mu^{(a)}\big| \qquad (1 \le a \le n+1),
\end{equation}
where we set $|\mu^{(0)}| = L$ and $|\mu^{(n+1)}| = 0$.
This is due to a known property of the KKR map $\phi$
and the fact that $\mu$ is a conserved quantity and
Conjecture \ref{k:con:pg}.

\begin{example}\label{t:ex:sept2_1}
Let $n=2, L=8$ and consider the same
conf\/iguration $\mu = ((211),(1))$ as
Example \ref{t:ex:sept1_2}.
The riggings are to obey the conditions
$0 \leq r^{(1)}_{1,1}, r^{(2)}_{1,1} \leq 1$
and $0 \leq r^{(1)}_{2,1} \leq
r^{(1)}_{2,2} \leq 3$.
\begin{equation*}
\begin{picture}(200,65)(-120,40)
\setlength{\unitlength}{0.5mm}

\put(5,68){$(\mu^{(1)}, r^{(1)})$}
\put(83,68){$(\mu^{(2)}, r^{(2)})$}
\multiput(0,0)(0,-10){2}{
\put(10,60){\line(1,0){20}}}

\multiput(0,0)(0,10){2}{
\put(10,30){\line(1,0){10}}}

\put(10,30){\line(0,1){30}}
\put(20,30){\line(0,1){30}}
\put(30,50){\line(0,1){10}}

\put(-20,52){$p^{(1)}_1=1$}\put(-20,38){$p^{(1)}_2=3$}
\put(35,53){$r^{(1)}_{1,1}$}
\put(25,41){$r^{(1)}_{2,2}$}
\put(25,28){$r^{(1)}_{2,1}$}

\multiput(-30,0)(0,0){1}{
\put(95,52){$p^{(2)}_1=1$}\put(140,52){$r^{(2)}_{1,1}$}

\multiput(25,0)(0,0){1}{
\put(100,60){\line(1,0){10}}
\put(100,50){\line(1,0){10}}
\put(100,50){\line(0,1){10}}
\put(110,50){\line(0,1){10}}}}

\end{picture}
\end{equation*}
Hence there are $2 \cdot 2 \cdot 10 = 40$ rigged conf\/igurations in
$\text{RC}(\mu)$.
The elements of ${\mathcal P}_+(\mu)$ and the riggings for
their images under the KKR map $\phi$ are given in the following table.
\footnotesize
\begin{center}
\begin{tabular}{c|c|c||c|c}
\hline
$i$ & highest paths in $X_i$ & riggings & highest paths in $Y_i$ & riggings
\\
\hline
1 & 11221123, 11231122, 12311221 & 0331, 1110, 0000 & 11221213, 12131122 &
0321, 1100 \\
2 & 11122123, 11221231, 12311122 & 1331, 0221, 1000 & 11121223, 11212231 &
1321, 0211 \\
3 & 11121232, 11212321 & 1221, 0111 & 11211232, 12112321,  11232112 & 1211,
0101, 0310 \\
4 & 11212132, 12121321, 12132112 & 1111, 0001, 0300 & 12112132, 11213212,
12132121 & 1101, 1310, 0200 \\
5 & 12121132, 12113212 & 1001, 1300 & 12113122 & 1200 \\
6 & $\varnothing$ & $\varnothing$ & 11212213 & 0311 \\
7 & 12112213, 11221312 & 0301, 0330
& 11211223, 12112231, 11223112
& 1311, 0201, 0320 \\
8 & 12111223, 11122312, 11223121 & 1301, 1330, 0220
& 12111232, 11123212,
11232121 & 1201, 1320, 0210 \\
9 & 11123122, 11231221 & 1220, 0110 & 12131221, 11213122
& 0100, 1210\\
\hline
\end{tabular}
\end{center}
\normalsize
Here the riggings denote $r^{(1)}_{1,1} r^{(1)}_{2,2} r^{(1)}_{2,1}
r^{(2)}_{1,1}$.
\end{example}

\section{Inverse scattering method}\label{k:sec:ist}

\subsection{Action-angle variables}\label{k:sec:aav}
By action-angle variables for the periodic $A^{(1)}_n$ SCA,
we mean the variables or combinatorial objects
that are conserved (action) or growing linearly (angle)
under the commuting family of time evolutions $\{T^{(r)}_l\}$.
They are scattering data in the context of
inverse scattering method~\mbox{\cite{GGKM, DT, DMN}}.
In our approach, the action-angle variables are constructed by
a suitable extension of rigged conf\/igurations,
which exploits a connection to the combinatorial Bethe ans{\"a}tze
both at \mbox{$q=1$}~\cite{KKR} and $q=0$ \cite{KN}.
These features have been fully worked out in~\cite{KTT} for $A^{(1)}_1$.
More recently it has been shown further that the set of angle variables
can be decomposed into connected components and
every such component is a torus \cite{T}.
Here we present a
conjectural generalization of these results to $A^{(1)}_n$ case.
It provides a conceptual explanation of the
dynamical period and the state counting formula proposed in \cite{KT1, KT2}.

Consider the level set ${\mathcal P}(\mu)$ (\ref{k:pmu})
with $\mu=(\mu^{(1)}, \ldots, \mu^{(n)})$
satisfying the condition (\ref{k:pp}).
The action variable for any path $p \in {\mathcal P}(\mu)$ is
def\/ined to be $\mu$ itself.
By Proposition \ref{k:pr:tte}, it is a~conserved quantity
under any time evolution.

{\sloppy Recall that each block $(a,i)$ of a rigged conf\/iguration is assigned with
the rigging $r^{(a)}_{i,\alpha}$ with $\alpha=1,\ldots, m^{(a)}_i$.
We extend $r^{(a)}_{i,\alpha}$ to $\alpha \in {\mathbb Z}$ uniquely
so that the quasi-periodicity
\begin{equation}\label{k:qpr}
r_{i,\alpha + m^{(a)}_i}^{(a)} = r_{i,\alpha}^{(a)} + p^{(a)}_i
\qquad (\alpha \in {\mathbb Z})
\end{equation}
is fulf\/illed.
Such a sequence ${\bf r}=(r^{(a)}_{i,\alpha})_{\alpha \in {\mathbb Z}}$
will be called the {\em quasi-periodic extension} of
$(r^{(a)}_{i,\alpha})_{1 \le i \le m^{(a)}_i}$.

}

Set
\begin{gather}
\tilde{\mathcal J}(\mu) =
\prod_{(ai) \in \overline{H}}
\tilde{\Lambda}\big(m^{(a)}_i, p^{(a)}_i\big),\nonumber\\
\tilde{\Lambda}(m, p) =
\{ (\lambda_\alpha)_{\alpha \in \mathbb{Z}} \mid
\lambda_\alpha \in {\mathbb Z}, \;
\lambda_\alpha \le \lambda_{\alpha+1}, \;
\lambda_{\alpha+m}\! =\! \lambda_\alpha\! + p \;
\hbox{ for all } \alpha \}\quad
(m\ge 1, p\ge 0),\!\!\!\label{k:lti}
\end{gather}
where $\prod$ denotes a direct product of sets.
Using Theorem \ref{k:th:rc}, we def\/ine an injection
\begin{equation*}
\begin{split}
\iota: {\mathcal P}_+(\mu) \overset{\phi}{\longrightarrow}
\text{RC}(\mu) \;\;&\longrightarrow \;\;\tilde{\mathcal J}(\mu)\\
p_+ \;\; \longrightarrow \;\;(\mu, {\bf r}) \;\;
&\longmapsto \;
\big(r^{(a)}_{i,\alpha}\big)_{(ai) \in \overline{H}, \alpha \in {\mathbb Z}}\,,
\end{split}
\end{equation*}
where the bottom right is the
quasi-periodic extension of the original rigging
${\bf r}=\big(r^{(a)}_{i,\alpha}\big)_{(ai\alpha) \in H}$
in $(\mu,{\bf r})$ as explained above.
Elements of $\tilde{\mathcal J}(\mu)$ will be called
{\em extended rigged configurations}.
One may still view an extended rigged conf\/iguration
as the right diagram in (\ref{k:rca0}).
The only dif\/ference from the original rigged conf\/iguration is that
the riggings $r^{(a)}_{i,\alpha}$ outside the range
$1\le \alpha \le m^{(a)}_i$ have also been f\/ixed by the quasi-periodicity
(\ref{k:qpr}).
In what follows, either the rigging
$\big(r^{(a)}_{i,\alpha}\big)_{(ai\alpha) \in H}$ or
its quasi-periodic extension
$\big(r^{(a)}_{i,\alpha}\big)_{(ai) \in \overline{H}, \alpha \in {\mathbb Z}}
\in \tilde{\mathcal J}(\mu)$
will be denoted by the same symbol ${\bf r}$.

Let ${\mathcal T}$ be the abelian group generated by all the time evolutions
$T^{(a)}_l$ with $1 \le a \le n$ and $l \ge 1$.
Let further ${\mathcal A}$ be a free abelian group generated by
the symbols $s^{(a)}_i$ with $(a,i) \in \overline{H}$.
We consider their commutative actions on $\tilde{\mathcal J}(\mu)$
as follows:
\begin{gather}
 T^{(a)}_l : \
\big(r^{(b)}_{j,\beta}\big)_{(bj) \in \overline{H}, \beta \in {\mathbb Z}}
\ \mapsto \
\big(r^{(b)}_{j,\beta}
+\delta_{ab}\min\big(l,l^{(b)}_j\big)\big)_{(bj) \in \overline{H}, \beta \in {\mathbb Z}},
\label{k:ter2}\\
s^{(a)}_i : \ \;
\big(r^{(b)}_{j,\beta}\big)_{(bj) \in \overline{H}, \beta \in {\mathbb Z}}
\  \mapsto \
\big(r^{(b)}_{j,\beta+\delta_{ab}\delta_{ij}}
+C_{ab}\min\big(l^{(a)}_i,l^{(b)}_j\big)\big)_{(bj) \in \overline{H}, \beta \in {\mathbb Z}}.
\label{k:sai}
\end{gather}
The generator $s^{(a)}_i$
is the $A^{(1)}_n$ version of the slide introduced in \cite{KTT}
for $n=1$.
Now we def\/ine
\begin{equation}\label{k:jdef}
{\mathcal J}(\mu) = \tilde{\mathcal J}(\mu)/ {\mathcal A},
\end{equation}
which is the set of all ${\mathcal A}$-orbits,
whose elements are
written as ${\mathcal A} \cdot {\bf r}$ with
${\bf r} \in \tilde{\mathcal J}(\mu)$.
Elements of ${\mathcal J}(\mu)$
will be called {\em angle variables}.
See also (\ref{k:waru}).
Since ${\mathcal T}$ and ${\mathcal A}$ act on
$\tilde{\mathcal J}(\mu)$ commutatively,
there is a natural action of ${\mathcal T}$ on
${\mathcal J}(\mu)$.
For $t \in {\mathcal T}$ and
$y = {\mathcal A} \cdot {\bf r} \in {\mathcal J}(\mu)$,
the action is given by
$t \cdot y = {\mathcal A} \cdot (t \cdot {\bf r})$.
In short, the time evolution of extended rigged conf\/igurations
(\ref{k:ter2}) naturally induces the time evolution of angle variables.

\begin{remark}\label{k:re:t11}
{}From (\ref{k:sai}) and (\ref{k:pdef}),
one can easily check
\begin{equation*}
\Bigg(\prod_{(ai)\in \overline{H}}s^{(a)\,m^{(a)}_i}_i\Bigg)({\bf r})
=\big(T^{(1)}_1\big)^L({\bf r}) =
\big(L\delta_{b1}+
r^{(b)}_{j,\beta}\big)_{(bj) \in \overline{H}, \beta \in {\mathbb Z}}
\end{equation*}
for ${\bf r}=\big(r^{(b)}_{j,\beta}\big)_{(bj) \in \overline{H}, \beta \in {\mathbb Z}}
\in \tilde{\mathcal J}(\mu)$.
Therefore ${\bf r}$ and $\big(T^{(1)}_1\big)^L({\bf r})$ def\/ine the
same angle variable.
\end{remark}

\subsection{Linearization of time evolution}\label{k:sec:lt}

Given the level set ${\mathcal P}(\mu)$ with $\mu$
satisfying the condition (\ref{k:pp}), we are going to introduce the
bijection $\Phi$ to the set of angle variables ${\mathcal J}(\mu)$.
Recall that $\Sigma(p)$ denotes the connected component
of ${\mathcal P}(\mu)$ that involves the path $p$.

\begin{conjecture}\label{t:conj:aug10_5}
Under the condition \eqref{k:pp},
$\Sigma (p) \cap {\mathcal P}_+(\mu) \ne \varnothing$ holds
for any $p \in {\mathcal P}(\mu)$.
\end{conjecture}

This implies that any path $p$ can be expressed in the form
$p= t\cdot p_+$ by some highest path
$p_+ \in {\mathcal P}_+(\mu)$ and time evolution
$t = \prod (T^{(r)}_{l})^{d^{(r)}_l}\in {\mathcal T}$.\footnote{We do not restrict $t$ to a specif\/ic subgroup of ${\mathcal T}$.}
Although such an expression is not unique in general,
one can go from ${\mathcal P}(\mu)$
to ${\mathcal J}(\mu)$ along the following scheme:
\begin{equation}\label{k:ds}
\begin{split}
\Phi: \;\;{\mathcal P}(\mu) & \; \longrightarrow\;
{\mathcal T} \times{\mathcal P}_+(\mu)
\;\overset{\text{id}\times \iota}{\longrightarrow}
\;{\mathcal T} \times \tilde{\mathcal J}(\mu) \;
\longrightarrow  \quad{\mathcal J}(\mu)\\
p \quad & \;\,\longmapsto \;\;\;(t,  p_+)
\;\;\quad\longmapsto\;\;\;\;
\bigl(t, {\bf r}\bigr) \quad\;\;\longmapsto \;
{\mathcal A}\cdot \bigl(t\cdot {\bf r}\bigr).
\end{split}
\end{equation}
Note that the abelian group
${\mathcal T}$ of time evolutions acts on paths by (\ref{k:te}) and
on extended rigged conf\/igurations by (\ref{k:ter2}).
The composition (\ref{k:ds}) serves as a def\/inition of a map $\Phi$
only if the non-uniqueness of the
decomposition $p \mapsto (t, p_+)$ is
`canceled' by regarding $t\cdot {\bf r}$
mod ${\mathcal A}$.

Our main conjecture in this subsection is the following.

\begin{conjecture}\label{t:conj:aug21_5}
Suppose the condition \eqref{k:pp} is satisfied.
Then the map $\Phi$ \eqref{k:ds} is well-defined,
bijective and commutative with the action of ${\mathcal T}$.
Namely, the following diagram is commutative
\begin{gather}\label{k:cd1}
\begin{CD}
{\mathcal P}(\mu) @>{\Phi}>> {\mathcal J}(\mu) \\
@V{\mathcal T}VV @VV{\mathcal T}V\\
{\mathcal P}(\mu) @>{\Phi}>> {\mathcal J}(\mu)
\end{CD}
\end{gather}
\end{conjecture}

This is consistent with the obvious
periodicity $(T^{(1)}_1)^L(p)=p$
of paths under the cyclic shift (\ref{k:shf}) and Remark \ref{k:re:t11}.
The time evolution of the angle variable (\ref{k:ter2}) is linear.
Thus the commutative diagram (\ref{k:cd1}) transforms
the nonlinear dynamics of paths into straight motions with
various velocity vectors.
The set ${\mathcal J}(\mu)$ is a tropical analogue of
Jacobi variety in the theory of quasi-periodic solutions of
soliton equations \cite{DT, DMN}.
The map $\Phi$ whose essential ingredient is
the KKR bijection $\phi$ plays the role of the Abel--Jacobi map.
Conceptually, the solution of the initial value problem is simply stated as
$t^{\mathcal N}(p)
= \Phi^{-1}\circ t^{\mathcal N} \circ \Phi(p)$ for any $t \in {\mathcal T}$,
where $t$ in the right hand side is linear.
Practical calculations can be found in
Examples \ref{k:ex:aa2} and \ref{k:ex:ivp}.
In Section \ref{k:sec:trt}, we shall present an explicit formula
for the tropical Jacobi-inversion $\Phi^{-1}$
in terms of a tropical Riemann theta function.
The result of this sort was f\/irst formulated and proved
in~\cite{KTT, KS1, KS2} for $n=1$.

\begin{example}\label{k:ex:aa}
Take an evolvable path
\begin{equation}\label{k:exp}
p= 211332111321133112221112 \in \big(B^{1,1}\big)^{\otimes 24},
\end{equation}
whose energies are given as
\begin{equation*}
E^{(1)}_1 = 5,  \quad E^{(1)}_2 = 10, \quad E^{(1)}_{\geq 3} = 12, \quad
E^{(2)}_1 = 2, \quad E^{(2)}_2 = 3, \quad E^{(2)}_3 = 4,  \quad
E^{(2)}_{\geq 4} = 5.
\end{equation*}
From this and (\ref{k:eee})
we f\/ind $p \in {\mathcal P}(\mu)$
with $\mu=((33222),(41))$.
The cyclic shift $(T^{(1)}_1)^j(p)$
is not highest for any $j$.
However it can be made into a highest path
$p_+$ and transformed to a~rigged configuration as follows:
\begin{equation}\label{k:rc1}
\begin{picture}(200,85)(-220,17)
\setlength{\unitlength}{0.5mm}

\put(-146,40)
{$p_+ = (T^{(1)}_1)^3(T^{(2)}_1)^3(p)$}
\put(-135,25)
{$ = 111221113221132113311322$}

\put(-22,32){$\overset{\phi}{\longmapsto}$}

\put(7,68){$(\mu^{(1)}, r^{(1)})$}
\put(72,68){$(\mu^{(2)}, r^{(2)})$}
\multiput(0,0)(0,-10){3}{
\put(10,60){\line(1,0){30}}}
\multiput(0,0)(0,10){4}{
\put(10,10){\line(1,0){20}}}
\multiput(0,0)(10,0){3}{
\put(10,60){\line(0,-1){50}}}
\put(40,60){\line(0,-1){20}}

\put(3,47){$4$}\put(3,22){$7$}
\put(42,52){$4$}
\put(42,42){$2$}
\put(32,32){$6$}
\put(32,22){$5$}
\put(32,12){$1$}

\multiput(-30,0)(0,0){1}{
\put(93,52){$2$}\put(93,42){$1$}
\put(143,52){$0$}\put(113,42){$0$}
\put(100,60){\line(1,0){40}}
\put(100,50){\line(1,0){40}}
\put(100,40){\line(1,0){10}}
\put(100,60){\line(0,-1){20}}
\put(110,60){\line(0,-1){20}}
\multiput(0,0)(10,0){3}{
\put(120,60){\line(0,-1){10}}}}

\end{picture}
\end{equation}
Thus action variable
$(\mu^{(1)}, \mu^{(2)})$ indeed coincides with
$\mu=((33222), (41))$ obtained from the energy
in agreement with Theorem \ref{k:th:rc}.
Although it is not necessary, we have exhibited the vacancy numbers
on the left of $\mu^{(1)}$ and $\mu^{(2)}$
for convenience.
We list the relevant data:

\vspace{2mm}
\begin{table}[h]
\begin{center}
\begin{tabular}{c|ccc|cc}
\hline
$(ai)$ & $l^{(a)}_i$ & $m^{(a)}_i$ & $p^{(a)}_i$ &
$r^{(a)}_{i,\alpha}$ & $\gamma^{(a)}_i$\\
\hline
(11) & 3 & 2 & 4 & $2,4$& 2\\
(12) & 2 & 3 & 7 & $1,5,6$ & 1\\
(21) & 4 & 1 & 2 & 0 & 1\\
(22) & 1 & 1 & 1 & 0 & 1\\
\hline
\end{tabular}
\end{center}
\end{table}

The index set $H$ (\ref{k:hdef})
corresponding to $\mu=(\mu^{(1)}, \mu^{(2)})$
in (\ref{k:rc1}) is
\begin{equation}\label{k:hido}
H = \{(111),(112),(121),(122),(123),(211), (221)\},
\end{equation}
which has the cardinality $G = 7$.
The quantity $\gamma^{(a)}_i$ is
the order of symmetry which will be explained
in Section \ref{k:sec:dcc}.
The matrix $F$ (\ref{k:f}) reads
\begin{gather}
F  =
\left(
\begin{array}{rrrr}
p^{(1)}_1 + 6 m^{(1)}_1 & 4 m^{(1)}_2 & -3 m^{(2)}_1& - m^{(2)}_2\\
4 m^{(1)}_1 & p^{(1)}_2 + 4 m^{(1)}_2 & -2 m^{(2)}_1 & - m^{(2)}_2\\
-3 m^{(1)}_1 & -2 m^{(1)}_2 & p^{(2)}_1 + 8 m^{(2)}_1 & 2 m^{(2)}_2\\
- m^{(1)}_1 & - m^{(1)}_2 & 2 m^{(2)}_1 & p^{(2)}_2 + 2 m^{(2)}_2
\end{array}
\right)\nonumber\\
\phantom{F}{} =
\begin{pmatrix}
16 & 12 & -3 & -1 \\
8 & 19 & -2 & -1 \\
-6 & -6 & 10 & 2 \\
-2 & -3 & 2 & 3
\end{pmatrix}.\label{t:eq:aug18_4}
\end{gather}

We understand that the riggings here are parts of
extended one obeying (\ref{k:qpr}).
Then from~(\ref{k:rc1}) and (\ref{k:ter2}),
the extended rigged conf\/iguration for $p$ (\ref{k:exp})
is the left hand side of
\begin{equation}\label{k:rc15}
\begin{picture}(200,75)(-220,13)
\setlength{\unitlength}{0.5mm}

\multiput(-150,0)(0,0){1}{
\multiput(0,0)(0,-10){3}{
\put(10,60){\line(1,0){30}}}
\multiput(0,0)(0,10){4}{
\put(10,10){\line(1,0){20}}}
\multiput(0,0)(10,0){3}{
\put(10,60){\line(0,-1){50}}}
\put(40,60){\line(0,-1){20}}

\put(42,52){$1$}
\put(42,42){$-1$}
\put(32,32){$3$}
\put(32,22){$2$}
\put(32,12){$-2$}

\multiput(-30,0)(0,0){1}{
\put(143,52){$-3$}\put(113,42){$-3$}
\put(100,60){\line(1,0){40}}
\put(100,50){\line(1,0){40}}
\put(100,40){\line(1,0){10}}
\put(100,60){\line(0,-1){20}}
\put(110,60){\line(0,-1){20}}
\multiput(0,0)(10,0){3}{
\put(120,60){\line(0,-1){10}}}}}

\put(-15,50){$s^{(1)}_2$}
\put(-15,40){$\simeq$}

\multiput(0,0)(0,-10){3}{
\put(10,60){\line(1,0){30}}}
\multiput(0,0)(0,10){4}{
\put(10,10){\line(1,0){20}}}
\multiput(0,0)(10,0){3}{
\put(10,60){\line(0,-1){50}}}
\put(40,60){\line(0,-1){20}}

\put(42,52){$5$}
\put(42,42){$3$}
\put(32,32){$9$}
\put(32,22){$7$}
\put(32,12){$6$}

\multiput(-30,0)(0,0){1}{
\put(143,52){$-5$}\put(113,42){$-4$}
\put(100,60){\line(1,0){40}}
\put(100,50){\line(1,0){40}}
\put(100,40){\line(1,0){10}}
\put(100,60){\line(0,-1){20}}
\put(110,60){\line(0,-1){20}}
\multiput(0,0)(10,0){3}{
\put(120,60){\line(0,-1){10}}}}

\end{picture}
\end{equation}
Here the transformation to the right hand side
demonstrates an example of identif\/ication by~${\mathcal A}$.
These two objects are examples of representative elements of
the angle variable for~$p$~(\ref{k:exp}).
We will calculate $(T^{(a)}_l)^{1000}(p)$ in Example~\ref{k:ex:ivp}
using the reduced angle variable
that will be introduced in the next subsection.
\end{example}

\begin{example}\label{k:ex:aa2}
Consider the evolvable path
\begin{equation}\label{k:exp2}
p = 321113211222111223331111
\in \big(B^{1,1}\big)^{\otimes 24},
\end{equation}
which is same as in Example~\ref{k:ex:int}.
The cyclic shift $(T^{(1)}_1)^j(p)$
is not highest for any $j$.
However it can be made into a highest path
$p_+$ and transformed to a rigged conf\/iguration as
\begin{equation}\label{k:rc3}
\begin{picture}(200,75)(-220,13)
\setlength{\unitlength}{0.5mm}

\put(-146,40)
{$p_+ = (T^{(1)}_1)^{21}T^{(2)}_1(p)$}
\put(-135,25)
{$ = 111221322111221332111331$}

\put(-22,32){$\overset{\phi}{\longmapsto}$}

\multiput(0,0)(0,-10){3}{
\put(10,60){\line(1,0){30}}}
\multiput(0,0)(0,10){4}{
\put(10,10){\line(1,0){20}}}
\multiput(0,0)(10,0){3}{
\put(10,60){\line(0,-1){50}}}
\put(40,60){\line(0,-1){20}}

\put(3,47){$4$}\put(3,22){$7$}
\put(42,52){$1$}
\put(42,42){$0$}
\put(32,32){$6$}
\put(32,22){$3$}
\put(32,12){$1$}

\multiput(-30,0)(0,0){1}{
\put(93,52){$2$}\put(93,42){$1$}
\put(143,52){$2$}\put(113,42){$0$}
\put(100,60){\line(1,0){40}}
\put(100,50){\line(1,0){40}}
\put(100,40){\line(1,0){10}}
\put(100,60){\line(0,-1){20}}
\put(110,60){\line(0,-1){20}}
\multiput(0,0)(10,0){3}{
\put(120,60){\line(0,-1){10}}}}

\end{picture}
\end{equation}
Thus $p$ (\ref{k:exp2})
belongs to the same level set ${\mathcal P}(\mu)$
with $\mu=((33222),(41))$
as Example \ref{k:ex:aa}, hence the matrix $F$ (\ref{k:f})
is again given by (\ref{t:eq:aug18_4}).
However,
it has a dif\/ferent (trivial) order of symmetry $\gamma^{(a)}_i$
given in the following table:
\vspace{2mm}

\begin{center}
\begin{tabular}{c|ccc|cc}
\hline
$(ai)$ & $l^{(a)}_i$ & $m^{(a)}_i$ & $p^{(a)}_i$ & $r^{(a)}_{i,\alpha}$ & $\gamma^{(a)}_i$\\
\hline
(11) & 3 & 2 & 4 & $0,1$& 1\\
(12) & 2 & 3 & 7 & $1,3,6$ & 1\\
(21) & 4 & 1 & 2 & 2 & 1\\
(22) & 1 & 1 & 1 & 0 & 1\\
\hline
\end{tabular}
\end{center}

Regard (\ref{k:rc3}) as a part of extended rigged conf\/iguration
as in Example \ref{k:ex:aa}.
Then from (\ref{k:ter2}),
the extended rigged conf\/iguration for $p$ (\ref{k:exp2})
is the left hand side of
\begin{equation}\label{k:rc35}
\begin{picture}(200,75)(-220,13)
\setlength{\unitlength}{0.5mm}

\multiput(-150,0)(0,0){1}{
\multiput(0,0)(0,-10){3}{
\put(10,60){\line(1,0){30}}}
\multiput(0,0)(0,10){4}{
\put(10,10){\line(1,0){20}}}
\multiput(0,0)(10,0){3}{
\put(10,60){\line(0,-1){50}}}
\put(40,60){\line(0,-1){20}}

\put(42,52){$-20$}
\put(42,42){$-21$}
\put(32,32){$-15$}
\put(32,22){$-18$}
\put(32,12){$-20$}

\multiput(-30,0)(0,0){1}{
\put(143,52){$1$}\put(113,42){$-1$}
\put(100,60){\line(1,0){40}}
\put(100,50){\line(1,0){40}}
\put(100,40){\line(1,0){10}}
\put(100,60){\line(0,-1){20}}
\put(110,60){\line(0,-1){20}}
\multiput(0,0)(10,0){3}{
\put(120,60){\line(0,-1){10}}}}}

\put(-15,50){$s^{(2)}_1$}
\put(-15,40){$\simeq$}

\multiput(0,0)(0,-10){3}{
\put(10,60){\line(1,0){30}}}
\multiput(0,0)(0,10){4}{
\put(10,10){\line(1,0){20}}}
\multiput(0,0)(10,0){3}{
\put(10,60){\line(0,-1){50}}}
\put(40,60){\line(0,-1){20}}

\put(42,52){$-23$}
\put(42,42){$-24$}
\put(32,32){$-17$}
\put(32,22){$-20$}
\put(32,12){$-22$}

\multiput(-30,0)(0,0){1}{
\put(143,52){$11$}\put(113,42){$1$}
\put(100,60){\line(1,0){40}}
\put(100,50){\line(1,0){40}}
\put(100,40){\line(1,0){10}}
\put(100,60){\line(0,-1){20}}
\put(110,60){\line(0,-1){20}}
\multiput(0,0)(10,0){3}{
\put(120,60){\line(0,-1){10}}}}

\end{picture}
\end{equation}
Here again the transformation to the right hand side
demonstrates an example of identif\/ication by~${\mathcal A}$.
These two objects are examples of representative elements of
the angle variable for~$p$~(\ref{k:exp2}).

To illustrate the solution of the initial value problem
along the inverse scheme (\ref{k:cd1}), we derive
\begin{gather}\label{k:rei}
\big(T^{(1)}_3\big)^{1000}(p) = 221132211331111321322111.
\end{gather}
By applying (\ref{k:ter2}) to the left hand side of (\ref{k:rc35}),
the angle variable for
$\big(T^{(1)}_3\big)^{1000}(p)$ is obtained as
\begin{equation*}
\begin{picture}(200,75)(-220,13)
\setlength{\unitlength}{0.5mm}

\multiput(-150,0)(0,0){1}{
\multiput(0,0)(0,-10){3}{
\put(10,60){\line(1,0){30}}}
\multiput(0,0)(0,10){4}{
\put(10,10){\line(1,0){20}}}
\multiput(0,0)(10,0){3}{
\put(10,60){\line(0,-1){50}}}
\put(40,60){\line(0,-1){20}}

\put(42,52){$2980$}
\put(42,42){$2979$}
\put(32,32){$1985$}
\put(32,22){$1982$}
\put(32,12){$1980$}

\multiput(-30,0)(0,0){1}{
\put(143,52){$1$}\put(113,42){$-1$}
\put(100,60){\line(1,0){40}}
\put(100,50){\line(1,0){40}}
\put(100,40){\line(1,0){10}}
\put(100,60){\line(0,-1){20}}
\put(110,60){\line(0,-1){20}}
\multiput(0,0)(10,0){3}{
\put(120,60){\line(0,-1){10}}}}}

\put(-15,40){$\simeq$}

\multiput(0,0)(0,-10){3}{
\put(10,60){\line(1,0){30}}}
\multiput(0,0)(0,10){4}{
\put(10,10){\line(1,0){20}}}
\multiput(0,0)(10,0){3}{
\put(10,60){\line(0,-1){50}}}
\put(40,60){\line(0,-1){20}}

\put(42,52){$9$}
\put(42,42){$8$}
\put(32,32){$11$}
\put(32,22){$9$}
\put(32,12){$7$}

\multiput(-30,0)(0,0){1}{
\put(143,52){$5$}\put(113,42){$1$}
\put(100,60){\line(1,0){40}}
\put(100,50){\line(1,0){40}}
\put(100,40){\line(1,0){10}}
\put(100,60){\line(0,-1){20}}
\put(110,60){\line(0,-1){20}}
\multiput(0,0)(10,0){3}{
\put(120,60){\line(0,-1){10}}}}

\end{picture}
\end{equation*}
where the equivalence is seen by using
$\big(s^{(1)}_1\big)^{-358}
\big(s^{(1)}_2\big)^{-136}
\big(s^{(2)}_1\big)^{-117}
\big(s^{(2)}_2\big)^{-86} \in {\mathcal A}$.
The right hand side is the angle variable of
$\big(T^{(1)}_1\big)^7 T^{(2)}_4(p'_+)$, where
$p'_+$ is the highest path obtained by the KKR map
\begin{equation*}
\begin{picture}(200,75)(-220,13)
\setlength{\unitlength}{0.5mm}

\multiput(-150,0)(0,0){1}{
\multiput(0,0)(0,-10){3}{
\put(10,60){\line(1,0){30}}}
\multiput(0,0)(0,10){4}{
\put(10,10){\line(1,0){20}}}
\multiput(0,0)(10,0){3}{
\put(10,60){\line(0,-1){50}}}
\put(40,60){\line(0,-1){20}}

\put(42,52){$2$}
\put(42,42){$1$}
\put(32,32){$4$}
\put(32,22){$2$}
\put(32,12){$0$}

\multiput(-30,0)(0,0){1}{
\put(143,52){$1$}\put(113,42){$0$}
\put(100,60){\line(1,0){40}}
\put(100,50){\line(1,0){40}}
\put(100,40){\line(1,0){10}}
\put(100,60){\line(0,-1){20}}
\put(110,60){\line(0,-1){20}}
\multiput(0,0)(10,0){3}{
\put(120,60){\line(0,-1){10}}}}}

\put(-22,35){$\overset{\phi^{-1}}{\longmapsto}$}

\put(0,35){$p'_+ =112211132112221331133211$}

\end{picture}
\end{equation*}
Calculating $\big(T^{(1)}_1\big)^7 T^{(2)}_4(p'_+)$, one f\/inds
(\ref{k:rei}).
\end{example}

\subsection{Decomposition into connected components}\label{k:sec:dcc}

As remarked under Conjecture \ref{t:conj:aug10_3},
the level set ${\mathcal P}(\mu)$ is a disjoint union of several
connected components $\Sigma(p)$
under the time evolution ${\mathcal T}$.
It is natural to decompose the commutative
diagram (\ref{k:cd1}) further into those connected components, i.e.,
${\mathcal T}$-orbits,
and seek the counterpart of $\Sigma(p)$ in ${\mathcal J}(\mu)$.
To do this one needs a precise description
of the internal symmetry of angle variables ref\/lecting
a certain commensurability of soliton conf\/iguration in a path.
In addition, it is necessary to separate
the angle variables into two parts,
one recording the internal symmetry
(denoted by $\boldsymbol{\lambda}$)\footnote{This $\boldsymbol{\lambda}$ may be viewed as an additional
conserved quantity.}
and the other accounting for the straight motions
(denoted by $\omega$).
We call the latter part {\em reduced angle variables}.
The point is that action of ${\mathcal T}$ becomes transitive
by switching to the reduced angle variables from
angle variables thereby allowing us to describe the
connected components and their multiplicity explicitly.
These results are generalizations of the $n=1$ case~\cite{T}.

Let $m\ge 1$ and $p\ge 0$.\footnote{This $p$ is not a path of course
and will set to be a vacancy number shortly.}
We introduce the following set by imposing
just a simple extra condition on $\tilde{\Lambda}(m,p)$ (\ref{k:lti}):
\begin{equation}\label{t:eq:may1_1}
\Lambda(m, p) =
\{ (\lambda_\alpha)_{\alpha \in \mathbb{Z}} \mid
\lambda_\alpha \in {\mathbb Z}, \;\lambda_1=0,\;
\lambda_\alpha \le \lambda_{\alpha+1}, \;
\lambda_{\alpha+m} = \lambda_\alpha + p \;
\hbox{ for all } \alpha \}.
\end{equation}
This is a f\/inite set with cardinality
\begin{equation}\label{k:cla}
|\Lambda (m,p)| = {p+m-1 \choose m-1}.
\end{equation}

Given a conf\/iguration $\mu$, specify the data
$m^{(a)}_i$, $l^{(a)}_i$, $p^{(a)}_i$, etc by (\ref{k:rca0})--(\ref{k:gdef}).
Set
\begin{equation*}
X= X(\mu) = X^1 \times X^2,\qquad
X^1 = {\mathbb Z}^g,\qquad
X^2 = \prod_{(ai) \in \overline{H}} \Lambda \big(m^{(a)}_i,p^{(a)}_i\big),
\end{equation*}
and consider the splitting of extended rigged conf\/igurations
into $X^1$ and $X^2$ as follows:
\begin{equation}\label{k:spl}
\begin{split}
\tilde{\mathcal J}(\mu)\ \ \ \quad&{\longleftrightarrow} \quad
X^1 \times X^2\\
\big(r^{(a)}_{i,\alpha}\big)_{(ai)\in \overline{H}, \alpha \in {\mathbb Z}}
\; &\longleftrightarrow  \quad
(\omega, \boldsymbol{\lambda})\quad\quad
\omega^{(a)}_i = r^{(a)}_{i,1},\qquad
\lambda^{(a)}_{i,\alpha} = r^{(a)}_{i,\alpha}-r^{(a)}_{i,1},
\end{split}
\end{equation}
where $\omega = \big(\omega^{(a)}_i\big)_{(ai)\in \overline{H}}$
and
$\boldsymbol{\lambda} = \big(\lambda^{(a)}_{i,\alpha}\big)_{
(ai) \in \overline{H}, \alpha \in {\mathbb Z}}$.

Recall that ${\mathcal T}$  and ${\mathcal A}$
are the abelian groups acting on
$\tilde{\mathcal J}(\mu)$ as in (\ref{k:ter2})--(\ref{k:sai}).
By the one to one correspondence (\ref{k:spl}),
they also act on $X$ commutatively as
\begin{gather}
T^{(a)}_l (\omega , \boldsymbol{\lambda})
= \big(\omega + h^{(a)}_l , \boldsymbol{\lambda}\big),
\qquad
h^{(a)}_l = \big(\delta_{ab} \min \big(l, l^{(b)}_j\big)\big)_{(bj) \in \overline{H}},
\label{k:ter}\\
s^{(a)}_i
(\omega , \boldsymbol{\lambda}) =
(\omega' , \boldsymbol{\lambda}'),\qquad
\begin{cases}
\omega' &= \omega + \big(
\lambda^{(b)}_{j,1+ \delta_{ij}\delta_{ab}} +
C_{ab} \min \big(l^{(a)}_i,l^{(b)}_j\big) \big)_{(bj) \in \overline{H} },\\
\boldsymbol{\lambda}' &=
\big(\lambda^{(b)}_{j,\alpha +\delta_{ij}\delta_{ab}} -
\lambda^{(b)}_{j,1 + \delta_{ij}\delta_{ab}}\big)_{(bj) \in \overline{H},
\alpha \in {\mathbb Z}}. \end{cases}\label{t:eq:aug21_1}
\end{gather}
In particular, ${\mathcal T}$  acts on $X^1$ transitively leaving
$X^2$ unchanged.

{}From (\ref{k:jdef}) and (\ref{k:spl}),
the angle variables may also be viewed as elements of $X/{\mathcal A}$.
It is the set of all ${\mathcal A}$-orbits, whose elements are
written as ${\mathcal A} \cdot x$ with $x \in X$.
Since ${\mathcal T}$ and ${\mathcal A}$ act on $X$ commutatively,
there is a natural action of ${\mathcal T}$ on $X / {\mathcal A}$.
For $t \in {\mathcal T}$ and
$y = {\mathcal A} \cdot x \in X / {\mathcal A}$,
the action is given by $t \cdot y = {\mathcal A} \cdot (t \cdot x)$.
For any representative element
$(\omega, \boldsymbol{\lambda}) \in X/{\mathcal A}$
of an angle variable,
we call the $X^1$ part $\omega$ a {\em reduced angle variable}.
Under the simple bijective correspondence (\ref{k:spl}),
the map (\ref{k:ds}) is rephrased as
\begin{equation*}
\begin{split}
\Phi: \;\;{\mathcal P}(\mu) & \; \longrightarrow\;
{\mathcal T} \times {\mathcal P}_+(\mu)
\;\overset{\text{id}\times \iota}{\longrightarrow}\;
\quad{\mathcal T} \times X \quad\;
\longrightarrow \quad \;\; X/{\mathcal A}\\
p \quad & \;\,\longmapsto \;\;\;(t,  p_+)
\;\;\quad\longmapsto\;\;
\bigl(t, (\omega, \boldsymbol{\lambda})\bigr)
\;\;\;\longmapsto \;
{\mathcal A}\cdot \bigl(t\cdot(\omega, \boldsymbol{\lambda})\bigr),
\end{split}
\end{equation*}
where we have used the same symbol $\Phi$.
The conjectural commutative diagram (\ref{k:cd1}) becomes
\begin{equation}\label{k:cd15}
\begin{CD}
{\mathcal P}(\mu) @>{\Phi}>> X/{\mathcal A} \\
@V{\mathcal T}VV @VV{\mathcal T}V\\
{\mathcal P}(\mu) @>{\Phi}>> X/{\mathcal A}
\end{CD}
\end{equation}

Now we are ready to describe the internal symmetry
of angle variables and the resulting decomposition.
We introduce a ref\/inement of $\Lambda(m,p)$ (\ref{t:eq:may1_1})
as follows:
\begin{equation}\label{t:eq:may12_1}
\Lambda_\gamma(m,p) = \left\{
\lambda \in \Lambda
\left(\frac{m}{\gamma}, \frac{p}{\gamma}\right) \Big|
\lambda \notin \Lambda
\left(\frac{m}{\gamma'}, \frac{p}{\gamma'}\right) \
\mbox{for any} \, \gamma' > \gamma \right\},
\end{equation}
where $\gamma$ is a (not necessarily largest)
common divisor of $m$ and $p$.
In other words $\Lambda_\gamma(m,p)$ is the set of all arrays
$(\lambda_\alpha)_{\alpha \in \mathbb{Z}} \in \Lambda (m,p)$
that satisfy the reduced quasi-periodicity
$\lambda_{\alpha+m/\gamma}
= \lambda_\alpha + p/\gamma$ but do not satisfy the same relation
when $\gamma$ is replaced by a larger $\gamma'$.
Such $\gamma$ is called the {\em order of symmetry}.
By the def\/inition one has the disjoint union decomposition:
\begin{equation*}
\Lambda(m,p) = \bigsqcup_\gamma \Lambda_\gamma(m,p),
\end{equation*}
where $\gamma$ extends over all the common divisors of $m$ and $p$.
Taking the cardinality of this relation
using (\ref{k:cla}) amounts to the identity
\begin{equation}\label{t:eq:aug19_2}
{ p+m-1 \choose m-1 } = \sum_\gamma |\Lambda_\gamma(m,p)|.
\end{equation}

Let $\boldsymbol{\gamma}
= \big(\gamma^{(a)}_i\big)_{(ai) \in \overline{H}}$
be the array of order of symmetry for all blocks.
Thus $\gamma^{(a)}_i$ is any common divisor of
$m^{(a)}_i$ and $p^{(a)}_i$.
We introduce the $g \times g$ matrix $F_{\boldsymbol{\gamma}}$ and the
subsets $X^2_{\boldsymbol{\gamma}} \subset X^2$ and
$X_{\boldsymbol{\gamma}} \subset X$ by
\begin{gather}
F_{\boldsymbol{\gamma}}  = \big(F_{ai, bj}/
\gamma^{(b)}_j\big)_{(ai), (bj) \in \overline{H}},
\label{k:ftdef}\\
X^2_{\boldsymbol{\gamma}}
 = \prod_{(ai) \in \overline{H}}
\Lambda_{\gamma^{(a)}_i}\big(m^{(a)}_i,p^{(a)}_i\big)
\subset X^2, \qquad
X_{\boldsymbol{\gamma}}
= X^1 \times X^2_{\boldsymbol{\gamma}} \subset X,
\nonumber
\end{gather}
where $F_{ai, bj}$ is specif\/ied in (\ref{k:f}).
The matrix $F$ hence $F_{\boldsymbol{\gamma}}$ are positive def\/inite if
$p^{(a)}_i\ge 0$ for any $(a,i) \in \overline{H}$.

\begin{example}\label{k:ex:ga}
Consider Example \ref{k:ex:aa}.
Take the left hand side of (\ref{k:rc15}) as a
representative element
of the angle variable of the path (\ref{k:exp}).
Then the splitting (\ref{k:spl}) gives
the reduced angle variable
$\omega = {}^t(-1,-2,-3,-3)$ and
$\boldsymbol{\lambda}$ is the
quasi-periodic extension of the rigging in
\begin{equation}\label{k:rc22}
\begin{picture}(200,75)(-120,13)
\setlength{\unitlength}{0.5mm}

\multiput(0,0)(0,-10){3}{
\put(10,60){\line(1,0){30}}}
\multiput(0,0)(0,10){4}{
\put(10,10){\line(1,0){20}}}
\multiput(0,0)(10,0){3}{
\put(10,60){\line(0,-1){50}}}
\put(40,60){\line(0,-1){20}}

\put(42,52){$2$}
\put(42,42){$0$}
\put(32,32){$5$}
\put(32,22){$4$}
\put(32,12){$0$}

\multiput(-30,0)(0,0){1}{
\put(143,52){$0$}\put(113,42){$0$}
\put(100,60){\line(1,0){40}}
\put(100,50){\line(1,0){40}}
\put(100,40){\line(1,0){10}}
\put(100,60){\line(0,-1){20}}
\put(110,60){\line(0,-1){20}}
\multiput(0,0)(10,0){3}{
\put(120,60){\line(0,-1){10}}}}

\end{picture}
\end{equation}
Thus $\boldsymbol{\lambda} =
\big(\lambda^{(1)}_1, \lambda^{(1)}_2,
\lambda^{(2)}_1, \lambda^{(2)}_2\big)$ reads
\begin{equation*}
\begin{split}
\lambda^{(1)}_1&=(\ldots, -4,-2, 0,2,4,6,8,10,\ldots),\quad
\lambda^{(1)}_2=(\ldots, -7,-3,-2, 0,4,5,7,11,12,\ldots),\\
\lambda^{(2)}_1&=(\ldots, -2, 0,2,4,\ldots),\quad\qquad\quad\qquad
\lambda^{(2)}_2=(\ldots, -1, 0,1,2,\ldots),
\end{split}
\end{equation*}
with $\lambda^{(a)}_{i,1}=0$.
There is order 2 symmetry $\gamma^{(1)}_1=2$
in the block $(a,i)=(1,1)$ since
$p^{(1)}_1=4$ and $m^{(1)}_1=2$ have a common divisor
$2$, and $\lambda^{(1)}_1$ satisf\/ies
$\lambda^{(1)}_{1,\alpha+m^{(1)}_1/2}
=\lambda^{(1)}_{1,\alpha}+p^{(1)}_1/2$.
Thus from~(\ref{t:eq:aug18_4}),
the matrix $F_{\boldsymbol{\gamma}}$ (\ref{k:ftdef}) becomes
\begin{gather}\label{k:ft}
F_{\boldsymbol{\gamma}}  = F \cdot {\rm diag}
\big(\tfrac12,1,1,1\big)
= \begin{pmatrix}
 8 & 12 & -3 & -1 \\
 4 & 19 & -2 & -1 \\
 -3 & -6 & 10 & 2 \\
 -1 & -3 & 2 & 3
\end{pmatrix}.
\end{gather}
\end{example}

\begin{example}\label{k:ex:ga2}
Consider Example \ref{k:ex:aa2}.
Take the left hand side of (\ref{k:rc35}) as a
representative element
of the angle variable of the path (\ref{k:exp2}).
Then the splitting (\ref{k:spl}) gives
the reduced angle variable
$\omega = {}^t(-21,-20,1,-1)$ and
$\boldsymbol{\lambda}$ is the
quasi-periodic extension of the rigging in
\begin{equation*}
\begin{picture}(200,75)(-120,13)
\setlength{\unitlength}{0.5mm}

\multiput(0,0)(0,-10){3}{
\put(10,60){\line(1,0){30}}}
\multiput(0,0)(0,10){4}{
\put(10,10){\line(1,0){20}}}
\multiput(0,0)(10,0){3}{
\put(10,60){\line(0,-1){50}}}
\put(40,60){\line(0,-1){20}}

\put(42,52){$1$}
\put(42,42){$0$}
\put(32,32){$5$}
\put(32,22){$2$}
\put(32,12){$0$}

\multiput(-30,0)(0,0){1}{
\put(143,52){$0$}\put(113,42){$0$}
\put(100,60){\line(1,0){40}}
\put(100,50){\line(1,0){40}}
\put(100,40){\line(1,0){10}}
\put(100,60){\line(0,-1){20}}
\put(110,60){\line(0,-1){20}}
\multiput(0,0)(10,0){3}{
\put(120,60){\line(0,-1){10}}}}

\end{picture}
\end{equation*}
Thus $\boldsymbol{\lambda} =
\big(\lambda^{(1)}_1, \lambda^{(1)}_2,
\lambda^{(2)}_1, \lambda^{(2)}_2\big)$ reads
\begin{equation*}
\begin{split}
\lambda^{(1)}_1&=(\ldots, -4,-3, 0,1,4,5,8,9,\ldots),\quad
\lambda^{(1)}_2=(\ldots, -7,-5,-2, 0,2,5,7,9,12,\ldots),\\
\lambda^{(2)}_1&=(\ldots, -2, 0,2,4,\ldots),\quad\qquad\quad\quad \ \
\lambda^{(2)}_2=(\ldots, -1, 0,1,2,\ldots),
\end{split}
\end{equation*}
with $\lambda^{(a)}_{i,1}=0$.
The order of symmetry is trivial
in that $\forall \gamma^{(a)}_i = 1$.
In this case one has $F_{\boldsymbol{\gamma}}= F$
(\ref{t:eq:aug18_4}).
\end{example}

An important property of
the subset $X_{\boldsymbol{\gamma}}\subset X$ is that
it is still invariant under the actions of both ${\mathcal T}$ and ${\mathcal A}$.
Let $X_{\boldsymbol{\gamma}}/{\mathcal A}$
be the set of all ${\mathcal A}$-orbits.
According to
$X = \bigsqcup_{\boldsymbol{\gamma}}X_{\boldsymbol{\gamma}}$,
we have the disjoint union decomposition:
\begin{equation}\label{t:eq:aug21_4}
X/{\mathcal A} = \bigsqcup_{\boldsymbol{\gamma}}\,
(X_{\boldsymbol{\gamma}}/ {\mathcal A}),
\end{equation}
which induces the disjoint union decomposition of the
level set ${\mathcal P}(\mu)$ according to (\ref{k:cd15}).
Writing the pre-image of
$X_{\boldsymbol{\gamma}}/ {\mathcal A}$
as ${\mathcal P}_{\boldsymbol{\gamma}}(\mu)$,
one has ${\mathcal P}(\mu) = \sqcup_{\boldsymbol{\gamma}}
{\mathcal P}_{\boldsymbol{\gamma}}(\mu)$ and
the conjectural commutative diagram (\ref{k:cd15})
splits into
\begin{equation}\label{k:cd16}
\begin{CD}
{\mathcal P}_{\boldsymbol{\gamma}}(\mu)
@>{\Phi}>> X_{\boldsymbol{\gamma}}/{\mathcal A} \\
@V{\mathcal T}VV @VV{\mathcal T}V\\
{\mathcal P}_{\boldsymbol{\gamma}}(\mu)
@>{\Phi}>> X_{\boldsymbol{\gamma}}/{\mathcal A}
\end{CD}
\end{equation}
${\mathcal P}_{\boldsymbol{\gamma}}(\mu)$ is the subset of
the level set ${\mathcal P}(\mu)$ characterized by
the order of symmetry  $\boldsymbol{\gamma}$.
We are yet to decompose it or equivalently
$X_{\boldsymbol{\gamma}}/{\mathcal A}$
further into ${\mathcal T}$-orbits.
Note that one can think of an ${\mathcal A}$-orbit
${\mathcal A} \cdot x$ either as an element of
$X_{\boldsymbol{\gamma}} / {\mathcal A}$ or
as a subset of $X_{\boldsymbol{\gamma}}$.
Similarly a ${\mathcal T}$-orbit
${\mathcal T} \cdot ({\mathcal A} \cdot x)$
in $X_{\boldsymbol{\gamma}}/{\mathcal A}$
can either be regarded as an element of
$(X_{\boldsymbol{\gamma}} / {\mathcal A}) / {\mathcal T}$ or
as a subset of $X_{\boldsymbol{\gamma}}/ {\mathcal A}$.
We adopt the latter interpretation.
Then as we will see shortly in Proposition \ref{k:pr:t},
the ${\mathcal T}$-orbits
${\mathcal T} \cdot ({\mathcal A} \cdot x)$
of $X_{\boldsymbol{\gamma}} / {\mathcal A}$
can be described explicitly in terms of
${\mathbb Z}^g/F_{\boldsymbol{\gamma}}\,{\mathbb Z}^g$, the set of integer points
on the torus equipped with the following action of ${\mathcal T}$:
\begin{equation}\label{t:eq:aug11_4}
T_l^{(r)} ({\mathcal I}) = {\mathcal I} + h_l^{(r)}
\mod{F_{\boldsymbol{\gamma}}\,\mathbb{Z}^g}\qquad \text{for}
\ \ {\mathcal I} \in {\mathbb Z}^g/F_{\boldsymbol{\gamma}}\,{\mathbb Z}^g,
\end{equation}
where $h^{(r)}_l$ is specif\/ied in (\ref{k:ter}).
In what follows, we shall refer to
${\mathbb Z}^g/F_{\boldsymbol{\gamma}}\,{\mathbb Z}^g$ simply as torus.

Given any $x =(\omega, \boldsymbol{\lambda})
\in X_{\boldsymbol{\gamma}}$,
consider the map
\begin{equation}\label{k:chi}
\begin{split}
\chi: \; {\mathcal T} \cdot ({\mathcal A} \cdot x)
&\longrightarrow {\mathbb Z}^g/F_{\boldsymbol{\gamma}}\,{\mathbb Z}^g\\
t\cdot ({\mathcal A}\cdot x) & \longmapsto \;\;
t\cdot \omega\mod F_{\boldsymbol{\gamma}}\,{\mathbb Z}^g.
\end{split}
\end{equation}

\begin{proposition}[\cite{T}]\label{k:pr:t}
The map $\chi$ is well-defined, bijective and
the following diagram is commutative:
\begin{equation*}
\begin{CD}
{\mathcal T} \cdot ({\mathcal A} \cdot x)
@>{\chi}>> {\mathbb Z}^g / F_{\boldsymbol{\gamma}}\, {\mathbb Z}^g \\
@V{\mathcal T}VV @VV{\mathcal T}V\\
{\mathcal T} \cdot ({\mathcal A} \cdot x)
@>{\chi}>>{\mathbb Z}^g / F_{\boldsymbol{\gamma}}\, {\mathbb Z}^g
\end{CD}
\end{equation*}
\end{proposition}

Combining Proposition \ref{k:pr:t} with
Conjecture \ref{t:conj:aug21_5}
or its ref\/ined form (\ref{k:cd16}),
we obtain an explicit description of
each connected component (${\mathcal T}$-orbit)
as a torus.

\begin{conjecture}\label{t:conj:aug10_8}
For any path $p \in {\mathcal P}_{\boldsymbol{\gamma}}(\mu)$,
the map $\Phi_\chi:=\chi\circ \Phi$ gives a bijection between
the connected component $\Sigma (p)$ and the torus
$\mathbb{Z}^g / F_{\boldsymbol{\gamma}}\, \mathbb{Z}^g$
making the following diagram commutative:
\begin{equation}\label{k:cd3}
\begin{CD}
\Sigma(p)
@>{\Phi_\chi}>> {\mathbb Z}^g / F_{\boldsymbol{\gamma}}\,{\mathbb Z}^g \\
@V{\mathcal T}VV @VV{\mathcal T}V\\
\Sigma(p)
@>{\Phi_\chi}>>{\mathbb Z}^g / F_{\boldsymbol{\gamma}}\, {\mathbb Z}^g
\end{CD}
\end{equation}
\end{conjecture}

The actions of ${\mathcal T}$
in Proposition \ref{k:pr:t} and Conjecture \ref{t:conj:aug10_8}
are both transitive.
Conjecture \ref{t:conj:aug10_8} is a principal claim in this paper.
It says that the reduced angle variables live in
the torus ${\mathbb Z}^g / F_{\boldsymbol{\gamma}}\, {\mathbb Z}^g$,
where time evolutions of paths become straight motions.
See Conjecture \ref{t:conj:aug11_6} for an analogous claim
in a more general case than (\ref{k:pp}).

\begin{remark}\label{t:rem:aug18_1}
Due to the transitivity of ${\mathcal T}$-action,
the commutative diagram (\ref{k:cd3}) persists
even if $\Phi_\chi$ is redef\/ined as $\Phi_\chi+ c$ with
any constant vector $c \in {\mathbb Z}^g$, meaning that the choice
of the inverse image of  $0 \in {\mathbb Z}^g / F_{\boldsymbol{\gamma}}\, {\mathbb Z}^g$
is at one's disposal.
With this option in mind, an explicit way to construct the bijection
$\Phi_\chi$ is as follows.
Fix an arbitrary highest path $p_+$ in $\Sigma (p)$
and set $\Phi_\chi(p_+) = 0$.
For any $p' \in \Sigma (p)$
let $t$ be an element of ${\mathcal T}$ such that $p' = t \cdot p_+$.
One can always f\/ind such $t$ since ${\mathcal T}$ acts on
$\Sigma (p)$ transitively.
If $t$ is written as $t = \prod_{r,l} \big(T_l^{(r)}\big)^{d_l^{(r)}}$,
then we set
$\Phi_\chi(p') = \sum_{r,l} d_l^{(r)} h_l^{(r)}
\mod{F_{\boldsymbol{\gamma}} \,\mathbb{Z}^g}$.
\end{remark}

\begin{remark}\label{t:rem:aug18_2}
For $(a, i) \in \overline{H}$ we def\/ine
$\eta^{(a)}_i = l^{(a)}_{i+1} + 1 \, (i < g_a), \; = 1 \, (i=g_a)$.
Let $\overline{{\mathcal T}}$ be the free abelian group generated by all
$T_{\eta^{(a)}_i}^{(a)}$'s.
We let it act on the torus
$\mathbb{Z}^g / F_{\boldsymbol{\gamma}} \,\mathbb{Z}^g$
by (\ref{t:eq:aug11_4}).
This action is transitive because the $g$-dimensional lattice
$\bigoplus_{(ai)\in \overline{H}} {\mathbb Z} h^{(a)}_{\eta^{(a)}_i}
\subset {\mathbb Z}^g$
coincides with $\mathbb{Z}^g$ itself.
If Conjecture \ref{t:conj:aug10_8} is valid, then
as a subgroup of ${\mathcal T}$,
$\overline{{\mathcal T}}$ also acts on $\Sigma (p)$ transitively.
Hence the image of the bijection $\Phi_\chi$ in the previous remark
can be written as a linear combination of
$h^{(r)}_{\eta^{(r)}_i}$ with $(r,i) \in \overline{H}$ only.
\end{remark}

\begin{example}\label{k:ex:ivp}
As an application of the inverse scheme (\ref{k:cd3}),
we take the path $p$ (\ref{k:exp})
in Example \ref{k:ex:aa} and
illustrate the solution of the initial value problem to derive
\begin{equation}\label{k:aim}
\begin{split}
\big(T^{(1)}_2\big)^{1000}(p)=122211122113321113211331,\\
\big(T^{(1)}_3\big)^{1000}(p) = 213321112211112221331113,\\
\big(T^{(2)}_1\big)^{1000}(p) =331113211332112211111222,\\
\big(T^{(2)}_2\big)^{1000}(p) = 311132113211122211221133,\\
\big(T^{(2)}_3\big)^{1000}(p) = 113221132111221132213311,\\
\big(T^{(2)}_4\big)^{1000}(p) = 221113221112211321333111.
\end{split}
\end{equation}
In Example \ref{k:ex:ga}, we have seen
that the reduced angle variable of $p$ is given by
\begin{equation*}
\omega=
\begin{pmatrix} -1 \\ -2 \\ -3 \\ -3
\end{pmatrix}
\in {\mathbb Z}^4/F_{\boldsymbol{\gamma}}\,{\mathbb Z}^4,
\end{equation*}
where $F_{\boldsymbol{\gamma}}$ is specif\/ied in (\ref{k:ft}).
The velocity vector $h^{(a)}_l$ (\ref{k:ter})
of the reduced angle variable corresponding to the time evolution
$T^{(a)}_l$ is given by
\begin{gather*}
h^{(1)}_1  =  \begin{pmatrix} 1 \\ 1 \\ 0 \\ 0\end{pmatrix},\quad
h^{(1)}_2  =  2h^{(1)}_1,\quad
h^{(1)}_3  =  \begin{pmatrix} 3 \\ 2 \\ 0 \\ 0\end{pmatrix},\quad
h^{(2)}_1  =  \begin{pmatrix} 0 \\ 0 \\ 1 \\ 1\end{pmatrix},\\
h^{(2)}_2  =  \begin{pmatrix} 0 \\ 0 \\ 2 \\ 1\end{pmatrix},\quad
h^{(2)}_3  =  \begin{pmatrix} 0 \\ 0 \\ 3 \\ 1\end{pmatrix},\quad
h^{(2)}_4  =  \begin{pmatrix} 0 \\ 0 \\ 4 \\ 1\end{pmatrix}.
\end{gather*}
The higher time evolutions
$h^{(1)}_{l\ge 3}$ and $h^{(2)}_{l \ge 4}$
coincide with $h^{(1)}_3$ and $h^{(2)}_4$, respectively.
Moreover from $h^{(1)}_2  = 2h^{(1)}_1$,
one has $\big(T^{(1)}_2\big)^{1000}(p) = \big(T^{(1)}_1\big)^{2000}(p)$,
and the latter is easily obtained since~$T^{(1)}_1$
is just a cyclic shift (\ref{k:shf}), yielding
$\big(T^{(1)}_2\big)^{1000}(p) =122211122113321113211331$
in agreement with~(\ref{k:aim}).

Next we consider $\big(T^{(1)}_3\big)^{1000}(p)$.
The reduced angle variable for this path is
\begin{equation*}
\omega+1000h^{(1)}_3 =
\begin{pmatrix} 2999 \\ 1998 \\ -3 \\ -3\end{pmatrix}
\equiv
\begin{pmatrix} 7 \\ 6 \\ 1 \\ 1\end{pmatrix}
\; \text{ mod } F_{\boldsymbol{\gamma}}\,{\mathbb Z}^4 \;=
6h^{(1)}_1+h^{(2)}_1+
\begin{pmatrix} 1 \\ 0 \\ 0 \\ 0\end{pmatrix}
\; \text{ mod }F_{\boldsymbol{\gamma}}\,{\mathbb Z}^4.
\end{equation*}
This means $\big(T^{(1)}_3\big)^{1000}(p) =
\big(T^{(1)}_1\big)^6T^{(2)}_1(p')$, where $p'$ is the path
corresponding to ${}^t(1,0,0,0)$.
In view of (\ref{k:rc22}),
this reduced angle variable corresponds to
\begin{equation*}
\begin{picture}(200,70)(-220,17)
\setlength{\unitlength}{0.5mm}

\put(10,32)
{$112211132211321133113221 = p'.$}

\put(-8,32){$\overset{\phi^{-1}}{\longmapsto}$}

\multiput(-135,0)(0,0){1}{
\multiput(0,0)(0,-10){3}{
\put(10,60){\line(1,0){30}}}
\multiput(0,0)(0,10){4}{
\put(10,10){\line(1,0){20}}}
\multiput(0,0)(10,0){3}{
\put(10,60){\line(0,-1){50}}}
\put(40,60){\line(0,-1){20}}

\put(3,47){$4$}\put(3,22){$7$}
\put(42,52){$3$}
\put(42,42){$1$}
\put(32,32){$5$}
\put(32,22){$4$}
\put(32,12){$0$}

\multiput(-30,0)(0,0){1}{
\put(93,52){$2$}\put(93,42){$1$}
\put(143,52){$0$}\put(113,42){$0$}
\put(100,60){\line(1,0){40}}
\put(100,50){\line(1,0){40}}
\put(100,40){\line(1,0){10}}
\put(100,60){\line(0,-1){20}}
\put(110,60){\line(0,-1){20}}
\multiput(0,0)(10,0){3}{
\put(120,60){\line(0,-1){10}}}}
}

\end{picture}
\end{equation*}
From $T^{(2)}_1(p')=112211112221331113213321$,
one gets
$\big(T^{(1)}_1\big)^6T^{(2)}_1(p')$ by a cyclic shift.
The result yields
$\big(T^{(1)}_3\big)^{1000}(p)
=213321112211112221331113$ in agreement with (\ref{k:aim}).
Of course, the choice of the representative
$\text{mod } F_{\boldsymbol{\gamma}}\,{\mathbb Z}^4$
is not unique.

The procedure is parallel for $\big(T^{(2)}_l\big)^{1000}(p)$
with $l=1,\ldots, 4$.
The reduced angle variable of $\big(T^{(2)}_l\big)^{1000}(p)$ is
\begin{equation*}
\omega+ 1000h^{(2)}_l=
\begin{pmatrix}-1 \\ -2 \\ 1000l-3 \\ 997 \end{pmatrix}
\overset{\text{mod }F_{\boldsymbol{\gamma}}\,{\mathbb Z}^4}
{\equiv}
\begin{cases}
4h^{(1)}_1 + 8h^{(2)}_1 +
{}^t(0, 2, 0, 0), &  l=1, \\
10h^{(1)}_1 + h^{(2)}_1 +
{}^t(0, 2, 0, 0), &  l=2,\\
15h^{(1)}_1 + 10h^{(2)}_1 +
{}^t(1, 0, 1, 0), & l=3,\\
21h^{(1)}_1 + 3h^{(2)}_1 +
{}^t(1, 0, 1, 0), & l=4.
\end{cases}
\end{equation*}
{}From this,  one has
\begin{equation}\label{k:tochu}
\begin{split}
\big(T^{(2)}_1\big)^{1000}(p) &= \big(T^{(1)}_1\big)^4\big(T^{(2)}_1\big)^8(p''),
\qquad\quad
\big(T^{(2)}_2\big)^{1000}(p) = \big(T^{(1)}_1\big)^{10}T^{(2)}_1(p''),\\
\big(T^{(2)}_3\big)^{1000}(p) &= \big(T^{(1)}_1\big)^{15}\big(T^{(2)}_1\big)^{10}(p'''),
\qquad
\big(T^{(2)}_4\big)^{1000}(p) = \big(T^{(1)}_1\big)^{21}\big(T^{(2)}_1\big)^3(p''').
\end{split}
\end{equation}
Here the paths $p''$ and $p'''$ are those corresponding to
the reduced angle variables
${}^t(0,2,0,0)$ and ${}^t(1,0,1,0)$, respectively.
They are obtained as
\begin{equation*}
\begin{picture}(200,70)(-220,17)
\setlength{\unitlength}{0.5mm}

\put(10,32)
{$111222132111332113311122 = p'',$}

\put(-8,32){$\overset{\phi^{-1}}{\longmapsto}$}

\multiput(-135,0)(0,0){1}{
\multiput(0,0)(0,-10){3}{
\put(10,60){\line(1,0){30}}}
\multiput(0,0)(0,10){4}{
\put(10,10){\line(1,0){20}}}
\multiput(0,0)(10,0){3}{
\put(10,60){\line(0,-1){50}}}
\put(40,60){\line(0,-1){20}}

\put(3,47){$4$}\put(3,22){$7$}
\put(42,52){$2$}
\put(42,42){$0$}
\put(32,32){$7$}
\put(32,22){$6$}
\put(32,12){$2$}

\multiput(-30,0)(0,0){1}{
\put(93,52){$2$}\put(93,42){$1$}
\put(143,52){$0$}\put(113,42){$0$}
\put(100,60){\line(1,0){40}}
\put(100,50){\line(1,0){40}}
\put(100,40){\line(1,0){10}}
\put(100,60){\line(0,-1){20}}
\put(110,60){\line(0,-1){20}}
\multiput(0,0)(10,0){3}{
\put(120,60){\line(0,-1){10}}}}
}

\end{picture}
\end{equation*}
\begin{equation*}
\begin{picture}(200,70)(-220,17)
\setlength{\unitlength}{0.5mm}

\put(10,32)
{$112211132211122133113321 = p'''$}

\put(-8,32){$\overset{\phi^{-1}}{\longmapsto}$}

\multiput(-135,0)(0,0){1}{
\multiput(0,0)(0,-10){3}{
\put(10,60){\line(1,0){30}}}
\multiput(0,0)(0,10){4}{
\put(10,10){\line(1,0){20}}}
\multiput(0,0)(10,0){3}{
\put(10,60){\line(0,-1){50}}}
\put(40,60){\line(0,-1){20}}

\put(3,47){$4$}\put(3,22){$7$}
\put(42,52){$3$}
\put(42,42){$1$}
\put(32,32){$5$}
\put(32,22){$4$}
\put(32,12){$0$}

\multiput(-30,0)(0,0){1}{
\put(93,52){$2$}\put(93,42){$1$}
\put(143,52){$1$}\put(113,42){$0$}
\put(100,60){\line(1,0){40}}
\put(100,50){\line(1,0){40}}
\put(100,40){\line(1,0){10}}
\put(100,60){\line(0,-1){20}}
\put(110,60){\line(0,-1){20}}
\multiput(0,0)(10,0){3}{
\put(120,60){\line(0,-1){10}}}}
}

\end{picture}
\end{equation*}

In (\ref{k:tochu}), the evolutions under $T^{(2)}_1$ is given by
\begin{gather*}
T^{(2)}_1(p'') = 111222112211333111321132,
\\
\big(T^{(2)}_1\big)^8(p'') = 132113321122111112223311,\\
\big(T^{(2)}_1\big)^3(p''')=111221113221112211321333,
\\
\big(T^{(2)}_1\big)^{10}(p''')=132213311113221132111221.
\end{gather*}
Applying the cyclic shift $T^{(1)}_1$ to these results, one can check
that (\ref{k:tochu}) leads to (\ref{k:aim}).
\end{example}

\subsection{Bethe ansatz formula from size
and number of connected components}\label{k:sec:baf}

The results in the previous section provide a beautiful interpretation of
the character formula derived from the Bethe ansatz at $q=0$ \cite{KN}
in terms of the size and number of orbits
in the periodic $A^{(1)}_n$ SCA.
Let us f\/irst recall the character formula:
\begin{gather}
(x_1+\cdots + x_{n+1})^L  =
\sum_\mu \Omega(\mu)\, x_1^{L-|\mu^{(1)}|}
x_2^{|\mu^{(1)}|- |\mu^{(2)}|}\cdots
x_{n+1}^{|\mu^{(n-1)}|-|\mu^{(n)}|}, \label{k:kn1}\\
\Omega(\mu)  =
(\det F)\prod_{(a i) \in \overline{H}} \frac{1}{m_i^{(a)}}
\binom{p_i^{(a)} + m_i^{(a)} - 1}{m_i^{(a)} - 1}
\quad (\in {\mathbb Z}).\label{t:eq:aug19_1}
\end{gather}
In (\ref{k:kn1}),
the sum is taken over $n$-tuple of Young diagrams
$\mu=(\mu^{(1)},\ldots, \mu^{(n)})$
{\em without any} constraint.
All the quantities appearing in (\ref{t:eq:aug19_1}) are determined
from $\mu$ by the left diagram in~(\ref{k:rca0})
and (\ref{k:hdef})--(\ref{k:f}).
The fact $\Omega(\mu) \in {\mathbb Z}$ can be easily seen
by expanding $\det F$.
It is vital that the binomial coef\/f\/icient here is
the extended one:
\begin{equation*}
{a \choose b} = \frac{a(a-1)\cdots (a-b+1)}{b!}
\qquad (a \in {\mathbb Z},\  b \in {\mathbb Z}_{\ge 0}),
\end{equation*}
which can be negative outside the range $0 \le b \le a$.
Due to such negative contributions, the inf\/inite sum (\ref{k:kn1})
cancels out except leaving
the f\/initely many positive contributions exactly when
$L \ge |\mu^{(1)}| \ge \cdots \ge |\mu^{(n)}|$.
For example when $L=6$, $n=2$, coef\/f\/icients of monomials in
the expansion (\ref{k:kn1}) of $(x_1+x_2+x_3)^6$
are calculated as
\begin{gather*}
60\, x_1^3x_2^2x_1  : \  \big(|\mu^{(1)}|, |\mu^{(2)}|\big)=(3,1),\\
\phantom{60\, x_1^3x_2^2x_1  :}{} \ \ \Omega((3),(1)) +
\Omega((21),(1))  +
\Omega((111),(1)) = 6+36+18 = 60,
\\
0\,x_1^5x_2^{-2}x_3^3   : \  \big(|\mu^{(1)}|, |\mu^{(2)}|\big)=(1,3),\\
\phantom{0\,x_1^5x_2^{-2}x_3^3   :}{} \  \ \Omega((1),(3)) +
\Omega((1),(21)) +
\Omega((1),(111)) = 6+(-18)+12=0.
\end{gather*}
In this way $\Omega(\mu)$ gives a decomposition of the multinomial
coef\/f\/icients according to $\mu$.
It is known \cite[Lemma 3.7]{KN} that $\Omega(\mu) \ge 1$
if $\mu$ is a conf\/iguration,  namely under the condition
$p^{(a)}_i \ge 0$ for any $(a, i)\in \overline{H}$.

The formula (\ref{t:eq:aug19_1}) originates in the
Bethe ansatz for the integrable vertex model
associated with
$U_q(A^{(1)}_n)$ under the periodic boundary condition.
In this context,  $\mu$ specif\/ies the string content.
Namely, $\mu^{(a)}$ signif\/ies that
there are $m^{(a)}_i$ strings of color $a$ and length $l^{(a)}_i$.
Under such a string hypothesis, the Bethe equation becomes
a linear congruence equation at $q=0$ (\ref{eq:sce}), and counting
its of\/f-diagonal solutions yields $\Omega(\mu)$.
In this sense, the identity (\ref{k:kn1}) implies a~formal
completeness of the Bethe ansatz and string hypothesis at $q=0$.
For more details, see Section~\ref{k:sec:rba}, especially Theorem~\ref{k:th:kn}.
A parallel story is known also at $q=1$ \cite{KKR}
as mentioned under Theorem~\ref{k:th:rc}.

Back to our $A^{(1)}_n$ SCA,
it is nothing but the integrable vertex model at $q=0$,
where $U_q\big(A^{(1)}_n\big)$ modules and row transfer matrices
are ef\/fectively replaced by the crystals and time evolution operators,
respectively.
In view of this, it is natural to link
the Bethe ansatz formula $\Omega(\mu)$
with the notions introduced in the previous subsection
like level set, torus, connected components (${\mathcal T}$-orbits) and so forth.
This will be done in this subsection,  providing
a conceptual explanation of the earlier observations \cite{KT1,KT2}.
Our main result is stated as
\begin{theorem}\label{k:th:bmain}
Assume the condition \eqref{k:pp} and Conjecture~{\rm \ref{t:conj:aug10_8}}.
Then the Bethe ansatz formula $\Omega(\mu)$ \eqref{t:eq:aug19_1}
counts the number of paths in the level set ${\mathcal P}(\mu)$
as follows:
\begin{equation}\label{k:kore}
|{\mathcal P}(\mu)| = \Omega(\mu)
= \sum_{\boldsymbol{\gamma}}
\underbrace{\det F_{\boldsymbol{\gamma}}}_{\text{size of a ${\mathcal T}$-orbit}}
\quad
\underbrace{
\prod_{(ai) \in \overline{H}}
\frac{\big|\Lambda_{\gamma^{(a)}_i}\big(m_i^{(a)},p_i^{(a)}\big)\big|}
{m_i^{(a)}/\gamma_i^{(a)}}}_{\text{number of ${\mathcal T}$-orbits}}.
\end{equation}
Here
$\boldsymbol{\gamma}
= \big(\gamma^{(a)}_i\big)_{(ai) \in \overline{H}}$
and the sum extends over all the orders of symmetry, i.e.,
each $\gamma^{(a)}_i$ runs over all the common divisors of
$m^{(a)}_i$ and $p^{(a)}_i$.
\end{theorem}

The result (\ref{k:kore}) uncovers the SCA meaning of the
Bethe ansatz formula (\ref{t:eq:aug19_1}).
It consists of the contributions from sectors
specif\/ied by the order of symmetry
$\boldsymbol{\gamma}$.
Each sector is an assembly of identical tori (connected components),
therefore its contribution is factorized into
its size and number (multiplicity).

For the proof we prepare a few facts.
Note that $s^{(a)}_i$ (\ref{t:eq:aug21_1}) sending
$\boldsymbol{\lambda}$ to $\boldsymbol{\lambda}'$ also
def\/ines an action of ${\mathcal A}$ on
$X^2_{\boldsymbol{\gamma}}$ part alone.

\begin{lemma}\label{t:lem:aug21_3}
The number of ${\mathcal T}$-orbits in
$X_{\boldsymbol{\gamma}}/{\mathcal A}$ is
$|X^2_{\boldsymbol{\gamma}}/{\mathcal A}|$,  i.e.,
the number of ${\mathcal A}$-orbits in $X^2_{\boldsymbol{\gamma}}$.
\end{lemma}
\begin{proof}
Let $y = {\mathcal A} \cdot x$
and $y' = {\mathcal A} \cdot x'$ be any two elements of
$X_{\boldsymbol{\gamma}}/ {\mathcal A}$.
They belong to a common ${\mathcal T}$-orbit
if and only if there exists $t \in {\mathcal T}$ such that $t \cdot y = y'$.
It is equivalent to saying that there exist
$t \in {\mathcal T}$ and $a \in {\mathcal A}$
such that $t \cdot (a \cdot x) = x'$.
Write $x,x'$ as $x = (\omega, \boldsymbol{\lambda})$
and $x' = (\omega', \boldsymbol{\lambda}')$.

Suppose $\boldsymbol{\lambda}$ and $\boldsymbol{\lambda}'$
belong to a common ${\mathcal A}$-orbit in
$X^2_{\boldsymbol{\gamma}}$.
Then there exits $a \in {\mathcal A}$ such that
$a \cdot \boldsymbol{\lambda} = \boldsymbol{\lambda}'$.
Hence we have $a \cdot x = (\tilde{\omega}, \boldsymbol{\lambda}')$
with some $\tilde{\omega} \in X^1$.
Since ${\mathcal T}$ acts on the $X^1$ part transitively
and leaves the $X^2_{\boldsymbol{\gamma}}$ part untouched,
there exits $t \in {\mathcal T}$ such that $t \cdot (a \cdot x) = x'$.

Suppose $\boldsymbol{\lambda}$ and $\boldsymbol{\lambda}'$
belong to dif\/ferent ${\mathcal A}$-orbits.
Then $a \cdot x$ and $x'$ have dif\/ferent
$X^2_{\boldsymbol{\gamma}}$ parts
for any $a \in {\mathcal A}$.
Hence for any $t \in {\mathcal T}$ and
$a \in {\mathcal A}$, we have $t \cdot (a \cdot x) \ne x'$.

Thus $y$ and $y'$ belong to a common ${\mathcal T}$-orbit if and only if
$\boldsymbol{\lambda}$ and $\boldsymbol{\lambda}'$
belong to a common ${\mathcal A}$-orbit.
The proof is completed.
\end{proof}

\begin{lemma}\label{k:le:xa}
The cardinality of the set
$X^2_{\boldsymbol{\gamma}}/{\mathcal A}$ is given by
\begin{equation*}
|X^2_{\boldsymbol{\gamma}}/{\mathcal A}|
= \prod_{(ai) \in \overline{H}}
\frac{\big|\Lambda_{ \gamma^{(a)}_i}\big(m^{(a)}_i, p^{(a)}_i\big)\big|}
{m^{(a)}_i/\gamma^{(a)}_i}.
\end{equation*}
\end{lemma}

\begin{proof}
It is suf\/f\/icient to show that each factor in the right hand side
gives the number of ${\mathcal A}$-orbits for each $(a,i)$ block.
So we omit all the
indices below and regard ${\mathcal A}$ as a free abelian
group generated by
a single element~$s$.

Let $\lambda = (\lambda_\alpha)_{\alpha \in {\mathbb Z}}$
be any element of $\Lambda_\gamma(m,p)$.
Then due to (\ref{t:eq:aug21_1}) we have $s^n \cdot \lambda =
(\lambda_{\alpha + n} - \lambda_{1 + n})_{\alpha \in {\mathbb Z}}$.
By (\ref{t:eq:may12_1}) and
(\ref{t:eq:may1_1}) this implies that $s^{m/\gamma} \cdot \lambda =
\lambda$ and
there is no $0 < k < m/\gamma$ such that $s^k \cdot \lambda = \lambda$.
Hence the number of ${\mathcal A}$-orbits in
$\Lambda_\gamma(m,p)$ is
$|\Lambda_\gamma(m,p)|/(m/\gamma)$.
\end{proof}

\begin{lemma}\label{k:le:tos}
For any path
$p \in {\mathcal P}_{\boldsymbol{\gamma}}(\mu)$,
the size of each ${\mathcal T}$-orbit $\Sigma(p)$ is
$|\Sigma(p)| = \det F_{\boldsymbol{\gamma}}$.
\end{lemma}

\begin{proof}
This is due to Conjecture \ref{t:conj:aug10_8}
and $\det F_{\boldsymbol{\gamma}}>0$ under the assumption (\ref{k:pp}).
See the remark under (\ref{k:ftdef}).
\end{proof}

Now we are ready to give

\begin{proof}[Proof of Theorem \ref{k:th:bmain}.]
\begin{align*}
|{\mathcal P}(\mu)| &\overset{(\ref{k:cd15}), \, (\ref{t:eq:aug21_4})}{=}
\sum_{\boldsymbol{\gamma}}
|X_{\boldsymbol{\gamma}}/{\mathcal A}|
=\sum_{\boldsymbol{\gamma}}
(\text{size of ${\mathcal T}$-orbits})\times
(\text{number of ${\mathcal T}$-orbits})\\
&\overset{\text{Lemmas} \, \ref{t:lem:aug21_3},\, \ref{k:le:tos}}{=}
\sum_{\boldsymbol{\gamma}}(\det F_{\boldsymbol{\gamma}})\,
|X^2_{\boldsymbol{\gamma}}/{\mathcal A}|\\
&\overset{(\ref{k:ftdef}), \, \text{Lemma}\, \ref{k:le:xa}}{=}
\sum_{\boldsymbol{\gamma}}
\frac{\det F}{\prod_{(ai) \in \overline{H}} \gamma^{(a)}_i}\,
\prod_{(ai) \in \overline{H}}
\frac{\big|\Lambda_{\gamma^{(a)}_i}\big(m_i^{(a)},p_i^{(a)}\big)\big|}
{m^{(a)}_i/\gamma^{(a)}_i}\\
&= \frac{\det F}{\prod_{(ai) \in \overline{H}}m^{(a)}_i}
\sum_{\boldsymbol{\gamma}}
\prod_{(ai) \in \overline{H}}
|\Lambda_{\gamma^{(a)}_i}(m_i^{(a)},p_i^{(a)})|\\
&= \frac{\det F}{\prod_{(ai) \in \overline{H}}m^{(a)}_i}
\prod_{(ai) \in \overline{H}}\Bigg(
\sum_{\gamma^{(a)}_i}
\big|\Lambda_{\gamma^{(a)}_i}\big(m_i^{(a)},p_i^{(a)}\big)\big|\Bigg)\\
&\overset{(\ref{t:eq:aug19_2})}{=}
(\det F)\prod_{(a i) \in \overline{H}} \frac{1}{m_i^{(a)}}
\binom{p_i^{(a)} + m_i^{(a)} - 1}{m_i^{(a)} - 1}
= \Omega(\mu). \tag*{\qed}
\end{align*}
  \renewcommand{\qed}{}
\end{proof}

One can give an explicit formula for
$|\Lambda_{\gamma}(m, p)|$
by solving (\ref{t:eq:aug19_2}) by the M{\"o}bius inversion:
\begin{equation}\label{k:lc}
|\Lambda_{\gamma}(m, p)| =
\sum_\beta \mu\left(\frac{\beta}{\gamma}\right)
{\frac{p+m}{\beta} -1 \choose \frac{m}{\beta} -1},
\end{equation}
where $\beta$ runs over
all the common divisors of $m$ and $p$ that is a multiple of $\gamma$.
Here $\mu$ is the
M{\"o}bius function in number theory  \cite{St} def\/ined by
$\mu(1)=1, \mu (k)=0$ if $k$
is divisible by the square of an integer greater than one,
and $\mu(k) = (-1)^j$ if $k$ is the product of $j$ distinct primes.
(This $\mu$ should not be confused with
the $n$-tuple of Young diagrams.)

\begin{example}\label{k:ex:deco}
Let us consider the decomposition of the level set
${\mathcal P}(\mu)$
with $\mu = ((33222),(41))$, which is the same as
Examples \ref{k:ex:aa} and \ref{k:ex:aa2}.
{}From the tables therein and
$\det F= 4656$, the formula (\ref{t:eq:aug19_1}) reads
\begin{gather*}
\Omega(\mu) =\det F\, \frac{1}{m^{(1)}_1}
\binom{p^{(1)}_1 + m^{(1)}_1-1}{m^{(1)}_1-1}
\frac{1}{m^{(1)}_2}
\binom{p^{(1)}_2 + m^{(1)}_2-1}{m^{(1)}_2-1}\\
\phantom{\Omega(\mu)}{}
=4656\cdot\frac{1}{2}\binom{5}{1}
\cdot\frac{1}{3}\binom{9}{2}=139680.
\end{gather*}
On the other hand, (\ref{k:lc}) gives
\begin{gather*}
\big|\Lambda_1\big(p^{(1)}_1,m^{(1)}_1\big)\big|
 = \mu\left(\frac{1}{1}\right)
{\frac{4+2}{1} -1 \choose \frac{2}{1} -1}+
\mu\left(\frac{2}{1}\right)
{\frac{4+2}{2} -1 \choose \frac{2}{2} -1}
= 5-1 = 4, \\
\big|\Lambda_2\big(p^{(1)}_1,m^{(1)}_1\big)\big|  =
\mu \left(\frac{2}{2}\right)
{\frac{4+2}{2} -1 \choose \frac{2}{2} -1}= 1,\\
\big|\Lambda_1\big(p^{(1)}_2,m^{(1)}_2\big)\big|  =
\mu \left(\frac{1}{1}\right)
{\frac{7+3}{1} -1 \choose \frac{3}{1} -1}= 36,\\
\big|\Lambda_1\big(p^{(2)}_1, m^{(2)}_1\big)\big|  =
\big|\Lambda_1\big(p^{(2)}_2, m^{(2)}_2\big)\big| = 1.
\end{gather*}
The formula (\ref{k:kore}) consists of
the two terms corresponding to
the order of symmetry
$\boldsymbol{\gamma}
=\big(\gamma^{(1)}_1,\gamma^{(1)}_2,
\gamma^{(2)}_1, \gamma^{(2)}_2\big)=(1,1,1,1)$ and $(2,1,1,1)$.
Let us write them as
$\boldsymbol{\gamma}_1$ and
$\boldsymbol{\gamma}_2$, respectively.
$F_{\boldsymbol{\gamma}_1}=F$ is given by
(\ref{t:eq:aug18_4}) and
$F_{\boldsymbol{\gamma}_2}$ by (\ref{k:ft}).
Their determinants are
$\det F_{\boldsymbol{\gamma}_1} = 4656$ and
$\det F_{\boldsymbol{\gamma}_2} = 2328$.
Thus (\ref{k:kore}) reads
\begin{equation*}
|{\mathcal P}(\mu)| =
\Omega (\mu) = 139680
=(\det F_{\boldsymbol{\gamma}_1})\cdot 24
+(\det F_{\boldsymbol{\gamma}_2}) \cdot 12.
\end{equation*}
Accordingly the level set ${\mathcal P}(\mu)$
is decomposed into 36 tori as
\begin{equation*}
{\mathcal P}(\mu)
= 24\big(\mathbb{Z}^4 /
F_{\boldsymbol{\gamma}_1}\mathbb{Z}^4\big)
\sqcup 12\big(\mathbb{Z}^4 /
F_{\boldsymbol{\gamma}_2}\mathbb{Z}^4\big).
\end{equation*}
The paths (\ref{k:exp}) and (\ref{k:exp2})
in Examples \ref{k:ex:aa} and \ref{k:ex:aa2} belong to one of
the tori
$\mathbb{Z}^4 / F_{\boldsymbol{\gamma}_2} \mathbb{Z}^4$
and
$\mathbb{Z}^4 / F_{\boldsymbol{\gamma}_1} \mathbb{Z}^4$,
respectively.
\end{example}

\subsection{Dynamical period}\label{k:sec:dp}
Given a path
$p \in {\mathcal P}_{\boldsymbol{\gamma}}(\mu)$, the
{\em smallest} positive integer satisfying
$(T_l^{(r)})^{\mathcal N} (p) = p$ is called
the {\em dynamical period} of $p$
under the time evolution $T_l^{(r)}$.
Here we derive an explicit formula of the dynamical period
as a simple corollary of Conjecture \ref{t:conj:aug10_8}.
It takes the precise account of the symmetry specif\/ied by
$\boldsymbol{\gamma}$ and ref\/ines the earlier
conjectures in \cite{KT1, KT2} for $n\ge 2$
which was obtained from the Bethe eigenvalues at $q=0$.

For nonzero rational numbers $r_1 , \ldots ,r_s$,
we def\/ine their least common multiple by
\begin{equation*}
{\rm LCM}(r_1 , \ldots ,r_s)= \min\bigl(|{\mathbb Z}\cap r_1{\mathbb Z}
\cap \cdots \cap r_s {\mathbb Z}|\setminus \{0\}\bigr).
\end{equation*}
Given $(r, l)$ with $1 \leq r \leq n$ and $l \ge 1$,
def\/ine $F_{\boldsymbol{\gamma}} [bj]$ (resp. $F[bj]$)
to be the $g\times g$ matrix obtained from
$F_{\boldsymbol{\gamma}}$ (\ref{k:ftdef}) (resp. $F$ (\ref{k:f}))
by replacing its $(bj)$th column by $h_l^{(r)}$ (\ref{k:ter}).
Set
\begin{equation}\label{t:eq:aug18_3}
{\mathcal N}_l^{(r)} = {\rm LCM}\!
\left(\frac{\det F_{\boldsymbol{\gamma}}}
{\det F_{\boldsymbol{\gamma}} [bj]}
\right)_{(bj) \in \overline{H}}
= {\rm LCM}\!
\left(\frac{\det{F}}{\gamma^{(b)}_j\det{F} [bj]}
\right)_{(bj) \in \overline{H}},
\end{equation}
where the $\text{LCM}$ should be taken over only those
$(bj)$ satisfying $\det F[bj] \ne 0$.

\begin{theorem}\label{k:th:dp}
Under Conjecture {\rm \ref{t:conj:aug10_8}},
the dynamical period under $T_l^{(r)}$
is equal to ${\mathcal N}_l^{(r)}$ \eqref{t:eq:aug18_3}
for all the paths in
${\mathcal P}_{\boldsymbol{\gamma}}(\mu)$.
\end{theorem}

\begin{proof}
Under Conjecture \ref{t:conj:aug10_8},
the relation $\big(T_l^{(r)}\big)^{\mathcal N} (p) = p$
is equivalent to
\begin{equation}\label{nhf}
{\mathcal N}h_l^{(r)}
\equiv 0 \mod{F_{\boldsymbol{\gamma}}\, \mathbb{Z}^g}
\end{equation}
for any path $p \in {\mathcal P}_{\boldsymbol{\gamma}}(\mu)$.
In other words there exists $z \in \mathbb{Z}^g$
such that
${\mathcal N}h_l^{(r)}  = F_{\boldsymbol{\gamma}}\,z$.
By Cramer's formula, the solution
$z = (z_j^{(b)})_{(bj) \in \overline{H}}$
of this equation is given by
$z_j^{(b)} = {\mathcal N}
\frac{\det F_{\boldsymbol{\gamma}} [bj]}
{\det F_{\boldsymbol{\gamma}}}$.
The quantity ${\mathcal N}_l^{(r)}$ (\ref{t:eq:aug18_3})
is the smallest positive integer ${\mathcal N}$ that
matches the condition $z \in \mathbb{Z}^g$ hence (\ref{nhf}).
\end{proof}

\begin{remark}\label{k:re:com}
To the dynamical period under a composite time evolution
$T=\prod_{r,l}\big(T^{(r)}_l\big)^{d^{(r)}_l}$,
the same formula (\ref{t:eq:aug18_3}) applies
by replacing the role of $h^{(r)}_l$ therein with
$\sum_{r,l}d^{(r)}_lh^{(r)}_l$.
\end{remark}

\begin{remark}
The formula (\ref{t:eq:aug18_3}) can be simplif\/ied when
$n=1$ (hence $r=1$) and the order of symmetry
is trivial, i.e., $\boldsymbol{\gamma}=(\gamma^{(1)}_i)$ with
$\forall\, \gamma^{(1)}_i = 1$
as explained in \cite[equation~(4.24)]{KTT}.
In particular if $l=\infty$ furthermore, one gets
($p^{(1)}_0 = l^{(1)}_0=0$)
\begin{equation*}
{\mathcal N}_\infty^{(1)} = {\rm LCM}\!
\left(\frac{p^{(1)}_jp^{(1)}_{j-1}}
{(l^{(1)}_{j}-l^{(1)}_{j-1})p^{(1)}_{g_1}}
\right)_{1 \le j \le g_1},
\end{equation*}
which reproduces the result obtained in \cite{YYT}.
\end{remark}

\begin{example}\label{k:ex:dp}
Let us take the path $p$ (\ref{k:exp}) and
compute its dynamical period ${\mathcal N}^{(r)}_l$
under the time evolution $T^{(r)}_l$.
The matrix $F_{\boldsymbol{\gamma}}$
is given in (\ref{k:ft}) and
$\det F_{\boldsymbol{\gamma}}= 2328$.
The relevant vectors~$h^{(a)}_l$ are listed in
Example \ref{k:ex:ivp}.
Consider $T^{(1)}_3$ for instance.
Then, $h^{(1)}_3={}^t(3,2,0,0)$ and
\begin{alignat*}{3}
& F_{\boldsymbol{\gamma}}[11] =
\begin{pmatrix}
3 & 12 & -3 & -1 \\
2 & 19 & -2 & -1 \\
0 & -6 & 10 & 2 \\
0 & -3 & 2 & 3
\end{pmatrix},\qquad &&
F_{\boldsymbol{\gamma}}[12] =
 \begin{pmatrix}
 8 & 3 & -3 & -1 \\
 4 & 2 & -2 & -1 \\
 -3 & 0 & 10 & 2 \\
 -1 & 0 & 2 & 3
\end{pmatrix}, &\\
& F_{\boldsymbol{\gamma}}[21]= \begin{pmatrix}
 8 & 12 &  3 & -1 \\
 4 & 19 & 2 & -1 \\
 -3 & -6 & 0 & 2 \\
 -1 & -3 & 0 & 3
\end{pmatrix},\qquad &&
F_{\boldsymbol{\gamma}}[22] = \begin{pmatrix}
 8 & 12 & -3 & 3 \\
 4 & 19 & -2 & 2 \\
 -3 & -6 & 10 & 0 \\
 -1 & -3 & 2 & 0
\end{pmatrix} &
\end{alignat*}
with
$\det F_{\boldsymbol{\gamma}}[11] = 840$,
$\det F_{\boldsymbol{\gamma}}[12] = 108$,
$\det F_{\boldsymbol{\gamma}}[22]= 276$ and
$\det F_{\boldsymbol{\gamma}}[22] = 204$.
Thus (\ref{t:eq:aug18_3}) is calculated as
\begin{gather*}
{\mathcal N}^{(1)}_3  = {\rm LCM}\left(
\frac{\det{F_{\boldsymbol{\gamma}}}}
{\det{F_{\boldsymbol{\gamma}}}[11]},
\frac{\det{F_{\boldsymbol{\gamma}}}}
{\det{F_{\boldsymbol{\gamma}}}[12]},
\frac{\det{F_{\boldsymbol{\gamma}}}}
{\det{F_{\boldsymbol{\gamma}}}[21]},
\frac{\det{F_{\boldsymbol{\gamma}}}}
{\det{F_{\boldsymbol{\gamma}}}[22]}
\right)\\
\phantom{{\mathcal N}^{(1)}_3}{} ={\rm LCM}\left(
\frac{2328}{840},
\frac{2328}{108},
\frac{2328}{276},
\frac{2328}{204}
\right) = {\rm LCM}\left(\frac{97}{35},\frac{194}{9},
\frac{194}{23},\frac{194}{17}\right) = 194.
\end{gather*}

Similarly, we get
\begin{gather*}
{\mathcal N}^{(1)}_1 =
{\rm LCM}\left(12, 24, 24, 24\right) = 24,\\
{\mathcal N}^{(1)}_2 =
{\rm LCM}\left(6, 12, 12, 12\right) = 12,\\
{\mathcal N}^{(2)}_1 =
{\rm LCM}\left(\frac{582}{23},\frac{1164}{17},
\frac{1164}{65},\frac{1164}{377}\right) = 1164,\\
{\mathcal N}^{(2)}_2 =
{\rm LCM}\left(\frac{388}{29},\frac{776}{13},
\frac{776}{141},\frac{776}{197}\right) = 776,\\
{\mathcal N}^{(2)}_3 =
{\rm LCM}\left(\frac{291}{32},\frac{582}{11},
\frac{582}{179},\frac{582}{107}\right) = 582,\\
{\mathcal N}^{(2)}_4 =
{\rm LCM}\left(\frac{1164}{169},\frac{2328}{49},
\frac{2328}{1009},\frac{2328}{265}\right) = 2328.
\end{gather*}
Despite the nontrivial order of symmetry
$\boldsymbol{\gamma}=(2,1,1,1)$, these f\/inal results coincide
with the trivial case $\boldsymbol{\gamma}=(1,1,1,1)$
treated in Example~\ref{k:ex:int}.
We have checked that these values agree with the actual
dynamical periods of $p$ by computer.
\end{example}

\subsection[Relation to Bethe ansatz at $q=0$]{Relation to Bethe ansatz at $\boldsymbol{q=0}$}\label{k:sec:rba}

Let us quickly recall the relevant results from the
Bethe ansatz at $q=0$.
For the precise def\/initions and statements,
we refer to \cite{KN}.
Consider the integrable $U_q\big(A^{(1)}_n\big)$ vertex model
on a periodic chain of length $L$.
If the quantum space is
the $L$-fold tensor product of the vector representation,
the Bethe equation takes the form:
\begin{equation}\label{eq:be}
\left(\frac{\sin{\pi \big(u^{(a)}_i + \sqrt{-1}\hbar\delta_{a1}\big)}}
{\sin{\pi\big(u^{(a)}_i - \sqrt{-1}\hbar\delta_{a1}\big)}}\right)^{L}
= - \prod_{b=1}^n\prod_{j=1}^{|\mu^{(b)}|}
\frac{\sin\pi\big(u^{(a)}_i - u^{(b)}_j + \sqrt{-1}\hbar C_{ab} \big)}
{\sin\pi\big(u^{(a)}_i - u^{(b)}_j - \sqrt{-1}\hbar C_{ab}\big)}
\end{equation}
for $1 \le a \le n$ and $1 \le i \le \big|\mu^{(a)}\big|$.
Here $L\ge \big|\mu^{(1)}\big| \ge \cdots  \ge \big|\mu^{(n)}\big| $
specif\/ies a sector (\ref{k:an}) preserved by row transfer matrices.
The parameter $\hbar$ is related to $q$ by $q=e^{-2\pi\hbar}$.

Fix an $n$-tuple of Young diagrams, i.e., the string content
$\mu=\big(\mu^{(1)},\ldots, \mu^{(n)}\big)$ as in (\ref{k:rca0}).
We keep the notations (\ref{k:hdef})--(\ref{k:gdef}).
By string solutions we mean the ones
in which the unknowns
$\big\{u^{(a)}_i\mid 1 \le a \le n, 1\le i  \le |\mu^{(a)}|\big\}$
are arranged as
\begin{equation*}
\bigcup_{(ai\alpha) \in H}
\bigcup_{u^{(a)}_{i,\alpha} \in {\mathbb R}}
\big\{ u^{(a)}_{i,\alpha} + \sqrt{-1}(i+1-2k)\hbar
+ \epsilon^{(a)}_{i\alpha k}
\mid 1 \le k \le i \big\},
\end{equation*}
where $\epsilon^{(a)}_{i\alpha k}$ stands for a small deviation.
$u^{(a)}_{i,\alpha}$ is the string center of the $\alpha$th
string of color~$a$ and length~$i$.
For a generic string solution,
the Bethe equation is linearized at $q=0$ into a~logarithmic form
called the {\em string center equation}:
\begin{equation}\label{eq:sce}
\sum_{(bj\beta)\in H}
A_{ai\alpha, bj\beta} \,u^{(b)}_{j,\beta} \equiv
\frac{1}{2}\big(p^{(a)}_i + m^{(a)}_i + 1\big) \quad \mathrm{mod}\ {\mathbb Z}
\end{equation}
for $(ai\alpha)\in H$.
Here the $G\times G$ coef\/f\/icient matrix
$A=(A_{ai\alpha, bj\beta})_{(ai\alpha), (bj\beta) \in H}$
is specif\/ied as
\begin{equation}\label{eq:a}
A_{ai\alpha, bj\beta}=\delta_{ab}\delta_{ij}\delta_{\alpha\beta}
\big(p^{(a)}_i+m^{(a)}_i\big) + C_{ab}
\min\big(l^{(a)}_i, l^{(b)}_j\big)-\delta_{ab}\delta_{ij}.
\end{equation}
It is known \cite{KN} that $A$ is positive def\/inite
if $\mu$ is a conf\/iguration, namely if~(\ref{k:pg}) is satisf\/ied.

There are a number of conditions which
the solutions of the string center equation~(\ref{eq:sce})
are to satisfy or to be identif\/ied thereunder.
First, the Bethe vector depends on $u^{(a)}_i$ only via
$e^{2\pi\sqrt{-1}u^{(a)}_i}$.
Therefore the string center
should be understood as
$u^{(a)}_{i,\alpha} \in \mathbb{R}/{\mathbb Z}$ rather than $\mathbb{R}$.
Second, the original Bethe equation (\ref{eq:be}) is symmetric
with respect to $u^{(a)}_i$ for $i=1, \ldots, \big|\mu^{(a)}\big|$, but
their permutation does not lead to a new Bethe vector.
Consequently, we should regard the string centers
of the $(a,i)$ block as
\begin{equation*}
\big(u^{(a)}_{i,1}, u^{(a)}_{i,2}, \ldots, u^{(a)}_{i, m^{(a)}_i}\big)
\in \left(\mathbb{R}/{\mathbb Z}\right)^{m^{(a)}_i}\!/{\frak S}_{m^{(a)}_i}.
\end{equation*}
Last, we prohibit the collision of string centers
$u^{(a)}_{i,\alpha}=u^{(a)}_{i,\beta}$ for
$1 \le \alpha \neq \beta \le m^{(a)}_i$
for any $(a,i) \in \overline{H}$.
This is a remnant of the well-known
constraint on the Bethe roots so that
the associated Bethe vector does not vanish.
To summarize, we consider {\em off-diagonal solutions}
${\bf u} = \big(u^{(a)}_{i,\alpha}\big)_{(ai\alpha)\in H}$
to the string center equation (\ref{eq:sce}) that live in
\begin{equation*}
(u^{(a)}_{i,1}, u^{(a)}_{i,2}, \ldots, u^{(a)}_{i, m^{(a)}_i}) \in
\bigl(\left(\mathbb{R}/{\mathbb Z}\right)^{m^{(a)}_i} -
\Delta_{m^{(a)}_i} \bigr) /{\frak S}_{m^{(a)}_i}
\quad\hbox{ for each }(a,i) \in \overline{H},
\end{equation*}
where $\Delta_m  = \{(v_1, \ldots, v_m)
\in (\mathbb{R}/{\mathbb Z})^m \mid v_\alpha = v_\beta\ \text{for some}\
1 \le \alpha \neq
\beta \le m\}$.
For simplicity we will often say
{\em $(q=0)$ Bethe roots} to mean
the of\/f-diagonal solutions to the string center equation.
Let ${\mathcal U}(\mu)$ be the set of the Bethe roots
having the string content $\mu$.
The following result is derived by counting the
Bethe roots of (\ref{eq:sce}) by the
M\"obius inversion formula.

\begin{theorem}[\protect{\cite[Theorem 3.2]{KN}}]\label{k:th:kn}
Under the condition \eqref{k:pg},
the formula \eqref{t:eq:aug19_1} gives the
number of Bethe roots, namely
$\Omega(\mu) = |{\mathcal U}(\mu)|$ is valid.
\end{theorem}

Given an extended rigged conf\/iguration
${\bf r}=\big(r^{(a)}_{i,\alpha}\big)_{(ai) \in \overline{H},\alpha \in {\mathbb Z}}
\in \tilde{\mathcal J}(\mu)$, we def\/ine a map
\begin{equation*}
\begin{split}
\Psi:  \qquad \tilde{\mathcal J}(\mu)\;\;
& \longrightarrow \;\; {\mathcal U}(\mu)\\
\big(r^{(a)}_{i,\alpha}\big)_{(ai) \in \overline{H},\alpha \in {\mathbb Z}}\;
& \longmapsto \;
\big(u^{(a)}_{i,\alpha}\big)_{(ai\alpha)\in H}
\end{split}
\end{equation*}
by
\begin{equation}\label{eq:sce2}
\sum_{(bj\beta)\in H}
A_{ai\alpha, bj\beta} \,u^{(b)}_{j,\beta}=
\frac{1}{2}\big(p^{(a)}_i + m^{(a)}_i + 1\big)
+r^{(a)}_{i,\alpha}+\alpha-1\quad \text{for }\;
(ai\alpha)\in H.
\end{equation}

\begin{theorem}\label{k:th:ur}
The map $\Psi$ induces the bijection
between the set ${\mathcal J}(\mu)$ of angle variables
and the set ${\mathcal U}(\mu)$ of Bethe roots.
\end{theorem}

The proof is parallel with \cite[Section 4.2]{KTT}
for the $n=1$ case.
The induced bijection will also be denoted by $\Psi$.
It also induces the action of the time evolutions
${\mathcal T}$ on ${\mathcal U}(\mu)$
by $T_l^{(r)}({\bf u}) = \Psi\big(T^{(r)}_l({\bf r})\big)$
for ${\bf u} = \Psi({\bf r})$.
Under the condition~(\ref{k:pp}),
Conjecture~\ref{t:conj:aug21_5} and Theorem~\ref{k:th:ur} lead to
the following commutative diagram among
the level set,  the angle variables and the Bethe roots:
\begin{equation}\label{eq:cd5}
\begin{CD}
{\mathcal P}(\mu) @>{\Phi}>> {\mathcal J}(\mu) @>{\Psi}>>
{\mathcal U}(\mu) \\
@V{{\mathcal T}}VV @V{{\mathcal T}}VV @V{{\mathcal T}}VV\\
{\mathcal P}(\mu) @>{\Phi}>> {\mathcal J}(\mu) @>{\Psi}>>
{\mathcal U}(\mu)
\end{CD}
\end{equation}
Their cardinality is given by
\begin{equation*}
\vert{\mathcal P}(\mu)\vert
= \vert{\mathcal J}(\mu)\vert
= \vert {\mathcal U}(\mu)\vert =
\Omega(\mu)
\end{equation*}
with $\Omega(\mu)$ def\/ined in (\ref{t:eq:aug19_1}).
We will argue the time evolutions of the Bethe roots further
in Section \ref{k:sec:bv}.

\subsection{General case}\label{k:sec:gen}
Let ${\mathcal P}(\mu)$ be a level set.
{}From $\mu$, specify the data like
$m^{(a)}_i$, $l^{(a)}_i$ and the vacancy number $p^{(a)}_i$
by (\ref{k:rca0})--(\ref{k:gdef}).
In Sections \ref{k:sec:lt}--\ref{k:sec:dp},
we have considered the case $\forall \, p^{(a)}_i \ge 1$ (\ref{k:pp}).
In this subsection we treat
the general conf\/iguration, namely we assume
\begin{equation}\label{k:pg}
\mu=\big(\mu^{(1)}, \ldots, \mu^{(n)}\big) \;
\text{ satisf\/ies }\; p^{(a)}_i \ge 0 \;\text{ for all }
(a,i) \in \overline{H}.
\end{equation}
It turns out that the linearization scheme
remains the same
provided one discards some time evolutions
and restricts the dynamics to a subgroup
${\mathcal T}'$ of ${\mathcal T}$.
We begin by preparing some notations about partitions.

Let $\lambda, \nu$ be partitions or equivalently Young diagrams.
We def\/ine $\lambda \cup \nu$ to be the partition
whose parts are those of $\lambda$ and $\nu$,
arranged in decreasing order \cite{M}.
For example, if $\lambda = (4221)$ and
$\nu = (331)$, then $\lambda \cup \nu = (4332211)$.

We denote by $\Box^{(a)}_i$ the $(a,i)$ block
in the diagram $\mu^{(a)}$ in the sense of (\ref{k:rca0}), and by
$2 \Box^{(a)}_i$ the corresponding
block in $\mu^{(a)} \cup \mu^{(a)}$, which is a
$\big(2 m^{(a)}_i\big) \times l^{(a)}_i$ rectangle.

Let $\lambda$ be a partition.
We say that $\lambda$ {\em covers}
the block $\Box^{(a)}_i$ in $\mu^{(a)}$
if $\Box^{(a)}_i$ is inside the diagram~$\lambda$
when the two diagrams $\lambda$ and $\mu^{(a)}$ are so placed as
their top-left corners coincide.

\begin{definition}
The block $\Box^{(a)}_i$ is {\em null} if $p_i^{(a)} =0$.
\end{definition}

\begin{definition}
The block $\Box^{(a)}_i$ is {\em convex} if $\big(\mu^{(a-1)} \cup \mu^{(a+1)}\big)$ covers
the block $2 \Box^{(a)}_i$ in $\mu^{(a)} \cup \mu^{(a)}$.
\end{definition}
Here we interpret $\mu^{(0)} = (1^L)$ and $\mu^{(n+1)} = \varnothing$.

Say that the time evolution $T^{(a)}_l$ is {\em inadmissible}
to $\mu$ if there exists $1 \leq i \leq g_a$ such that
$\Box^{(a)}_i$ is null and convex, and $l^{(a)}_i > l > l^{(a)}_{i+1}$.
Otherwise we say that $T^{(a)}_l$ is {\em admissible} to $\mu$.
Here $l^{(a)}_{g_a+1}$ is to be understood as $0$.
Note that when $l^{(a)}_i - l^{(a)}_{i+1} = 1$,
the block $\Box^{(a)}_i$ does not cause an
inadmissible time evolution even if it is null and convex.
We will also say that $T^{(a)}_l$ is admissible or inadmissible
to a path $p \in {\mathcal P}(\mu)$
depending on whether $T^{(a)}_l$ is admissible or
inadmissible to~$\mu$.

As a generalization of Conjecture \ref{t:conj:aug10_3},
we propose

\begin{conjecture}\label{t:conj:aug11_1}
$T^{(r)}_l({\mathcal P}(\mu)) = {\mathcal P}(\mu)$
holds if and only if $T^{(r)}_l$ is admissible to $\mu$.
Thus, $T^{(r)}_l({\mathcal P}(\mu)) = {\mathcal P}(\mu)$ is valid
for any $r, l$ if and only if there is no block
$\Box^{(a)}_i$ that is null, convex and of length
$l^{(a)}_i$ such that $l^{(a)}_i - l^{(a)}_{i+1} \ge 2$.
\end{conjecture}

\begin{example}
Consider the following conf\/igurations corresponding to
the evolvable paths
$p_I=12 11 22 333 111 11 22 22 233 333$ (upper one)
and $p_{II}=111 222 133 321 111 222 233 333$ (lower one).
Hatched blocks are null and convex.
They show that
$T^{(2)}_4$ is inadmissible to $p_I$, and
$T^{(2)}_1$ and $T^{(2)}_2$ are inadmissible to $p_{II}$.
$$
\includegraphics{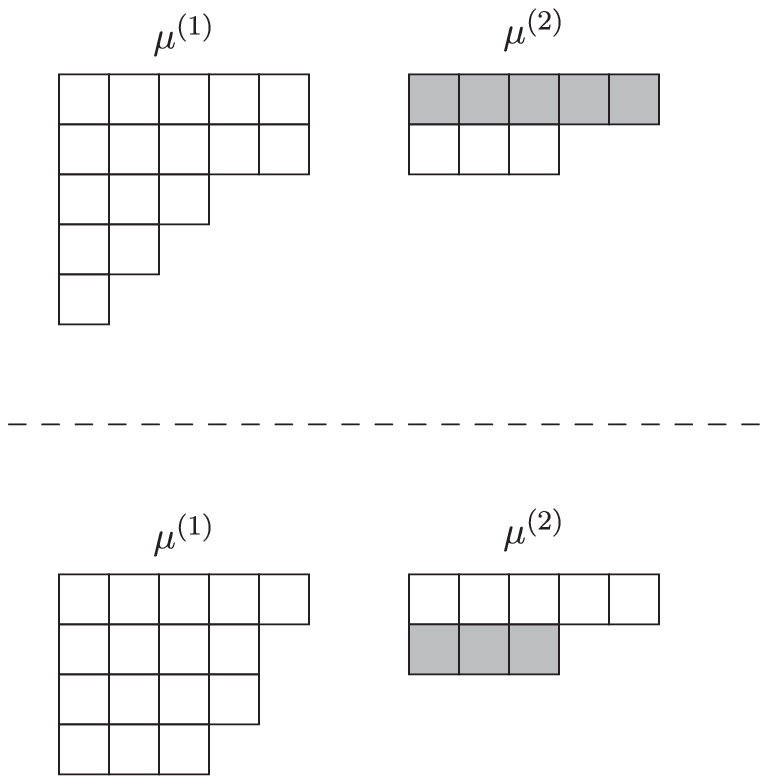}
$$
\end{example}

Let ${\mathcal T}'$ be the abelian group generated by all
the time evolutions $T_l^{(a)}$ admissible to $\mu$.
Then Conjecture \ref{t:conj:aug11_1}
implies that ${\mathcal T}'$ acts on the level set ${\mathcal P}(\mu)$.
We def\/ine the connected component
$\Sigma' (p)$ to be the ${\mathcal T}'$-orbit in
${\mathcal P}(\mu)$ that contains $p$.

As it turns out,
the restriction of the dynamics from ${\mathcal T}$ to
${\mathcal T}'$ will be matched by introducing a sub-lattice
$\mathbb{L}$ of ${\mathbb Z}^g$ as follows:
\begin{equation*}
\mathbb{L} = \bigoplus_{(ai)\in \overline{H}}
{\mathbb Z} h^{(a)}_{\xi^{(a)}_i} \subset {\mathbb Z}^g,
\qquad
\xi^{(a)}_i =
\begin{cases}
l^{(a)}_i & \mbox{if $\Box^{(a)}_i$ is null and convex,}\\
l^{(a)}_{i+1} + 1 & \mbox{otherwise.}
\end{cases}
\end{equation*}
See (\ref{k:ter}) for the def\/inition of $h^{(a)}_l$.
The subgroup ${\mathcal T}' \subset {\mathcal T}$
and $\xi^{(a)}_i$ here are analogues of~$\overline{\mathcal T}$ and~$\eta^{(a)}_i$
in Remark~\ref{t:rem:aug18_2} in the present setting.
It is easy to check
\begin{lemma}\label{t:lem:aug25_1}
Let $\mu = \big(\mu^{(1)}, \ldots, \mu^{(n)}\big)$
be a configuration, i.e., \eqref{k:pg} is satisfied.
Suppose the $(a,i)$ block of $\mu$ is null and convex.
Then:
\begin{enumerate}\itemsep=0pt
\item[$(i)$] $p^{(a)}_{i+1} = 0$.

\item[$(ii)$] For $b = a \pm 1$ there is no $1 \leq j \leq g_b$ such that
$l^{(a)}_i > l^{(b)}_j > l^{(a)}_{i+1}$.
\end{enumerate}
\end{lemma}

Recall that the matrix $F_{\boldsymbol{\gamma}}$ is def\/ined by
(\ref{k:f}) and (\ref{k:ftdef}).
\begin{proposition}
$F_{\boldsymbol{\gamma}}\, \mathbb{Z}^g$
is a sub-lattice of $\mathbb{L}$.
\end{proposition}

\begin{proof}
Let $\tilde{F}_a$ be the
$g_a \times (g_{a-1} + g_a + g_{a+1})$
sub-matrix of $F_{\boldsymbol{\gamma}}$
for color $a$ row and colors~$a$,~$a \pm 1$ column indices.
Its $(bj)$th column $(b=a, a\pm 1)$
is the array $\big(F_{ai,bj}/\gamma^{(b)}_j\big)_{1\le i \le g_a}$.
Let $\bar{h}_{\xi^{(a)}_i}^{(a)} \in {\mathbb Z}^{g_a}$
be the sub-vector obtained from $h_{\xi^{(a)}_i}^{(a)}$
by extracting the color $a$ rows.
We are to show that for any~$a$,
all the column vectors of $\tilde{F}_a$ can be expressed as
a linear combination of $\bar{h}_{\xi^{(a)}_i}^{(a)}$ with integer coef\/f\/icients.
For simplicity we demonstrate it by an example.
The general case is similar.

Suppose $g_a = 7$ and the $(a,i)$ block is
null and convex if and only if $i=1,2,6,7$.
Then the vectors $\bar{h}_{\xi^{(a)}_i}^{(a)} \; (1 \leq i \leq 7)$ are given by
the columns of the matrix in the left hand side of the following equation:
\begin{equation}\label{t:eq:aug25_2}
\begin{pmatrix}
l_1 & l_2 & l_4+1 & l_5+1 & l_6+1 & l_6 & l_7 \\
l_2 & l_2 & l_4+1 & l_5+1 & l_6+1 & l_6 & l_7 \\
l_3 & l_3 & l_4+1 & l_5+1 & l_6+1 & l_6 & l_7 \\
l_4 & l_4 & l_4   & l_5+1 & l_6+1 & l_6 & l_7 \\
l_5 & l_5 & l_5   & l_5   & l_6+1 & l_6 & l_7 \\
l_6 & l_6 & l_6   & l_6   & l_6   & l_6 & l_7 \\
l_7 & l_7 & l_7   & l_7   & l_7   & l_7 & l_7
\end{pmatrix}
\simeq
\begin{pmatrix}
l_1 & l_2 & 1 &   &    &   &   \\
l_2 & l_2 & 1 &   &    &   &   \\
l_3 & l_3 & 1 &   &    &   &   \\
  &   &       & 1 &    &   &   \\
  &   &       &   &  1 &   &   \\
  &   &       &   &    & l_6 & l_7 \\
  &   &       &   &    & l_7 & l_7
\end{pmatrix},
\end{equation}
where $l_i = l_i^{(a)}$.
The right hand side is obtained from the left
(and vice versa) by elementary transformations
adding integer multiple of one column to other columns.

By $(i)$ of Lemma \ref{t:lem:aug25_1}
we have $p^{(a)}_3 = 0$, and by $(ii)$ of the same lemma
any column of $\tilde{F}_a$ is an integer multiple
(and possible integer addition to the 4-th and 5-th rows) of
one of the following vectors:
\begin{gather*}
\begin{pmatrix}
l_1  \\
l_2  \\
l_3  \\
l_4  \\
l_5  \\
l_6  \\
l_7
\end{pmatrix}
\mbox{for $l^{(b)}_j \geq l_1$},\quad
\begin{pmatrix}
l_2  \\
l_2  \\
l_3  \\
l_4  \\
l_5  \\
l_6  \\
l_7
\end{pmatrix}
\mbox{for $l^{(b)}_j = l_2$},\quad
\begin{pmatrix}
l  \\
l  \\
l  \\
x  \\
y  \\
l_6  \\
l_7
\end{pmatrix}
\mbox{for $l_3 \geq l^{(b)}_j = l \geq l_6$},\quad
\begin{pmatrix}
l_7  \\
l_7  \\
l_7  \\
l_7  \\
l_7  \\
l_7  \\
l_7
\end{pmatrix}
\mbox{for $l^{(b)}_j = l_7$},
\end{gather*}
where $x$ and $y$ are some integers.
Clearly these vectors can be expressed
by linear combinations of the columns
of the matrices in~(\ref{t:eq:aug25_2}) with integer coef\/f\/icients.
\end{proof}

\begin{proposition}\label{k:pr:ll}
The volume of the unit cell of the lattice
$\mathbb{L}$ is equal to $\prod \big(l^{(a)}_i - l^{(a)}_{i+1}\big)$,
where the product is taken over all the null and convex
blocks $\Box^{(a)}_i$.
\end{proposition}

\begin{proof}
By using the above example we have
\begin{equation*}
\begin{pmatrix}
l_1 & l_2 & 1 &   &    &   &   \\
l_2 & l_2 & 1 &   &    &   &   \\
l_3 & l_3 & 1 &   &    &   &   \\
  &   &       & 1 &    &   &   \\
  &   &       &   &  1 &   &   \\
  &   &       &   &    & l_6 & l_7 \\
  &   &       &   &    & l_7 & l_7
\end{pmatrix}
\simeq
\begin{pmatrix}
l_1 -l_2& l_2 - l_3& 1 &   &    &   &   \\
        & l_2 - l_3& 1 &   &    &   &   \\
                &  & 1 &   &    &   &   \\
  &   &                & 1 &    &   &   \\
  &   &                &   &  1 &   &   \\
  &   &      &   &    & l_6 - l_7 & l_7 \\
  &   &       &   &    &          & l_7
\end{pmatrix}.
\end{equation*}
The general case is similar.
\end{proof}

Let the abelian group ${\mathcal T}'$ act
on the torus $\mathbb{L} /
F_{\boldsymbol{\gamma}}\,\mathbb{Z}^g$
by (\ref{t:eq:aug11_4}).
This action is transitive, since~${\mathcal T}'$ contains
the free abelian subgroup generated by all
$T_{\xi^{(a)}_i}^{(a)}$'s.

\begin{conjecture}\label{k:con:pg}
Under the condition \eqref{k:pg},
$\Sigma'(p) \cap {\mathcal P}_+(\mu) \ne \varnothing$ holds
for any $p \in {\mathcal P}(\mu)$.
\end{conjecture}

Owing to this property, one can introduce $\Phi'$ by the same scheme
as (\ref{k:ds}) with ${\mathcal T}$
replaced by ${\mathcal T}'$.
Similarly, let
$\chi': \; {\mathcal T}' \cdot ({\mathcal A} \cdot x)
\longrightarrow \mathbb{L}/F_{\boldsymbol{\gamma}}\,{\mathbb Z}^g$
be def\/ined as in (\ref{k:chi}).
As the generalization of Conjecture \ref{t:conj:aug10_8},
we have

\begin{conjecture}\label{t:conj:aug11_6}
For any path
$p \in {\mathcal P}_{\boldsymbol{\gamma}}(\mu)$,
the map $\Phi'_{\chi} :=\chi' \circ \Phi'$
gives a bijection between
the connected component $\Sigma'(p)$ and the torus
$\mathbb{L} / F_{\boldsymbol{\gamma}}\,\mathbb{Z}^g$
making the following diagram commutative:
\begin{gather*}
\begin{CD}
\Sigma'(p)
@>{\Phi'_\chi}>> \mathbb{L} / F_{\boldsymbol{\gamma}}\,{\mathbb Z}^g \\
@V{\mathcal T}'VV @VV{\mathcal T}'V\\
\Sigma'(p)
@>{\Phi'_\chi}>>\mathbb{L} / F_{\boldsymbol{\gamma}}\,{\mathbb Z}^g
\end{CD}
\end{gather*}
\end{conjecture}

By Proposition \ref{k:pr:ll}, this conjecture implies
the relation
$| \Sigma' (p) | = \det F_{\boldsymbol{\gamma}}/\prod \big(l^{(a)}_i - l^{(a)}_{i+1}\big)$,
where the product is taken over all null and convex blocks
$\Box^{(a)}_i$.
Note that this factor is common to all the connected components
in the level set.
Thus Conjecture~\ref{t:conj:aug11_6} also tells that
\begin{equation*}
|{\mathcal P}(\mu)|
= \frac{\Omega (\mu)}{\prod (l^{(a)}_i - l^{(a)}_{i+1})},
\end{equation*}
which ref\/ines \cite[Conjecture 4.2]{KT2}.

\begin{example}
Consider the path $p_{III}=121212343434121212$ of length $L=18$.
Its soliton content is given in the following table.
\begin{center}
\begin{tabular}{l|ccc|lc|cc}
\hline
\tsep{1.5mm}\bsep{2mm} $(ai)$ & $l^{(a)}_i$ & $m^{(a)}_i$ & $p^{(a)}_i$ & $r^{(a)}_{i,\alpha}$ & $\gamma^{(a)}_i$ & $\xi^{(a)}_i$ & ${}^t  h^{(a)}_{\xi^{(a)}_i}$\\
\hline
$(11)^*$ & 3 & 1 & 0 & $0$& 1 &3 & (3,1,0,0)\\
$(12)^*$ & 1 & 9 & 0 & $0,\ldots,0$ & 9 & 1 & (1,1,0,0)\\
$(21)$ & 3 & 2 & 3 & $0,0$ & 1 & 1 & (0,0,1,0)\\
$(31)^*$ & 3 & 1 & 0 & 0 & 1 & 3 & (0,0,0,1) \\
\hline
\end{tabular}
\end{center}
Its $(1,1)$, $(1,2)$, $(3,1)$ blocks marked by $\ast$ are null and convex.
The soliton content is depicted as follows.
(Null and convex blocks are hatched.)
$$
\includegraphics{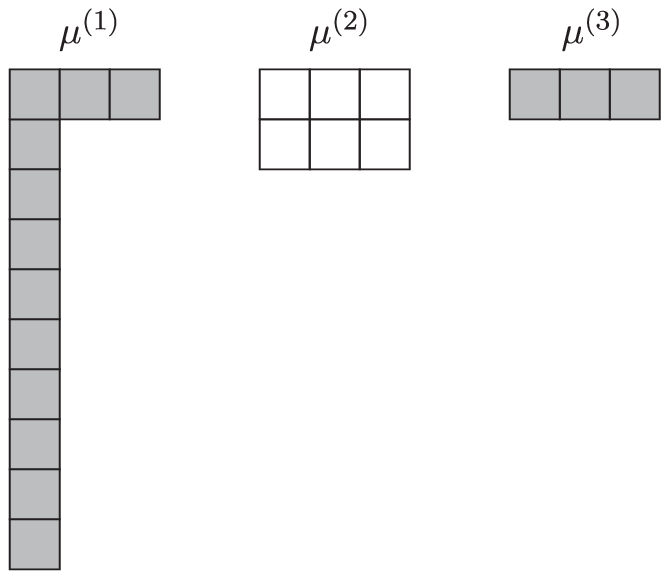}
$$
$T^{(1)}_2$, $T^{(3)}_1$ and $T^{(3)}_2$ are inadmissible to $p_{III}$, while all the other $T^{(a)}_l$'s are admissible.
The matrix~$F_{\boldsymbol{\gamma}}$ for $p_{III}$ is given by
\begin{displaymath}
F_{\boldsymbol{\gamma}}=
\left(
\begin{array}{rrrr}
 6 & 2 & -6 & 0 \\
 2 & 2 & -2 & 0 \\
 -3 & -1 & 15 & -3 \\
 0 & 0 & -6 & 6
\end{array}
\right).
\end{displaymath}
The cardinality of the connected component
in the level set is calculated as
\[
| \Sigma' (p_{III}) | = \det F_{\boldsymbol{\gamma}}
/\big(l^{(1)}_1 - l^{(1)}_{2}\big)\big(l^{(1)}_2 - l^{(1)}_{3}\big)
\big(l^{(3)}_1 - l^{(3)}_{2}\big)
= 432/(2 \cdot 1 \cdot 3) = 72.
\]
\end{example}

\begin{example}
Consider the path $p_{IV}=1122331142233444$ of length $L=16$.
Its soliton content is given in the following table.
\begin{center}
\begin{tabular}{l|ccc|lc|cc}
\hline
\tsep{1.5mm}\bsep{2mm} $(ai)$ & $l^{(a)}_i$ & $m^{(a)}_i$ & $p^{(a)}_i$
& $r^{(a)}_{i,\alpha}$ & $\gamma^{(a)}_i$
& $\xi^{(a)}_i$ & ${}^t  h^{(a)}_{\xi^{(a)}_i}$\\
\hline
$(11)$   & 4 & 1 & 0 & $0$        & 1 & 3 & (3,2,0,0,0)\\
$(12)$   & 2 & 4 & 2 & $1,1,0,0$  & 2 & 1 & (1,1,0,0,0)\\
$(21)^*$ & 4 & 1 & 0 & $0$        & 1 & 4 & (0,0,4,2,0)\\
$(22)^*$ & 2 & 2 & 0 & $0,0$      & 2 & 2 & (0,0,2,2,0)\\
$(31)$   & 4 & 1 & 0 & $0$        & 1 & 1 & (0,0,0,0,1)\\
\hline
\end{tabular}
\end{center}
Its $(2,1)$, $(2,2)$ blocks marked by $\ast$ are null and convex.
The soliton content is depicted in the following
with the null and convex blocks hatched,
showing that $T^{(2)}_1$ and $T^{(2)}_3$ are inadmissible to $p_{IV}$.
$$
\includegraphics{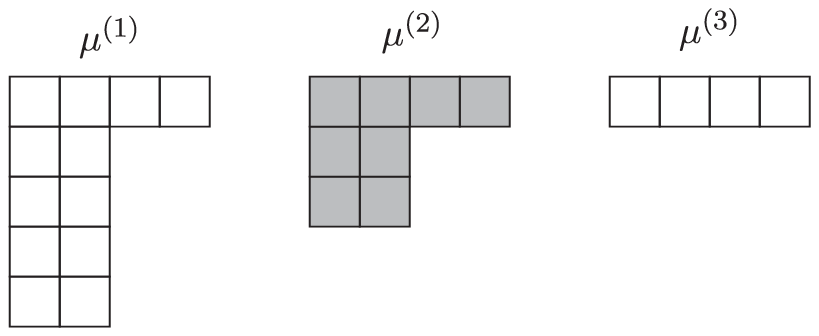}
$$
The matrix $F_{\boldsymbol{\gamma}}$ for $p_{IV}$ is given by
\begin{displaymath}
F_{\boldsymbol{\gamma}}=\left(
\begin{array}{rrrrr}
 8 & 8 & -4 & -2 & 0 \\
 4 & 9 & -2 & -2 & 0 \\
 -4 & -4 & 8 & 4 & -4 \\
 -2 & -4 & 4 & 4 & -2 \\
 0 & 0 & -4 & -2 & 8
\end{array}
\right)
\end{displaymath}
The cardinality of the connected component is calculated as
\[
| \Sigma' (p_{IV}) | =
\det F_{\boldsymbol{\gamma}}/
\big(l^{(2)}_1 - l^{(2)}_{2}\big)\big(l^{(2)}_2 - l^{(2)}_{3}\big)
= 2048/(2 \cdot 2) = 512.
\]
\end{example}

\subsection{Summary of conjectures}\label{sec:suma}

Here is a summary of our conjectures.
Each one is based on those in the preceding lines.
Thus the principal ones are the linearizations.
\begin{center}
\begin{tabular}{c|c|c}
Condition on soliton/string content $\mu$
& $\forall\,  p^{(a)}_i \ge 1$  \;(\ref{k:pp})
&$\forall \, p^{(a)}_i \ge 0$  \;(\ref{k:pg}) \\
\hline
\tsep{1mm} Stability of dynamics
&Conjecture \ref{t:conj:aug10_3} \;\;${\mathcal T}$
&Conjecture \ref{t:conj:aug11_1} \;\;${\mathcal T}'$ \\
Intersection with highest paths
&Conjecture \ref{t:conj:aug10_5} \;\;\;\;\;
& Conjecture \ref{k:con:pg}\;\;\;\;\;\;\; \\
Linearization
&Conjecture \ref{t:conj:aug21_5} $\simeq$
\ref{t:conj:aug10_8} \;$\Phi_\chi$
&Conjecture \ref{t:conj:aug11_6} \;\;$\Phi'_\chi$ \bsep{1mm} \\
\hline
\tsep{1mm} Relevant torus
& ${\mathbb Z}^g/F_{\boldsymbol{\gamma}}\,{\mathbb Z}^g$
& $\mathbb{L} /F_{\boldsymbol{\gamma}}\,\mathbb{Z}^g$
\end{tabular}
\end{center}

\section{Tropical Riemann theta from combinatorial Bethe ansatz}
\label{k:sec:trt}

\subsection{From tropical tau function to tropical Riemann theta function}

Let $\big(\mu^{(a)}, r^{(a)}\big)$ be the color $a$ part of the
rigged conf\/iguration $(\mu, {\bf r})$
depicted in the right diagram of (\ref{k:rca0}).
We keep the notations in (\ref{k:hdef})--(\ref{k:gdef}).
An explicit formula for the image path of the KKR map $\phi^{-1}$
(\ref{k:kkr}) is known in terms of the
{\em tropical tau function} \cite{KSY}.
It is related to the charge function on  rigged conf\/igurations
and is actually obtained from the tau function in the KP hierarchy \cite{JM}
by the ultradiscretization
with an elaborate adjustment of parameters from the KP and
rigged conf\/igurations:
\begin{gather}
\tau_{k,d} = -\min\Bigg\{
\frac{1}{2}\sum_{(a i \alpha),\,(b j \beta)}\!\!C_{ab}
\min\big(l^{(a)}_i, l^{(b)}_j\big)N^{(a)}_{i,\alpha}N^{(b)}_{j,\beta}
+\sum_{(a i \alpha)}r^{(a)}_{i,\alpha}N^{(a)}_{i \alpha} \nonumber\\
\qquad\qquad\qquad{}-k\sum_{(i\alpha)}N^{(1)}_{i,\alpha}
+\sum_{(j \beta)}l^{(d)}_jN^{(d)}_{j,\beta}\Bigg\}\qquad
(1 \le d \le n+1).\label{k:tau1}
\end{gather}
Here the sums range over $(ai\alpha) \in H$ wherever
$N^{(a)}_{i,\alpha}$ is involved.
Thus in the second line of (\ref{k:tau1}),
$(i\alpha)$ runs over
$1 \le i \le g_1,\, 1 \le \alpha \le m^{(1)}_i$ and so does
$(j\beta)$ over  $1 \le j \le g_d,\, 1 \le \beta \le m^{(d)}_j$.
The last term $\sum_{(j \beta)}l^{(d)}_jN^{(d)}_{j,\beta}$
is to be understood as 0 when $d=n+1$.
In (\ref{k:tau1}), $\min$ is taken over
$N^{(a)}_{i,\alpha}\in \{0,1\}$ for all $(a i \alpha) \in H$.
Thus it consists of $2^G$ candidates.
Obviously, $\tau_{k,d} \in {\mathbb Z}_{\ge 0}$ holds.

In the present case,  the paths are taken from
$(B^{1,1})^{\otimes L}$.
We parameterize the set
$B^{1,1}=\{1,2,\ldots, n+1\}$
as
$\{(x_1,x_2,\ldots, x_{n+1}) \in \{0,1\}^{n+1}\mid
x_1+\cdots + x_{n+1}=1\}$.

\begin{theorem}[\cite{KSY}]\label{k:th:ksy}
The image of the KKR map $\phi^{-1}$ is expressed as follows:
\begin{gather}
\phi^{-1}: \
(\mu, {\bf r})  \mapsto x_1\otimes x_2 \otimes \cdots \otimes x_L
\in (B^{1,1})^{\otimes L}, \nonumber\\
\phantom{\phi^{-1}: }{} \ \ x_k  = (x_{k,1}, x_{k,2}, \ldots, x_{k,n+1}) \in B^{1,1},
\nonumber\\
\phantom{\phi^{-1}: }{} \ \ x_{k,a} = \tau_{k,a}-\tau_{k-1,a}-\tau_{k,a-1}+\tau_{k-1,a-1}
\qquad (2 \le a \le n+1). \label{k:db}
\end{gather}
\end{theorem}
In the context of the box-ball system,
$x_{k,a}$ represents the number of balls of color $a$
in the $k$th box from the left for $2 \le a \le n+1$.
The remaining $x_{k,1}$ is determined from this by
$x_{k,1} = 1-x_{k,2}-\cdots - x_{k,n+1}$, which
also takes values $x_{k,1}=0,1$ in the present case.

The following is a special case ($\forall\, \lambda_k=1$)
of \cite[Proposition 5.1]{KSY}.
\begin{proposition}[\cite{KSY}]\label{k:pr:tbl}
The tropical tau function satisfies the tropical Hirota equation:
\begin{equation*}
\tau_{k-1,d} + \overline{\tau}_{k,d-1} =
\max(\overline{\tau}_{k,d} + \tau_{k-1,d-1},\,
\tau_{k,d}+\overline{\tau}_{k-1,d-1}-1)\qquad
(2 \le d \le n+1),
\end{equation*}
where $\overline{\tau}_{k,d}$ is obtained from
$\tau_{k,d}$ by replacing
$r^{(a)}_{i,\alpha}$ with $r^{(a)}_{i,\alpha}+l^{(a)}_i\delta_{a1}$.
\end{proposition}

\begin{lemma}\label{k:lem:ext}
Let $p, p' \in \big(B^{1,1}\big)^{\otimes L}$ be the highest paths
whose rigged configurations are
$(\mu, {\bf r})$ and $(\mu', {\bf r}')$, respectively.
Then the rigged configuration of the highest path
$p\otimes p' \in \big(B^{1,1}\big)^{\otimes 2L}$ is
$(\mu, {\bf r}) \sqcup (\mu', {\tilde {\bf r}}')$, where
$({\tilde {\bf r}}')
= \big(r^{(a)\prime}_{i,\alpha}+p^{(a)}_i\big)_{(ai\alpha)\in H}$ and
$p^{(a)}_i$ is the vacancy number of $(\mu, {\bf r})$.
\end{lemma}

In the lemma, $\sqcup$ stands for the sum (union)
as the multisets of strings,
namely, the rows of $\mu$ attached with riggings.
More formally a string is a triple $\big(a,l^{(a)}_i, r^{(a)}_{i,\alpha}\big)$
consisting of its color~$a$, length $l^{(a)}_i$ and rigging
$r^{(a)}_{i,\alpha}$ that is labeled with $H$.
Let $p \in \big(B^{1,1}\big)^{\otimes L}$ be the highest path
corresponding to the rigged conf\/iguration $(\mu, {\bf r})$.
From Lemma \ref{k:lem:ext}, the rigged conf\/iguration of
the highest path $p^{\otimes M}\in
\big(B^{1,1}\big)^{\otimes ML}$ is
$\big(\mu, {{\bf r}}^1\big)
\sqcup \big(\mu, {{\bf r}}^2\big) \sqcup \cdots
\sqcup \big(\mu, {{\bf r}}^{M}\big)$,
where
${{\bf r}}^k
= \big(r^{(a)}_{i,\alpha}+(k-1)p^{(a)}_i\big)_{(ai\alpha)\in H}$.
Pictorially, this corresponds to extending
the right diagram in (\ref{k:rca0})
vertically $M$ times adding
$p^{(a)}_i, 2p^{(a)}_i, \ldots, (M-1)p^{(a)}_i$
to the riggings in the $(a,i)$ block for each $(a, i)$.
This reminds us of (\ref{k:qpr}), and is in fact the origin of
the extended rigged conf\/iguration in Section~\ref{k:sec:aav}.

We proceed to the calculation of the tropical tau function
$\tau^M_{k,d}$ associated with the above rigged conf\/iguration
$\big(\mu, {{\bf r}}^1\big)
\sqcup \big(\mu, {{\bf r}}^2\big) \sqcup \cdots
\sqcup \big(\mu, {{\bf r}}^{M}\big)$.
In (\ref{k:tau1}), the variable $N^{(a)}_{i,\alpha}$
is to be replaced by the $M$ replicas
$N^{(a)}_{i,\alpha,1},\ldots, N^{(a)}_{i,\alpha, M}$
to cope with the $M$-fold extension:
\begin{gather}
\tau^M_{k,d} = -\min\Bigg\{
\frac{1}{2}\sum_{(a i \alpha),\,(b j \beta)}
\sum_{1\le s, t \le M}\!\!C_{ab}
\min\big(l^{(a)}_i, l^{(b)}_j\big)N^{(a)}_{i,\alpha,s}N^{(b)}_{j,\beta,t}\nonumber\\
 \qquad\qquad\qquad{}
+\sum_{(a i \alpha)}\sum_{1\le s \le M}
\big(r^{(a)}_{i,\alpha}+(s-1)p^{(a)}_i\big)N^{(a)}_{i,\alpha,s}\nonumber\\
 \qquad\qquad\qquad{}
-k\sum_{(i \alpha)}\sum_{1\le s \le M}N^{(1)}_{i,\alpha,s}
+\sum_{(j \beta)}\sum_{1 \le t \le M}l^{(d)}_jN^{(d)}_{j,\beta,t}\Bigg\}.\label{k:tau2}
\end{gather}
All the summands here are invariant under permutations within
$N^{(a)}_{i,\alpha,1},\ldots, N^{(a)}_{i,\alpha, M}$
for each $(a,i)$ except the replica symmetry breaking term
$\sum_{(a i \alpha)}
\sum_{1\le s \le M}(s-1)p^{(a)}_iN^{(a)}_{i,\alpha,s}$.
Due to $p^{(a)}_i \ge 0$, one can reduce the minimizing variables
$N^{(a)}_{i,\alpha, s}\in \{0,1\}$ to
$n^{(a)}_{i,\alpha}\in \{0,1,\ldots, M-1\}$ such that
\begin{equation}\label{emv}
N^{(a)}_{i,\alpha,1}= N^{(a)}_{i,\alpha,2}
= \cdots = N^{(a)}_{i,\alpha, n^{(a)}_{i,\alpha}}=1,
\qquad
N^{(a)}_{i,\alpha, n^{(a)}_{i,\alpha}+1}=\cdots =
N^{(a)}_{i,\alpha, M}=0.
\end{equation}
As the result, (\ref{k:tau2}) becomes
\begin{gather*}
\tau^M_{k,d} = -\min\Bigg\{
\frac{1}{2}\sum_{(a i \alpha),\,(b j \beta)}\!\!C_{ab}
\min\big(l^{(a)}_i, l^{(b)}_j\big)n^{(a)}_{i,\alpha}n^{(b)}_{j,\beta}\\
 \qquad\qquad\qquad{}
+\sum_{(a i \alpha)}\Bigg(
r^{(a)}_{i,\alpha}n^{(a)}_{i,\alpha}
+\frac{n^{(a)}_{i,\alpha}\big(n^{(a)}_{i,\alpha}-1\big)}{2}p^{(a)}_i\Bigg)-k\sum_{(i \alpha)}n^{(1)}_{i,\alpha}
+\sum_{(j \beta)}l^{(d)}_jn^{(d)}_{j,\beta}\Bigg\},
\end{gather*}
where the minimum is now taken over
$n^{(a)}_{i,\alpha}\in \{0,1,\ldots, M-1\}$ for all the blocks $(a,i)$.
The notation can be eased considerably by
introducing a quadratic form of
${\bf n}=\big(n^{(a)}_{i,\alpha}\big)_{(a i \alpha)\in H}$ as follows:
\begin{gather}
\tau^M_{k,d} = -\min\left\{
\frac{1}{2}{}^t{\bf n}B{\bf n}+
{}^t\Bigl({\bf r}-\frac{\bf p}{2}
-k{\bf h}^{(1)}_1+{\bf h}^{(d)}_\infty
\Bigr){\bf n}\right\}\quad (1\le d \le n+1, 1 \le k \le ML),\!\!\!\label{k:tau4}\\
B  =\big(\delta_{ab}\delta_{ij}\delta_{\alpha\beta}p^{(a)}_i+
C_{ab}\min\big(l^{(a)}_i, l^{(b)}_j\big)\big)_{(a i \alpha), (b j \beta) \in H},
\label{k:B}\\
{\bf r}  =\big(r^{(a)}_{i,\alpha}\big)_{(a i \alpha)\in H},\qquad
{\bf p} =\big(p^{(a)}_i\big)_{(a i \alpha)\in H},\label{k:rp}\\
{\bf h}^{(c)}_l
=\big(\delta_{ac}\min\big(l, l^{(c)}_i\big)\big)_{(a i \alpha)\in H}\quad
(1 \le c \le n, l\ge 1), \qquad {\bf h}^{(n+1)}_l=0.\label{k:hv}
\end{gather}
The $G$-dimensional vector ${\bf h}^{(c)}_l$ here
should be distinguished from the $g$-dimensional one \linebreak $h^{(c)}_l$~(\ref{k:ter}).
It is the velocity vector of the time evolution $T^{(c)}_l$
in $\tilde{\mathcal J}(\mu)$ and ${\mathcal J}(\mu)$
in the light of~(\ref{k:ter2}).
The matrix $B$ is symmetric and
positive def\/inite. See for example  \cite[Lemma 3.8]{KN}.
We call $B$ (\ref{k:B}) the tropical period matrix
although a connection to a tropical curve is yet to be clarif\/ied.
We also introduce the $G$-dimensional vector
${\bf 1} = (1)_{(a i \alpha)\in H}$.
Then from~(\ref{k:pdef}) and~(\ref{k:hv}) we get
\begin{gather}
B{\bf 1}  = \Bigg(\sum_{(b j \beta)\in H}
\big(\delta_{ab}\delta_{ij}\delta_{\alpha\beta}p^{(a)}_i+
C_{ab}\min\big(l^{(a)}_i, l^{(b)}_j\big)\big)\Bigg)_{(a i \alpha)}
=(L\delta_{a1})_{(a i \alpha)} = L{\bf h}^{(1)}_1,
\label{k:bh}\\
{}^t {\bf 1}{\bf h}^{(c)}_\infty  = \big|\mu^{(c)}\big|.\label{new}
\end{gather}

\begin{example}\label{k:ex:Brp}
With the ordering of indices in $H$ (\ref{k:hido}),
the tropical period matrix $B$ (\ref{k:B}) and
the vectors ${\bf r}$, ${\bf p}$ (\ref{k:rp})
for the rigged conf\/iguration (\ref{k:rc1}) are given
as follows:
\begin{equation}\label{k:Bex}
B =
\begin{pmatrix}
10 & 6 & 4 & 4 & 4 & -3 & -1\\
6 & 10 & 4 & 4 & 4 & -3 & -1\\
4 & 4 & 11 & 4 & 4 & -2 & -1\\
4 & 4 & 4 & 11 & 4 & -2 & -1\\
4 & 4 & 4 & 4 & 11 & -2 & -1\\
-3 & -3 & -2 & -2 & -2 & 10 & 2\\
-1 & -1 & -1 & -1 & -1 & 2 & 3
\end{pmatrix},
\qquad
{\bf r} = \begin{pmatrix}
4 \\ 2 \\ 6 \\ 5 \\ 1 \\ 0 \\ 0
\end{pmatrix},
\qquad
{\bf p} = \begin{pmatrix}
4 \\ 4 \\ 7 \\ 7 \\ 7 \\ 2 \\ 1
\end{pmatrix}.
\end{equation}
The velocity vector ${\bf h}^{(a)}_l$ (\ref{k:hv})
for the time evolution $T^{(a)}_l$ is specif\/ied as
\begin{gather*}
{\bf h}^{(1)}_1 = \frac{1}{2}{\bf h}^{(1)}_2
=\begin{pmatrix}
1 \\ 1 \\ 1 \\ 1 \\ 1 \\ 0 \\ 0
\end{pmatrix},\qquad
{\bf h}^{(1)}_{l\ge 3} = \begin{pmatrix}
3 \\ 3 \\ 2 \\ 2 \\ 2 \\ 0 \\ 0
\end{pmatrix},\qquad
{\bf h}^{(2)}_1 = \begin{pmatrix}
0 \\ 0 \\ 0 \\ 0 \\ 0 \\ 1 \\ 1
\end{pmatrix},\\
{\bf h}^{(2)}_2 = \begin{pmatrix}
0 \\ 0 \\ 0 \\ 0 \\ 0 \\ 2 \\ 1
\end{pmatrix},\qquad
{\bf h}^{(2)}_3 = \begin{pmatrix}
0 \\ 0 \\ 0 \\ 0 \\ 0 \\ 3 \\ 1
\end{pmatrix},\qquad
{\bf h}^{(2)}_{l\ge 4} = \begin{pmatrix}
0 \\ 0 \\ 0 \\ 0 \\ 0 \\ 4 \\ 1
\end{pmatrix}.
\end{gather*}
\end{example}

Now we take $M$ to be even and make the shifts
${\bf  n} \rightarrow {\bf n}+\frac{M}{2}{\bf 1}$ and
$k \rightarrow k'=k+\frac{ML}{2}$.
Using~(\ref{k:bh}), we f\/ind
\begin{equation}\label{k:tau5}
\tau^M_{k', \,d}=-\min\left\{
\frac{1}{2}{}^t{\bf n}B{\bf n}+
{}^t\left({\bf r}-\frac{\bf p}{2}
-k{\bf h}^{(1)}_1+{\bf h}^{(d)}_\infty
\right){\bf n}\right\} + u_M\,{}^t{\bf h}^{(d)}_\infty{\bf 1}+v_Mk' + w_M.
\end{equation}
The f\/irst term is formally identical with (\ref{k:tau4}) but
now the minimum extends over
$-\frac{M}{2} {\le} n^{(a)}_{i, \alpha}{<}\frac{M}{2}$.
The scalars $u_M$, $v_M$, $w_M$
are independent of $k'$ and $d$, therefore
these terms are irrelevant when taking the
double dif\/ference
$\tau^M_{k', \,a}-\tau^M_{k'-1, \,a}
-\tau^M_{k', \,a-1}+\tau^M_{k'-1, \,a-1}$
as in (\ref{k:db}).

Passing to the limit $M \rightarrow \infty$,
we see that the relevant term in (\ref{k:tau5}) tends to
\begin{equation}\label{k:taut}
\lim_{M:\text{ even}\rightarrow \infty}
\left(\tau^M_{k', \,d}-
u_M\,{}^t{\bf h}^{(d)}_\infty{\bf 1}-v_Mk' - w_M\right)
= \Theta\left({\bf r}-\frac{\bf p}{2}
-k{\bf h}^{(1)}_1+{\bf h}^{(d)}_\infty\right).
\end{equation}
Here $\Theta$ denotes the tropical Riemann theta function
\begin{equation}\label{k:trt}
\Theta({\bf z}) = -\min_{{\bf n} \in {\mathbb Z}^G}
\left\{\frac{1}{2}{}^t{\bf n}B{\bf n}+ {}^t{\bf z}{\bf n}\right\},
\end{equation}
which enjoys the quasi-periodicity
\begin{equation}\label{k:qp}
\Theta({\bf z}+{\bf v}) = \Theta({\bf z})
+ {}^t{\bf v}B^{-1}\left({\bf z}+\frac{\bf v}{2}\right)\qquad
\text{ for } {\bf v} \in B{\mathbb Z}^G.
\end{equation}
In the context of the periodic box-ball system,
the tropical Riemann theta function was f\/irstly
obtained in this way in \cite{KS1} for rank $n=1$ case.
See \cite{IT0,IT,MZ} for an account from the
tropical geometry point of view.
Another remark is that Proposition \ref{k:pr:tbl} and (\ref{k:taut})
directly lead to the tropical Hirota equation for our $\Theta$:
\begin{gather*}
\Theta\big({\bf J}+{\bf h}^{(1)}_\infty+{\bf h}^{(d-1)}_\infty\big)+
\Theta\big({\bf J}+{\bf h}^{(1)}_1+{\bf h}^{(d)}_\infty\big)\\
\qquad{} =\max\bigl\{
\Theta\big({\bf J}+{\bf h}^{(1)}_\infty+{\bf h}^{(d)}_\infty\big)
+\Theta\big({\bf J}+{\bf h}^{(1)}_1+{\bf h}^{(d-1)}_\infty\big),\\
\qquad\quad \ \ \Theta\big({\bf J}+{\bf h}^{(1)}_1+{\bf h}^{(1)}_\infty
+{\bf h}^{(d-1)}_\infty\big)
+ \Theta\big({\bf J}+{\bf h}^{(d)}_\infty\big)-1\bigr\},
\end{gather*}
where we have set ${\bf J} = {\bf r}-\frac{\bf p}{2}-k{\bf h}^{(1)}_1$.
For $n=1$, see also \cite{IT}.

From (\ref{k:db}),  we arrive at our main formula in this section.

\begin{theorem}\label{k:th:xt}
The highest path
$x_1\otimes x_2 \otimes \cdots \otimes x_L
\in {\mathcal P}_+(\mu)$
corresponding to the rigged con\-figuration $(\mu,{\bf r})$ is given by
$x_k = (x_{k,1}, x_{k,2}, \ldots, x_{k,n+1})$ with
$x_{k,1} + \cdots + x_{k, n+1}=1$ and
\begin{gather}
x_{k,a}= \Theta\left({\bf r}-\frac{\bf p}{2}
-k{\bf h}^{(1)}_1+{\bf h}^{(a)}_\infty\right)
- \Theta\left({\bf r}-\frac{\bf p}{2}
-(k-1){\bf h}^{(1)}_1+{\bf h}^{(a)}_\infty\right)\label{k:xt4}\\
\phantom{x_{k,a}}{}- \Theta\left({\bf r}-\frac{\bf p}{2}
-k{\bf h}^{(1)}_1+{\bf h}^{(a-1)}_\infty\right)
+ \Theta\left({\bf r}-\frac{\bf p}{2}
-(k-1){\bf h}^{(1)}_1+{\bf h}^{(a-1)}_\infty\right)\quad
(2\le a \le n+1).\nonumber
\end{gather}
\end{theorem}

The data ${\bf r}$, ${\bf p}$ and ${\bf h}^{(c)}_l$
are determined from $(\mu,{\bf r})$ by
(\ref{k:rca0}), (\ref{k:pdef}), (\ref{k:rp}) and (\ref{k:hv}).

Compared with the previous expression (\ref{k:db}) with (\ref{k:tau1}),
the alternative formula (\ref{k:xt4}) enjoys an extra symmetry under the shifts
${\bf r} \rightarrow {\bf r} + B{\mathbb Z}^G$
and $k \rightarrow k+L{\mathbb Z}$ due to
the quasi-periodicity~(\ref{k:qp}) and (\ref{k:bh}).
In other words, (\ref{k:xt4})
is def\/ined for ${\bf r} \in {\mathbb Z}^G
/B{\mathbb Z}^G$ and $k \in {\mathbb Z}/L{\mathbb Z}$.
In view of this, it is natural to relate Theorem \ref{k:th:xt}
with our periodic $A^{(1)}_n$ SCA.

\begin{theorem}\label{k:con:ivp}
Assume the condition \eqref{k:pp}.
Let $p=x_1\otimes x_2 \otimes \cdots \otimes x_L
\in {\mathcal P}(\mu)$ be the path of the periodic $A^{(1)}_n$ SCA
whose action variable
is the quasi-periodic extension \eqref{k:qpr} of
${\bf r} =\big(r^{(a)}_{i,\alpha}\big)_{(a i \alpha)\in H}$.
Under Conjecture~{\rm \ref{t:conj:aug21_5}},
$x_k = (x_{k,1}, x_{k,2}, \ldots, x_{k,n+1})$
is given by \eqref{k:xt4}.
\end{theorem}

Note that this reduces a solution of the initial value problem
$p\rightarrow T^{(r_1)}_{l_1}\cdots T^{(r_N)}_{l_N}(p)$
to a simple substitution ${\bf r} \rightarrow
{\bf r} + {\bf h}^{(r_1)}_{l_1} + \cdots
+ {\bf h}^{(r_N)}_{l_N}$.

\begin{remark}\label{k:re:jb}
In addition to ${\bf r} \rightarrow {\bf r} + B{\mathbb Z}^G$,
each term in (\ref{k:xt4}) is invariant
under ${\bf r}=\big(r^{(a)}_{i,\alpha}\big)_{(ai\alpha)\in H} \mapsto
\big(r^{(a)}_{i,\sigma(\alpha)}\big)_{(ai\alpha)\in H}$
for any permutation $\sigma \in {\frak S}_{m^{(a)}_i}$ that can
depend on $(ai)\in \overline{H}$.
Here ${\frak S}_m$ denotes the symmetric group of degree $m$.
In fact the set ${\mathcal J}(\mu)$ of angle variables is naturally
described as
\begin{equation}\label{k:waru}
{\mathcal J}(\mu) \simeq {\mathbb Z}^G/B{\mathbb Z}^G/
\prod_{(ai)\in \overline{H}}{\frak S}_{m^{(a)}_i}.
\end{equation}
To go to the right hand side, one just forgets
the inequality $r^{(a)}_{i,\alpha}\le r^{(a)}_{i,\alpha+1}$
within each block $(ai) \in \overline{H}$ identifying all
the re-orderings. Compare also the matrix $B$ (\ref{k:B})
with (\ref{k:sai}).
\end{remark}

When $n=1$,
one can remove the assumption `Under Conjecture \ref{t:conj:aug21_5}'
in Theorem \ref{k:con:ivp} due to~\cite{KTT}.
If further $\forall\,  m^{(1)}_i=1$,
the formula (\ref{k:xt4}) coincides with \cite[equation~(3.8)]{KS1}.
Theorem~\ref{k:con:ivp} here tells that it also holds throughout $m^{(1)}_i>1$
by a natural extrapolation
whose concrete form is given in~(\ref{k:xsl2}) and~(\ref{k:bhp}).

The case $n=1$ and general $m^{(1)}_i$ in (\ref{k:xt4})
takes a dif\/ferent form from the corresponding formula in \cite{KS2}
expressed by a higher characteristic tropical Riemann theta function.
This apparent dif\/ference in guises
is caused by a dif\/ferent choices of ef\/fective
minimizing variables as $n^{(a)}_{i,\alpha}$ in (\ref{emv}).
In fact,  introducing the $M$ replicas
$N^{(a)}_{i,\alpha,1},\ldots, N^{(a)}_{i,\alpha,M}$ as in (\ref{k:tau2})
is not the unique way of handling  the
tropical tau function for the $M$-fold extended rigged conf\/iguration
$\big(\mu, {{\bf r}}^1\big)
\sqcup \big(\mu, {{\bf r}}^2\big) \sqcup \cdots
\sqcup \big(\mu, {{\bf r}}^{M}\big)$.
Another natural option is to stick to the form (\ref{k:tau1}) but
extend~$N^{(a)}_{i,\alpha}$ from $1 \le \alpha \le m^{(a)}_i$ to
$1\le \alpha \le Mm^{(a)}_i$ assuming (\ref{k:qpr}) for the rigging.
Then the minimizing variable $n^{(a)}_i$ can be introduced simply via
\begin{equation*}
N^{(a)}_{i,1}= N^{(a)}_{i,2}
= \cdots = N^{(a)}_{i, n^{(a)}_i}=1,
\qquad
N^{(a)}_{i, n^{(a)}_i+1}=\cdots =
N^{(a)}_{i, Mm^{(a)}_i}=0.
\end{equation*}
What makes the calculation tedious after this is that one has to classify
$n^{(a)}_i$ by $\text{mod}\, m^{(a)}_i$.
This was done in \cite{KS2} for $n=1$.
Here we omit the detail and only mention that
the result is expressed in terms of
a rational characteristic tropical Riemann theta function
with the $g\times g$ reduced period matrix:
\begin{equation*}
B^{\rm red} =
\text{diag}(m^{(a)}_i)_{(ai)\in \overline{H}}\,F
=\big(\delta_{ab}\delta_{ij}m^{(a)}_ip^{(a)}_i
+C_{ab}\min\big(l^{(a)}_i, l^{(b)}_j\big)m^{(a)}_im^{(b)}_j\big)_{
(ai), (bj) \in \overline{H}}.
\end{equation*}
For $n=1$, this coincides with  \cite[equation~(5.1)]{KS2}.

Let us proceed to an explicit $\Theta$ formula
for a carrier of type $B^{1,l}$.
Consider the highest path
$p=x_1\otimes x_2 \otimes \cdots \otimes x_L
\in {\mathcal P}_+(\mu)$ in Theorem \ref{k:th:xt}
expressed as (\ref{k:xt4}).
We consider the calculation of the time evolution
$T^{(1)}_l(p)=p'=x'_1\otimes \cdots \otimes x'_L$
according to (\ref{k:hone}).
Locally it is depicted as

\begin{picture}(50,60)(-170,-20)
\put(12,26){$x_k$}\put(59,26){$x_{k+1}$}
\put(20,0){\line(0,1){20}}\put(5,10){\line(1,0){30}}
\put(70,0){\line(0,1){20}}\put(85,10){\line(-1,0){30}}
\put(-20,8){$y_{k-1}$}\put(40,8){$y_k$}
\put(13,-13){$x_k'$}
\put(60,-13){$x_{k+1}'$}
\put(-45,10){$\ldots$}\put(97,10){$\ldots$}
\end{picture}

From the proof of Lemma \ref{k:le:ee},
a carrier satisfying
the periodic boundary condition
$y_0 = y_L \in B^{1,l}$ and thereby inducing the
time evolution $T^{(1)}_l$
can be constructed by
$u^{1,l}\otimes p \simeq \hat{p} \otimes y_0$.
Here $u^{1,l}\in B^{1,l}$ is the
semistandard tableau of length $l$ row shape
whose entries are all $1$.
This f\/ixes the carriers in the intermediate
stage $y_1, y_2, \ldots, y_{L-1} \in B^{1,l}$ by
$y_0\otimes (x_1\otimes \cdots \otimes x_k)
\overset{\sim}{\mapsto}
(x'_1\otimes \cdots \otimes x'_k)\otimes y_k$
under the isomorphism
$B^{1,l}\otimes \big(B^{1,1}\big)^{\otimes k}
\simeq \big(B^{1,1}\big)^{\otimes k}\otimes B^{1,l}$.

\begin{theorem}\label{th:car}
The carrier $y_k \in B^{1,l}$ corresponding to (\ref{k:xt4})
in the above sense is given by
$y_k = (y_{k,1}, y_{k,2}, \ldots, y_{k,n+1})$ with
$y_{k,1} + \cdots + y_{k, n+1}=l$ and
\begin{gather}
y_{k,a}  = \Theta\left({\bf r}-\frac{\bf p}{2}
-k{\bf h}^{(1)}_1+{\bf h}^{(a)}_\infty\right)
- \Theta\left({\bf r}-\frac{\bf p}{2}
-k{\bf h}^{(1)}_1+{\bf h}^{(1)}_l+{\bf h}^{(a)}_\infty\right)\label{carn}\\
\phantom{y_{k,a}}{}- \Theta\left({\bf r}-\frac{\bf p}{2}
-k{\bf h}^{(1)}_1+{\bf h}^{(a-1)}_\infty\right)
+ \Theta\left({\bf r}-\frac{\bf p}{2}
-k{\bf h}^{(1)}_1+{\bf h}^{(1)}_l+{\bf h}^{(a-1)}_\infty\right)\quad\!\!
(2\le a \le n+1).\nonumber
\end{gather}
\end{theorem}

This can be verif\/ied by modifying the derivation of
Theorem~\ref{k:th:xt} slightly by using
\cite[Theo\-rem~2.1]{KSY} and
\cite[Lemma~8.5]{KSS}.
The periodicity $y_0=y_L$ is easily checked by (\ref{k:qp}).
Naturally~(\ref{carn}) gives the
carrier for general paths in Theorem \ref{k:con:ivp}
under Conjecture \ref{t:conj:aug21_5}.
For $n=1$ such a result for a higher spin case was obtained in \cite{KS3}
under the condition $\forall\,  m^{(1)}_i=1$.
We will present a~modest application of
the formula (\ref{carn}) in Section \ref{k:sec:mis}.

\subsection[Bethe vector at $q=0$ from tropical Riemann theta function]{Bethe vector at $\boldsymbol{q=0}$ from tropical Riemann theta function}
\label{k:sec:bv}

We assume the condition (\ref{k:pp})
and Conjecture \ref{t:conj:aug21_5} in this subsection.
Write the string center equation (\ref{eq:sce2}) simply as
\begin{equation}\label{k:sce3}
A{\bf u} = {\bf c}+{\bf r}+\boldsymbol{\rho}
\end{equation}
using the $G$-dimensional vectors
${\bf c} = \big(\frac{1}{2}
\big(p^{(a)}_i + m^{(a)}_i + 1\big)\big)_{(ai\alpha) \in H}$ and
$\boldsymbol{\rho} = (\alpha-1)_{(ai\alpha) \in H}$.
Note that the time evolution of the angle variable
(\ref{k:ter2}) is written as
$T^{(r)}_l({\bf r}) = {\bf r}+{\bf h}^{(r)}_l$
with ${\bf h}^{(r)}_l$ def\/ined by (\ref{k:hv}).
In view of (\ref{k:sce3}),
the time evolution of the Bethe roots introduced in
(\ref{eq:cd5}) is expressed as
$T^{(r)}_l({\bf u}) = {\bf u} + A^{-1}{\bf h}^{(r)}_l$.
At f\/irst sight, this appears contradictory, because
$T^{(r)}_l$ is a transfer matrix at $q=0$, which
should leave the $q=0$ Bethe eigenvectors invariant
up to an overall scalar hence the relevant Bethe roots as well.
The answer to this puzzle is that the path
$p\in {\mathcal P}(\mu)$ that we are associating
with ${\bf u}$ or ${\bf r}$ by
$p = \Phi^{-1}(\mu, {\bf r})$
is a {\em monomial} in
$(\mathbb{C}^{n+1})^{\otimes L}$,
which is {\em not} a Bethe vector at $q=0$ in general.

It is easy to remedy this.
In fact, for each Bethe root ${\bf u}$ or equivalently
the angle variable
${\bf r}=A{\bf u} -{\bf c}-\boldsymbol{\rho}
\in {\mathcal J}(\mu)$,
one can construct a vector
$\vert {\bf r} \rangle  \in ({\mathbb C}^{n+1})^{\otimes L}$
that possesses every aspect as a $q=0$ Bethe vector as follows:
\begin{gather}
\vert {\bf r}\rangle
= \sum_{{\bf s} \in {\mathcal J}(\mu)}
c_{{\bf s}, {\bf r}}\, {\bf v}({\bf s}),\label{eq:jp}\\
c_{{\bf s}, {\bf r}}  = \exp\left(-2\pi\sqrt{-1}\;
{}^t{\bf s}\left(A^{-1}({\bf r}+{\bf c}+\boldsymbol{\rho})
+\frac{{\bf 1}}{2}\right)\right), \label{k:c}\\
{\bf v}({\bf s})   =
\begin{pmatrix}x_{1,1}({\bf s})\\ \vdots \\
x_{1,n+1}({\bf s}) \end{pmatrix}
\otimes \cdots \otimes
\begin{pmatrix}x_{L, 1}({\bf s})\\ \vdots \\
x_{L, n+1}({\bf s}) \end{pmatrix}
\in {\mathcal P}(\mu) \subseteq
\big({\mathbb C}^{n+1}\big)^{\otimes L}. \nonumber
\end{gather}
Here ${\bf 1}=(1)_{(ai\alpha)\in H}$ in (\ref{k:c}),
which has also appeared in (\ref{k:bh}).
The component $x_{k,a}({\bf s}) \in \{0,1\}$
is specif\/ied by Theorem \ref{k:con:ivp}.
Namely,
\begin{gather*}
x_{k,a}({\bf s})  = \Theta\left({\bf s}-\frac{\bf p}{2}
-k{\bf h}^{(1)}_1+{\bf h}^{(a)}_\infty\right)
- \Theta\left({\bf s}-\frac{\bf p}{2}
-(k-1){\bf h}^{(1)}_1+{\bf h}^{(a)}_\infty\right)\\
\phantom{x_{k,a}({\bf s})  =}{} - \Theta\left({\bf s}-\frac{\bf p}{2}
-k{\bf h}^{(1)}_1+{\bf h}^{(a-1)}_\infty\right)
+ \Theta\left({\bf s}-\frac{\bf p}{2}
-(k-1){\bf h}^{(1)}_1+{\bf h}^{(a-1)}_\infty\right)
\end{gather*}
for $2\le a \le n+1$
and $x_{k,1}({\bf s}) + \cdots + x_{k,n+1}({\bf s})=1$.
The vectors ${\bf p}$ and ${\bf h}^{(a)}_l$ are def\/ined in
(\ref{k:rp}) and (\ref{k:hv}).
We embed ${\mathcal P}(\mu)$
into $({\mathbb C}^{n+1})^{\otimes L}$
naturally and extend the time evolutions
to the latter by ${\mathbb C}$-linearity.
The vector ${\bf v}({\bf s})$ here is the
path corresponding to the angle variable
${\bf s} \in {\mathcal J}(\mu)$
in Theorem \ref{k:con:ivp}.
It follows that
$T^{(r)}_l({\bf v}({\bf s}))
= {\bf v}\big({\bf s}+{\bf h}^{(r)}_l\big)$.
Then using
${\mathcal J}(\mu) + {\bf h}^{(r)}_l
= {\mathcal J}(\mu)$,
it is elementary to check
\begin{gather}\label{k:eie}
T^{(r)}_l(\,\vert {\bf r}\rangle\,)
= c_{-{\bf h}^{(r)}_l, \,{\bf r}}\,\vert {\bf r}\rangle,\\
c_{-{\bf h}_l^{(r)}, \,{\bf r}}
= \exp\left(2\pi\sqrt{-1}\;
{}^t{\bf h}^{(r)}_l
\left(A^{-1}({\bf r}+{\bf c}+\boldsymbol{\rho})
+\frac{{\bf 1}}{2}\right)\right)
= \exp\left(2\pi\sqrt{-1}\;
{}^t{\bf h}^{(r)}_l
\left({\bf u}+\frac{{\bf 1}}{2}\right)\right).\nonumber
\end{gather}
The last expression of the
quantity $c_{-{\bf h}_l^{(r)}, \,{\bf r}}$
reproduces the
$q=0$ Bethe eigenvalue
$\Lambda^{(r)}_l
=\prod_{j\alpha}\big(-z^{(r)}_{j\alpha}\big)^{\min(j,l)}$
given in \cite[equation~(3)]{KT2}
with the string center $z^{(r)}_{j\alpha}
= \exp\big(2\pi\sqrt{-1}u^{(r)}_{j\alpha}\big)$.
Note further that the transition relation (\ref{eq:jp}) is inverted as
\begin{equation*}
{\bf v}({\bf s}) = \frac{1}{\vert {\mathcal J}(\mu) \vert}
\sum_{{\bf r} \in {\mathcal J}}
\overline{c}_{{\bf s}, {\bf r}}
\vert {\bf r}\rangle,
\end{equation*}
where $\overline{c}_{{\bf s}, {\bf r}}$
denotes the complex conjugate of $c_{{\bf s}, {\bf r}}$.
It follows that the space of the `$q=0$ Bethe vectors'
$\vert {\bf r} \rangle$
coincides with the space spanned
by the evolvable paths of the prescribed soliton content $\mu$.
\begin{equation}\label{k:ss}
\bigoplus_{{\bf r}  \in {\mathcal J}(\mu)}{\mathbb C}
\vert {\bf r} \rangle
= \bigoplus_{p \in {\mathcal P}(\mu)} {\mathbb C} \,p.
\end{equation}
Thus the approach here bypasses
a formidable task of computing the $q\rightarrow 0$ limit of the
Bethe vectors directly, and leads to the joint eigenvectors
$\vert {\bf r} \rangle$ of the commuting family of
time evolutions.
They form a basis of the space having the
prescribed soliton content and possess
the eigenvalues anticipated
from the Bethe ansatz at $q=0$.
These features together with Theorem~\ref{k:th:rc}
constitute the quantitative background of the identity~(\ref{k:sost}).

\subsection[Bethe eigenvalue at $q=0$ and dynamical period]{Bethe eigenvalue at $\boldsymbol{q=0}$ and dynamical period}
\label{k:sec:bedp}

Let us remark on the relation of the
Bethe eigenvalue (\ref{k:eie}) to the dynamical period.
In \cite{KT1,KT2}, the formula
(\ref{t:eq:aug18_3}) for the dynamical period of periodic $A^{(1)}_n$
SCA with $n>1$ without the order of symmetry $\boldsymbol{\gamma}$
was found by demanding the condition
\begin{equation}\label{k:njou}
\bigl(c_{-{\bf h}_l^{(r)}, \,{\bf r}}\,
\bigr)^{{\mathcal N}^{(r)\prime}_l} = \pm 1.
\end{equation}
In fact, from the middle expression in (\ref{k:eie})
this condition is satisf\/ied if
\begin{equation*}
{\mathcal N}^{(r)\prime}_l {\bf h}^{(r)}_l =0 \mod A{\mathbb Z}^G.
\end{equation*}
By the same argument as the proof of Theorem \ref{k:th:dp},
one f\/inds that the smallest positive integer
satisfying this condition is
\begin{equation*}
{\mathcal N}^{(r)\prime}_l
= {\rm LCM}\!
\left(\frac{\det A}
{\det A[bj\beta]}\right)_{(bj\beta) \in H}
= {\rm LCM}\!
\left(\frac{\det{F}}{\det{F} [bj]}
\right)_{(bj) \in \overline{H}},
\end{equation*}
where $A[bj\beta]$ denotes the $G\times G$
matrix obtained from $A$ (\ref{eq:a})
by replacing its $(bj\beta)$th column
by ${\bf h}^{(r)}_l$ (\ref{k:hv}).
The $\text{LCM}$ in the f\/irst (second) expression
should be taken over only those
$(bj\beta)$ such that
$\det A[bj\beta] \neq 0$
($(bj)$ such that $\det F[bj] \ne 0$).
The second equality is due to the
structures of matrices $A$ and
$F$ (\ref{k:f}) and verif\/ied by a direct calculation.
The expression ${\mathcal N}^{(r)\prime}_l $
derived in the heuristic approach \cite{KT1,KT2}
captures the main structure of the full formula (\ref{t:eq:aug18_3}).
However it neither f\/ixes the sign in (\ref{k:njou}) nor takes
the order of symmetry $\boldsymbol{\gamma}$ into account.

These shortcoming are f\/ixed of course by
ref\/ining the construction of joint eigenvectors
along the connected components.
Instead of trying to split the sum (\ref{eq:jp}) into them,
we simply introduce
\begin{equation}\label{k:obr}
|\omega, \boldsymbol{\lambda}\rangle
=\sum_{\phi \in {\mathbb Z}^g/F_{\boldsymbol{\gamma}}{\mathbb Z}^g}
\exp\big(-2\pi\sqrt{-1}\,{}^t\phi\, F_{\boldsymbol{\gamma}}^{-1}
(\omega + d_{\boldsymbol{\lambda}})\big)
{\bf v}(\phi,\boldsymbol{\lambda})
\in \big({\mathbb C}^{n+1}\big)^{\otimes L}
\end{equation}
for each angle variable $(\omega, \boldsymbol{\lambda})
\in X_{\boldsymbol{\gamma}}/{\mathcal A}$
having the order of symmetry $\boldsymbol{\gamma}$.
Here ${\bf v}(\phi,\boldsymbol{\lambda})
\in {\mathcal P}_{\boldsymbol{\gamma}}(\mu)$ is the path
corresponding to $(\phi,\boldsymbol{\lambda})
\in X_{\boldsymbol{\gamma}}/{\mathcal A}$, and
$d_{\boldsymbol{\lambda}} \in {\mathbb Z}^g$ can be chosen arbitrarily.
Then by noting $T^{(r)}_l {\bf v}(\phi,\boldsymbol{\lambda})
={\bf v}(\phi+h^{(r)}_l, \boldsymbol{\lambda})$,
it is straightforward to check
\begin{equation*}
T^{(r)}_l(\, |\omega, \boldsymbol{\lambda}\rangle\,) =
\exp\big(2\pi\sqrt{-1}\,{}^th^{(r)}_l\, F_{\boldsymbol{\gamma}}^{-1}
(\omega + d_{\boldsymbol{\lambda}})\big)
 |\omega, \boldsymbol{\lambda}\rangle,
\end{equation*}
This eigenvalue is indeed an ${\mathcal N}^{(r)}_l$th root
of unity due to the proof of Theorem \ref{k:th:dp}.

Given $\boldsymbol{\lambda}$,
there are $\det F_{\boldsymbol{\gamma}}$
independent vectors $|\omega, \boldsymbol{\lambda}\rangle$
(\ref{k:obr}).
On the other hand, the number of the choices of
$\boldsymbol{\lambda}
\in X^2_{\boldsymbol{\gamma}}/{\mathcal A}$
is given by Lemma \ref{k:le:xa}.
Obviously these vectors are all independent
because the set of monomials involved in
$|\omega, \boldsymbol{\lambda}\rangle$ and
$|\omega', \boldsymbol{\lambda}'\rangle$ are
the set of paths that are
disjoint if $\boldsymbol{\lambda} \neq \boldsymbol{\lambda}'$.
{}From this fact and (\ref{k:kore}), we obtain the ref\/inement of
(\ref{k:ss}) according to the order of symmetry $\boldsymbol{\gamma}$
and further (bit tautologically) according to the
connected components:
\begin{equation*}
\bigoplus_{(\omega, \boldsymbol{\lambda})
\in X_{\boldsymbol{\gamma}}/{\mathcal A}}{\mathbb C}
|\omega, \boldsymbol{\lambda}\rangle
= \bigoplus_{p \in {\mathcal P}_{\boldsymbol{\gamma}}(\mu)}
 {\mathbb C} \,p,\qquad
\bigoplus_{\omega \in {\mathbb Z}^g/F_{\boldsymbol{\gamma}}{\mathbb Z}^g}{\mathbb C}
|\omega, \boldsymbol{\lambda}\rangle
= \bigoplus_{p \in \Sigma(p_0)}{\mathbb C} \,p,
\end{equation*}
where for example
$p_0 = \Phi^{-1}((\omega = 0, \boldsymbol{\lambda}))
\in {\mathcal P}_{\boldsymbol{\gamma}}(\mu)$.

\subsection{Miscellaneous
calculation of time average}\label{k:sec:mis}

As a modest application of the formulas
by the tropical Riemann theta function,
we f\/irst illustrate a calculation of some time average
along the periodic box-ball system ($n=1$).
Analogous results will be stated for general $n$ in the end.
We use the terminology in the periodic box-ball system.

Let $p \in (B^{1,1})^{\otimes L}$ be the path
with the angle variable ${\bf I}$.
The time evolution $T^{(1)}_l$ will simply be denoted by $T_l$.
Set $T^t_l(p) = x_1(t)\otimes \cdots \otimes x_L(t)$.
For $n=1$, one can label $x_k(t)=(x_{k,1}(t), x_{k,2}(t)) \in B^{1,1}$
just by $x_{k,2}(t)$, which from now on will simply be denoted by
$x_k(t) (=0,1)$.
It represents the number of balls in the $k$th box.
Then (\ref{k:xt4}) reads
\begin{gather}
x_k(t) = \Theta({\bf J}-k{\bf h}_1+t{\bf h}_l)
-\Theta({\bf J}-(k-1){\bf h}_1+t{\bf h}_l)\nonumber\\
\phantom{x_k(t) =}{} -\Theta({\bf J}-k{\bf h}_1+t{\bf h}_l+{\bf h}_\infty)
+\Theta({\bf J}-(k-1){\bf h}_1+t{\bf h}_l+{\bf h}_\infty)\label{k:xsl2}
\end{gather}
with ${\bf J} = {\bf I}-\frac{\bf p}{2}$.
The notations
(\ref{k:B})--(\ref{k:hv}) are simplif\/ied hereafter as
\begin{gather}
B  = (\delta_{ij}p_i + 2\min(l_i,l_j))_{(i\alpha), (j\beta) \in H},
\qquad
H=\{(i, \alpha)\mid 1 \le i \le g, 1 \le \alpha \le m_i\},\nonumber\\
{\bf p}  =(p_i)_{(i\alpha) \in H},\qquad
p_i = L-2\sum_{j=1}^g\min(l_i,l_j)m_j,\qquad
{\bf h}_l = (\min(l,l_i))_{(i\alpha) \in H}.\label{k:bhp}
\end{gather}
Here $l_i$, $m_i$ are the shorthand of  $l^{(1)}_i$, $m^{(1)}_i$
that specify the action variable (single partition)
$\mu=\mu^{(1)}$ as in (\ref{k:rca0}).
The relation (\ref{k:bh}) reads
\begin{equation}\label{k:bh2}
B{\bf h}_1= L{\bf h}_1.
\end{equation}
Note that the total number of balls at any time is
$|\mu | = \sum_{i=1}^gl_im_i$.
This also follows immediately from (\ref{k:xsl2}) as
\begin{gather*}
\sum_{k=1}^L x_k(t)  =
\Theta({\bf J}-L{\bf h}_1+t{\bf h}_l)
-\Theta({\bf J}+t{\bf h}_l)
-\Theta({\bf J}-L{\bf h}_1+t{\bf h}_l+{\bf h}_\infty)
+\Theta({\bf J}+t{\bf h}_l+{\bf h}_\infty)\\
\hphantom{\sum_{k=1}^L x_k(t)}{} =-L{}^t{\bf h}_1B^{-1}
\left({\bf J}+t{\bf h}_l-\frac{L{\bf h}_1}{2}\right)
+L{}^t{\bf h}_1B^{-1}
\left({\bf J}+t{\bf h}_l+{\bf h}_\infty-\frac{L{\bf h}_1}{2}\right)\\
\hphantom{\sum_{k=1}^L x_k(t)}{}
=L{}^t{\bf h}_1B^{-1}{\bf h}_\infty={}^t{\bf h}_1{\bf h}_\infty
=|\mu|,
\end{gather*}
where we have used (\ref{k:bh2}) and the
quasi-periodicity (\ref{k:qp}).
We introduce the density of balls:
\begin{equation*}
\rho = {}^t{\bf h}_1B^{-1}{\bf h}_\infty = \frac{|\mu|}{L}.
\end{equation*}

Let ${\mathcal N}_l$ be the dynamical period under $T_l$.
Thus $x_k({\mathcal N}_l)=x_k(0)$ holds for any $1\le k \le L$.

\begin{proposition}\label{k:pr:av1}
The time average of $x_k(t)$ under $T_\infty$
over the period ${\mathcal N}_\infty$ is given by
\begin{gather}\label{k:baka}
\frac{1}{{\mathcal N}_\infty}\sum_{t=0}^{{\mathcal N}_\infty-1}x_k(t) = \rho.
\end{gather}
\end{proposition}
\begin{proof}
Using (\ref{k:xsl2}) with $l=\infty$,
we f\/ind that $\sum_{t=0}^{{\mathcal N}_\infty-1}x_k(t)$
is equal to
\begin{equation*}
\Theta({\bf J}-k{\bf h}_1)
-\Theta({\bf J}-k{\bf h}_1+{\mathcal N}_\infty{\bf h}_\infty)
-\Theta({\bf J}-(k-1){\bf h}_1)
+\Theta({\bf J}-(k-1){\bf h}_1+{\mathcal N}_\infty{\bf h}_\infty).
\end{equation*}
Note that the assumption implies
${\mathcal N}_\infty{\bf h}_\infty \in B{\mathbb Z}^G$.
Thus one can reduce this by applying the quasi-periodicity
(\ref{k:qp}), obtaining
${\mathcal N}_\infty{}^t{\bf h}_\infty B^{-1}{\bf h}_1={\mathcal N}_\infty\rho$.
\end{proof}
Actually (\ref{k:baka}) follows at once without this sort of calculation
if the left hand side is assumed to be independent of $k$.
However, such homogeneity under spatial translation is not always valid
for some time evolution $T_l$ with $l<\infty$
and the initial condition that possess a special commensurability.

Denote the time average
$\frac{1}{{\mathcal N}_l}\sum_{t=0}^{{\mathcal N}_l-1}Q(t)$
of a quantity $Q(t)$ under $T_l$
by $\langle Q \rangle_l$.
Then a trivial corollary of Proposition \ref{k:pr:av1} is
\begin{equation*}
\Bigl\langle \sum_{k=1}^L\omega^{k-1}x_k \Bigr\rangle_\infty
= \begin{cases} \rho L & \omega=1,\\
\rho\frac{1-\omega^L}{1-\omega} & \omega \neq 1.
\end{cases}
\end{equation*}
In particular, this average vanishes if $\omega$ is a nontrivial $L$th root of unity.

Let us proceed to a less trivial example.
In~\cite{KS3}, the number of balls in the carrier
for any time evolution $T_l$ was expressed also
in terms of the tropical Riemann theta function.
The carrier that induces the time evolution $p \mapsto T_l(p)$
is the element $v_l \in B^{1,l}$ such that
$v_l \otimes p \simeq T_l(p)\otimes v_l$
under the isomorphism of crystals
$B^{1,l}\otimes \big(B^{1,1}\big)^{\otimes L}
\simeq \big(B^{1,1}\big)^{\otimes L} \otimes B^{1,l}$.
See the explanation around~(\ref{k:hone}).
There uniquely exists such $v_l$ for any $p$ with
density $\rho<1/2$ as shown in
\cite[Proposition~2.1]{KTT}.
Let us consider the intermediate stage of sending the
carrier $y_0(t):=v_l$ to the right by repeated applications of
combinatorial $R$:
\begin{gather*}
v_l  \otimes    x_1(t)  \otimes   \cdots   \otimes  x_k(t)
 \otimes x_{k+1}(t)
\cdots \otimes x_L(t)\\
\qquad{}\simeq
x_1(t+1)\otimes \cdots \otimes x_k(t+1)\otimes
y_k(t)  \otimes x_{k+1}(t)
\cdots \otimes x_L(t).
\end{gather*}
This is depicted locally as

\begin{picture}(50,60)(-170,-20)
\put(12,26){$x_k(t)$}\put(59,26){$x_{k+1}(t)$}
\put(20,0){\line(0,1){20}}\put(0,10){\line(1,0){30}}
\put(70,0){\line(0,1){20}}\put(90,10){\line(-1,0){30}}
\put(35,8){$y_k(t)$}
\put(-5,-13){$x_k(t+1)$}
\put(48,-13){$x_{k+1}(t+1)$}
\put(-22,10){$\ldots$}\put(100,10){$\ldots$}
\end{picture}

We identify $y_k(t)=(y_{k,1}(t), y_{k,2}(t))\in B^{1,l}$
with the number of balls in the capacity $l$ carrier.
Namely, $y_{k,2}(t) \in \{0,1,\ldots, l\}$
will simply be denoted by $y_k(t)$.
(Hence $y_{k,1}(t) = l-y_k(t)$.)
Then for $x_k(t)$ given in (\ref{k:xsl2}),  one has \cite{KS3}
\begin{gather*}
y_k(t)  = \Theta({\bf J}-k{\bf h}_1+t{\bf h}_l)
-\Theta({\bf J}-k{\bf h}_1+(t+1){\bf h}_l)\\
\phantom{y_k(t)  =}{} -\Theta({\bf J}-k{\bf h}_1+t{\bf h}_l+{\bf h}_\infty)
+\Theta({\bf J}-k{\bf h}_1+(t+1){\bf h}_l+{\bf h}_\infty).
\end{gather*}
The preceding result (\ref{carn}) reduces to this upon
setting $n=1$, $a=2$ and
${\bf r}-\frac{\bf p}{2} = {\bf J} + t{\bf h}_l$.

\begin{proposition}\label{k:pr:av2}
The time average of the number of balls in the carrier with
capacity $l$ at position $k$ under $T_l$ is given by
\begin{equation*}
\langle y_k \rangle_l = {}^t{\bf h}_lB^{-1}{\bf h}_\infty.
\end{equation*}
\end{proposition}

\begin{proof}
The proof is parallel with the one for Proposition \ref{k:pr:av1}
by using ${\mathcal N}_l{\bf h}_l \in B{\mathbb Z}^G$
and the quasi-periodicity (\ref{k:qp}).
\end{proof}

Note that the result is independent
of the order of symmetry $\boldsymbol{\gamma}$
as well as $k$.

\begin{example}\label{k:ex:last}
Take the path
$p=111222221111222222221111111111111112222111211$
of length $L=45$ from \cite[Example 5.2]{KS2}.
The action variable is
\begin{equation*}
\begin{picture}(100,56)(-20,-10)
\put(0,40){\line(1,0){90}}
\put(0,30){\line(1,0){90}}
\put(0,20){\line(1,0){40}}
\put(0,10){\line(1,0){40}}
\put(0,0){\line(1,0){10}}

\put(0,0){\line(0,1){40}}
\put(10,0){\line(0,1){40}}
\put(20,10){\line(0,1){30}}
\put(30,10){\line(0,1){30}}
\put(40,10){\line(0,1){30}}
\put(50,30){\line(0,1){10}}
\put(60,30){\line(0,1){10}}
\put(70,30){\line(0,1){10}}
\put(80,30){\line(0,1){10}}
\put(90,30){\line(0,1){10}}

\end{picture}
\end{equation*}
The data from (\ref{k:bhp}) read
\begin{equation*}
B = \begin{pmatrix}
 27 & 8 & 8 & 2 \\
 8 & 27 & 8 & 2 \\
 8 & 8 & 27 & 2 \\
 2 & 2 & 2 & 39
\end{pmatrix},\qquad
{\bf p} = \begin{pmatrix}
9 \\
19 \\
19\\
37
\end{pmatrix},\qquad
{\bf h}_l = \begin{pmatrix}
\min(l,9)\\
\min(l,4)\\
\min(l,4)\\
\min(l,1)
\end{pmatrix}.
\end{equation*}
We list average $\langle y_k \rangle_l$ with the dynamical period
${\mathcal N}_l$.
\begin{center}
\begin{tabular}{c|c|c|c|c|c|c|c|c|c}
$l$ & 1 & 2& 3& 4& 5& 6& 7& 8& $\infty$ \\
\hline
$\overset{\phantom{}}{{\mathcal N}_l}$
& 45 & 1665 & 333 & 1665 & 1665 & 31635
& 3515 & 31635 & 6327 \\\hline
$\langle y_k \rangle_l$ &
$\overset{\phantom{}}{\rho=\frac{2}{5}=0.4}$
& $\frac{147}{185}$ & $\frac{44}{37}$
& $\frac{293}{185}$     & $\frac{6646}{3515}$
& $\frac{1545}{703}$   & $\frac{8804}{3515}$
& $\frac{9883}{3515}$ & $\frac{10962}{3515}=3.11863\dots$
\end{tabular}
\end{center}
\end{example}

Now it requires little explanation
to present the analogous result for $A^{(1)}_n$ case.
Denote the time average
$\frac{1}{\mathcal N^{(1)}_l}
\sum_{t=0}^{{\mathcal N}^{(1)}_l-1}Q(t)$
of a quantity $Q(t)$ under the time evolution $T^{(1)}_l$
by $\langle Q \rangle_l$.
Here ${\mathcal N}^{(1)}_l$ is the dynamical period
under $T^{(1)}_l$.
Then, the same calculation as
Proposition \ref{k:pr:av2} using~(\ref{carn}) leads to

\begin{proposition}
Under Conjecture {\rm \ref{t:conj:aug21_5}},
the time average of the number of balls of color $a$ in the carrier
of type $B^{1,l}$ at position $k$ is given by
\begin{equation*}
\langle y_{k,a} \rangle_l
= {}^t{\bf h}^{(1)}_lB^{-1}
\big({\bf h}^{(a-1)}_\infty - {\bf h}^{(a)}_\infty\big)\qquad
(2 \le a \le n+1).
\end{equation*}
\end{proposition}
Again this depends neither on $k$ nor
on the order of symmetry $\boldsymbol{\gamma}$.
It reduces to Proposition~\ref{k:pr:av2} for $n=1$
due to ${\bf h}^{(n+1)}_l = 0$.
See~(\ref{k:hv}).
In particular at $l=1$ we have
\begin{equation*}
\langle y_{k,a} \rangle_1
= \frac{1}{L}
{}^t{\bf 1}\big({\bf h}^{(a-1)}_\infty - {\bf h}^{(a)}_\infty\big)
=\frac{|\mu^{(a-1)}| - |\mu^{(a)}|}{L}
\end{equation*}
by means of (\ref{k:bh}) and (\ref{new}).
From (\ref{k:an}), this is the density of the color $a$ balls
in the path $\big(B^{1,1}\big)^{\otimes L}$, which is
natural in view of (\ref{k:shf}).
The average of the total number of balls is
$\sum_{a=2}^{n+1}\langle y_{k,a} \rangle_l
= {}^t{\bf h}^{(1)}_lB^{-1}{\bf h}^{(1)}_\infty$,
which again reduces to Proposition \ref{k:pr:av2} for $n=1$.

\begin{example}\label{ex:owari}
Consider the paths in Example \ref{k:ex:int} ($n=2$).
The conf\/iguration $\mu$ is given in
Example \ref{k:ex:aa2}, which is the same as the one in
Example \ref{k:ex:aa}.
The matrix $B$ and ${\bf h}^{(a)}_l$
are listed in Example \ref{k:ex:Brp}, which says
${\bf h}^{(1)}_{3}={\bf h}^{(1)}_\infty$.
Denote the density by
$\rho_a = \frac{|\mu^{(a-1)}| - |\mu^{(a)}|}{L}$.
The average $\langle y_{k,a} \rangle_l$ is
shown in the following table with the dynamical period
given in the end of Example~\ref{k:ex:int}.
\begin{center}
\begin{tabular}{c|c|c|c}
$l$ & 1 & 2& $\infty$ \\
\hline
${\mathcal N}^{(1)}_l$ & 24 & 12 & 194\\
\hline
$\langle y_{k,2} \rangle_l$ &
$\underset{\phantom{}}
{\overset{\phantom{}}{\rho_2=\frac{7}{24}}
=0.291\dots}$
& ${\frac{7}{12}}=0.583\dots$
& $\frac{155}{194} =0.798\dots$\vspace{-0.1cm}\\
\hline
$\langle y_{k,3} \rangle_l$ &
$\overset{\phantom{}}{\rho_3=\frac{5}{24}=0.208\dots}$
& $\frac{5}{12}=0.416\dots$ & $\frac{109}{194}=0.561\dots$
\end{tabular}
\end{center}
\end{example}

\appendix

\section{Row and column insertions}\label{k:app:ins}

Let $T$ be a semistandard tableau.
The row insertion of the number $r$ into $T$ is
denoted by $T \leftarrow r$ and def\/ined recursively as follows:
$$
\includegraphics{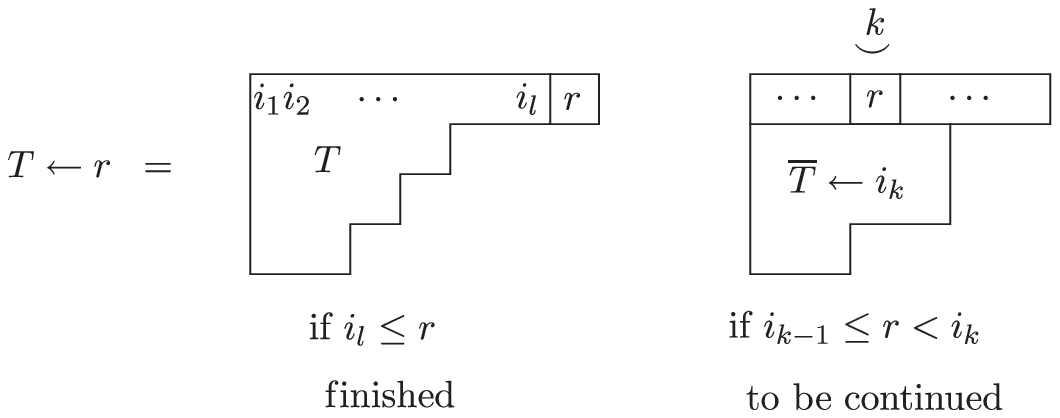}
$$
Here we have denoted the entries in the f\/irst row of
$T$ by $i_1\le  i_2 \le \cdots \le i_l$, and the other part of $T$ by
$\overline{T}$.
The f\/inished case includes the situation $T = \varnothing$.
In the to be continued case, $r$ bumps out
the smallest number $i_k$ that is larger than $r$
(row bumping).
The tableau $T\leftarrow r$ is obtained by repeating the row bumping
until f\/inished.

The column insertion of the number $r$ into $T$ is
denoted by $r \rightarrow T$
and def\/ined recursively as follows:
$$
\includegraphics{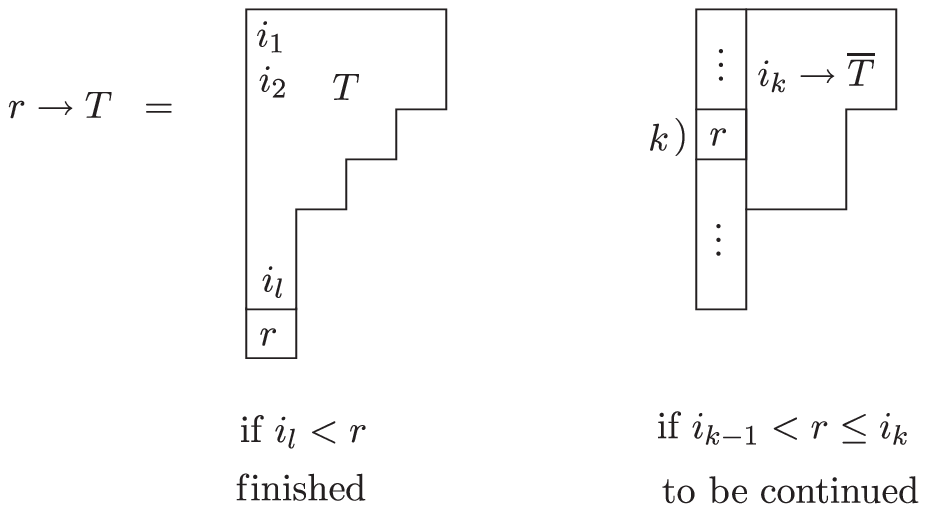}
$$
Here we have denoted the entries of the f\/irst column of $T$ by
$i_1<i_2< \cdots < i_l$, and the other part of $T$ by $\overline{T}$.
The f\/inished case includes the situation $T = \varnothing$.
In the to be continued case, $r$ bumps out
the smallest number $i_k$ that is not less than $r$ (column bumping).
The tableau
$r \rightarrow T$ is obtained by repeating the column bumping
until f\/inished.

\section{Proof of Theorem \ref{k:th:fac}}\label{k:app:fac}

In this appendix, we use the notions from crystal theory.
See \cite{K, KMN2,NY} and \cite[Section 2.3]{KTT} for the notations
$\varphi_i, \varepsilon_i$, the Kashiwara operator $\tilde{e}_i$
and the Weyl group simple ref\/lection $S_i$.
The crystal $B^{1,l}$ can  be parameterized as
$B^{1,l}= \{u=(u_1,\ldots, u_{n+1})\in ({\mathbb Z}_{\ge 0})^{n+1} \mid
u_1 + \cdots + u_{n+1} = l \}$,
where $u_i$ denotes the number of letter $i$ contained in a
length $l$ row shape semistandard tableau $u$ \cite{NY,HKT}.
The cyclic automorphism $\sigma$ acts as
$\sigma((u_1,\ldots, u_{n+1})) = (u_2,\ldots, u_{n+1},u_1)$,
or equivalently $u_i \mapsto u_{i+1}$.
The automorphism $\sigma$ acts on a tensor product component-wise.
The indices are to be considered in ${\mathbb Z}/(n+1){\mathbb Z}$.
Recall that $B=(B^{1,1})^{\otimes L}$ and $B_1$ is def\/ined in
(\ref{k:B1}).
We invoke the following factorization property
of the combinatorial $R$ into simple ref\/lections.

\begin{theorem}[\protect{\cite[Theorem 2]{HKT}}]
For any path $p \in B$ and
$u=(u_1, u_2,\ldots, u_{n+1}) \in B^{1,l}$ with
$u_1 \gg u_2,\ldots, u_{n+1}$, let
$B^{1,l} \otimes B \ni u \otimes p
\simeq {\tilde p} \otimes {\tilde u} \in B \otimes B^{1,l}$
be the image of the isomorphism.
Then ${\tilde u} \otimes {\tilde p} $ is
expressed by a product of simple reflections as
\begin{equation}\label{hkt2}
{\tilde u} \otimes {\tilde p} =
(\sigma \otimes \sigma)S_{2} \cdots S_{n}S_{n+1}(u \otimes p).
\end{equation}
\end{theorem}
Of course the actual image ${\tilde p} \otimes {\tilde u}$ is
obtained by swapping the two tensor components of (\ref{hkt2}).

For $l$ suf\/f\/iciently large, f\/ix the element
$u = (u_1,\ldots, u_{n+1}) \in B^{1,l}$ according to $p \in B_1$ by
\begin{gather}
u_i = \varphi_i(S_{i+1}S_{i+2}\cdots S_{n+1}(p))\qquad  (2\le i \le n+1), \nonumber\\
u_1 =l-(u_2+\cdots + u_{n+1}).\label{k:uc}
\end{gather}
This is possible, i.e., $u_1\ge 0$, if $l$ is suf\/f\/iciently large.
We are going to show that substitution of~(\ref{k:uc}) into~(\ref{hkt2}) leads to
\begin{equation}\label{k:moku}
{\tilde u}  = u,\qquad
{\tilde p} = \sigma S_{2}S_{3}\cdots S_{n+1}(p).
\end{equation}
Since $p$ is $T^{(1)}_\infty$-evolvable
by the assumption of Theorem \ref{k:th:fac},
${\tilde u}  = u$ means that $u$ (\ref{k:uc}) is a~proper carrier and
\begin{equation}\label{k:tss}
T^{(1)}_\infty(p) = \sigma S_{2}S_{3}\cdots S_{n+1}(p).
\end{equation}
Let $r_i$ be the Weyl group simple ref\/lection acting on~$B$ component-wise.
Then, by using $r_i^2=\text{id}$ and
$\sigma r_2r_3\cdots r_{n+1}=\text{id}$,
the formula (\ref{k:tss}) is rewritten as
\begin{equation*}
T^{(1)}_\infty(p) = \overline{K}_{2}\overline{K}_3\cdots
\overline{K}_{n+1}(p), \quad
\overline{K}_a = r_{n+1}\cdots r_{a+1}r_a S_a r_{a+1}\cdots r_n r_{n+1}.
\end{equation*}
It is straightforward to check that the gauge transformed simple ref\/lection
$\overline{K}_a$ is equal to $K_a$ def\/ined by the procedure $(i)$--$(iv)$
explained before Theorem \ref{k:th:fac} by
using the signature rule explained in \cite[Section 2.3]{KTT}.

It remains to verify (\ref{k:moku}).
We illustrate the calculation for $n=3$,
for the general case is completely parallel.
Then (\ref{k:uc}) reads
(note $S_4=S_0, \varphi_4=\varphi_0, \varepsilon_4=\varepsilon_0$, etc.)
\begin{equation*}
u_4=\varphi_0(p),\qquad
u_3=\varphi_3(S_0(p)),\qquad
u_2 = \varphi_2(S_3S_0(p)).
\end{equation*}
We set
\begin{gather*}
u^{(1)}\otimes p^{(1)} = S_0(u \otimes p),\qquad
u^{(2)}\otimes p^{(2)} = S_3S_0(u \otimes p),\qquad
u^{(3)}\otimes p^{(3)} = S_2S_3S_0(u \otimes p),\\
\alpha_4=\varepsilon_0(p),\qquad
\alpha_3=\varepsilon_3(p^{(1)}),\qquad
\alpha_2=\varepsilon_2(p^{(2)}).
\end{gather*}
We know $u_1$ is suf\/f\/iciently large and
from $p \in B_1$, $\alpha_4 \ge u_4$ is valid.
To f\/irst compute $S_0(u\otimes p)$, we consider the 0-signature:
\begin{equation*}
u\otimes p = (u_1,u_2,u_3,u_4)\otimes p\quad 0:
\overbrace{--\cdots --}^{u_1}\overbrace{+++}^{u_4}
\overbrace{----}^{\alpha_4}\overbrace{+++}^{u_4}.
\end{equation*}
Thus $S_0$ actually acts as
$S_0 = \tilde{e}_0^{u_1-u_4}\otimes
\tilde{e}_0^{\alpha_4-u_4} = \tilde{e}_0^{u_1-u_4}\otimes S_0$,
leading to
\begin{equation*}
u^{(1)}\otimes p^{(1)} = (u_4,u_2,u_3,u_1)\otimes S_0(p)\quad 3:
\overbrace{--\cdots --}^{u_1}\overbrace{+++}^{u_3}
\overbrace{----}^{\alpha_3}\overbrace{+++}^{u_3},
\end{equation*}
where we have depicted the 3-signature on account of
$\varphi_3(p^{(1)}) = \varphi_3(S_0(p)) = u_3$.
{}From $p \in B_1$ we see that
$\alpha_3 = \varepsilon_3(p^{(1)}) = \varepsilon_3(S_0(p))
\ge \varphi_3(S_0(p)) = u_3$.
Thus the next $S_3$ actually acts as
$S_3 = \tilde{e}_3^{u_1-u_3}\otimes \tilde{e}_3^{\alpha_3-u_3}
= \tilde{e}_3^{u_1-u_3}\otimes S_3$, leading to
\begin{equation*}
u^{(2)}\otimes p^{(2)} = (u_4,u_2,u_1,u_3)\otimes S_3S_0(p)\quad 2:
\overbrace{--\cdots --}^{u_1}\overbrace{+++}^{u_2}
\overbrace{----}^{\alpha_2}\overbrace{+++}^{u_2},
\end{equation*}
where we have depicted the 2-signature on account of
$\varphi_2(p^{(2)}) = \varphi_2(S_3S_0(p)) = u_2$.
{}From $p \in B_1$ we see that
$\alpha_2 = \varepsilon_2(p^{(2)}) = \varepsilon_2(S_3S_0(p))
\ge \varphi_2(S_3S_0(p)) = u_2$.
Thus the next $S_2$ actually acts as
$S_2 = \tilde{e}_2^{u_1-u_2}\otimes \tilde{e}_2^{\alpha_2-u_2}
= \tilde{e}_2^{u_1-u_2}\otimes S_2$, leading to
\begin{equation*}
u^{(3)}\otimes p^{(3)} = (u_4,u_1,u_2,u_3)\otimes S_2S_3S_0(p).
\end{equation*}
Thus (\ref{hkt2}) gives
\begin{equation*}
{\tilde u}\otimes {\tilde p}=
\sigma((u_4,u_1,u_2,u_3))\otimes \sigma S_2S_3S_0(p)
=(u_1,u_2,u_3,u_4) \otimes \sigma S_2S_3S_0(p),
\end{equation*}
yielding (\ref{k:moku}).

\section{KKR bijection}\label{k:app:kkr}

\subsection{General remarks}
The original KKR bijection \cite{KKR, KR} is the one
between rigged conf\/igurations and Littlewood--Richardson tableaux.
Its ultimate generalization in type $A^{(1)}_n$
corresponding to
$B^{r_1,l_1}\otimes \cdots \otimes B^{r_1,l_L}$
is available in \cite{KSS}.
In the simple setting of this paper,
the Littlewood--Richardson tableaux are
in one to one correspondence with highest paths
in $B=(B^{1,1})^{\otimes L}$.
The KKR bijection in this paper means the one (\ref{k:kkr}) between
rigged conf\/igurations and those highest paths.
See \cite[Appendix C]{KSY}
for an exposition in a slightly more general setting
$B^{1, l_1}\otimes \cdots \otimes B^{1, l_L}$.

Recall that a rigged conf\/iguration
$(\mu, {\bf r})$ is a multiset of strings,
where a string is a triple
$\big(a, l^{(a)}_i, r^{(a)}_{i,\alpha}\big)$ consisting of
color, length and rigging.
See the right diagram in (\ref{k:rca0}).
A string $\big(a, l^{(a)}_i, r^{(a)}_{i,\alpha}\big)$ is
{\em singular} if $r^{(a)}_{i,\alpha} = p^{(a)}_i$, namely
if the rigging attains the allowed maximum (\ref{k:rcon}).
We also recall that one actually has to attach the data $L$ with
$(\mu, {\bf r})$ to specify the vacancy numbers
$p^{(a)}_i$ (\ref{k:pdef}).
Thus we write a rigged conf\/iguration as $(\mu, {\bf r})_L$.
We regard a highest path
$b_1 \otimes \cdots \otimes b_L \in (B^{1,1})^{\otimes L}$
as a word $b_1b_2\ldots b_L \in \{1,2,\ldots, n+1\}^L$.
The algorithms explained below obviously
satisfy the property (\ref{k:an}).

\subsection[Algorithm for $\phi$]{Algorithm for $\boldsymbol{\phi}$}

Given a highest path $b_1\ldots b_L$,
we construct the
rigged conf\/iguration $\phi(b_1\ldots b_L) = (\mu, {\bf r})_L$
inductively with respect to $L$.
When $L=0$, we understand that $\phi(\cdot)$
is the array of empty partitions.
Suppose that
$\phi(b_1\ldots b_L) = (\mu, {\bf r})_L$ has been obtained.
Denote $b_{L+1} \in \{1,\ldots, n+1\}$ simply by $d$.
We are to construct
$(\mu', {\bf r}')_{L+1}=\phi(b_1\ldots b_L d)$ from
$(\mu, {\bf r})_L$ and $d$.
If $d=1$, then
$(\mu', {\bf r}')_{L+1}=(\mu, {\bf r})_{L+1}$ meaning that
nothing needs to be done except increasing the
vacancy numbers $p^{(1)}_i$ by one.
(Recall in (\ref{k:pdef})
that only $p^{(1)}_i$ depends on $L$.)
Suppose $d \ge 2$.
\begin{enumerate}\itemsep=0pt

\item[$(i)$]
Set $\ell^{(d)}=\infty$.
For $c=d-1, d-2, \ldots, 1$ in this order, proceed as follows.
Find the color $c$ singular string whose length $\ell^{(c)}$
is largest within the condition $\ell^{(c)}\le \ell^{(c+1)}$.
If there are more than one such strings, pick any one of them.
If there is no such string with color~$c$, set $\ell^{(c)} = 0$.
Denote these selected strings by
$\big(c, \ell^{(c)}, r^{(c)}_\ast\big)$ with $c=d-1,d-2,\ldots, 1$,
where it is actually void when $\ell^{(c)}=0$.

\item[$(ii)$]
 Replace the selected string
$\big(c, \ell^{(c)}, r^{(c)}_\ast\big)$ by
$\big(c, \ell^{(c)}+1, r^{(c)}_\bullet\big)$
for all $c=d-1,d-2,\ldots, 1$ leaving the other strings unchanged.
Here the new rigging $r^{(c)}_\bullet$ is to be f\/ixed so that
the enlarged string $\big(c, \ell^{(c)}+1, r^{(c)}_\bullet\big)$ becomes
singular with respect to
the resulting new conf\/iguration $\mu'$ and $L+1$.
\end{enumerate}

The algorithm is known to be well-def\/ined and
the resulting object gives the sought rigged conf\/iguration
$(\mu', {\bf r}')_{L+1} =\phi(b_1\ldots b_L d)$.

\subsection[Algorithm for $\phi^{-1}$]{Algorithm for $\boldsymbol{\phi^{-1}}$}

Given a rigged conf\/iguration $(\mu, {\bf r})_L$,
we construct a highest path
$b_1\ldots b_L = \phi^{-1}((\mu, {\bf r})_L)$
inductively with respect to $L$.
We are to determine $d (=b_L) \in \{1,\ldots, n+1\}$
and $(\mu', {\bf r}')_{L-1}$ such that
$\phi^{-1}((\mu, {\bf r})_L) =
\phi^{-1}((\mu', {\bf r}')_{L-1}) \,d$.
\begin{enumerate}\itemsep=0pt

\item[$(i)$]
Set $\ell^{(0)}=1$.
For $c=1,2, \ldots, n$ in this order, proceed as follows until stopped.
Find the color $c$ singular string whose length $\ell^{(c)}$
is smallest within the condition $\ell^{(c-1)}\le \ell^{(c)}$.
If there are more than one such strings, pick any one of them.
If there is no such string with color $c$,  set $d=c$ and stop.
If $c=n$ and such a color $n$ string still exists, set $d=n+1$
and stop.
Denote these selected strings by
$\big(c, \ell^{(c)}, r^{(c)}_\ast\big)$ with $c=1,2, \ldots, d-1$.

\item[$(ii)$]
Replace the selected string
$\big(c, \ell^{(c)}, r^{(c)}_\ast\big)$ by
$\big(c, \ell^{(c)}-1, r^{(c)}_\bullet\big)$
for all $c=1,2, \ldots, d-1$ leaving the other strings unchanged.
When $\ell^{(c)}=1$, this means that the length one string
is to be eliminated.
The new rigging $r^{(c)}_\bullet$ is to be f\/ixed so that
the shortened string $\big(c, \ell^{(c)}-1, r^{(c)}_\bullet\big)$ becomes
singular with respect to the resulting new conf\/iguration $\mu'$ and $L-1$.
\end{enumerate}
For empty rigged conf\/iguration, we understand that
$\phi^{-1}((\varnothing, \varnothing)_L) =
\phi^{-1}((\varnothing, \varnothing)_{L-1}) \,1
= \cdots = \overbrace{11\ldots 1}^L$.
The algorithm is known well-def\/ined and ends up with
the empty rigged conf\/iguration at $L=0$.
The resulting sequence gives the sought highest path
$b_1\ldots b_L = \phi^{-1}((\mu, {\bf r})_{L})$.

\begin{example}\label{k:ex:kkr}
Let us demonstrate $\phi^{-1}$ along an $L=8$ example.
\begin{equation}\label{k:samp}
\begin{picture}(200,45)(-130,42)
\setlength{\unitlength}{0.5mm}

\multiput(0,0)(0,-10){2}{
\put(10,60){\line(1,0){20}}}

\multiput(0,0)(0,10){2}{
\put(10,30){\line(1,0){10}}}

\put(10,30){\line(0,1){30}}
\put(20,30){\line(0,1){30}}
\put(30,50){\line(0,1){10}}

\put(2,53){1}\put(2,38){3}
\put(35,52){1}
\put(23,42){2}
\put(23,32){1}

\multiput(-65,0)(0,0){1}{
\put(117,52){1}\put(140,52){0}

\multiput(25,0)(0,0){1}{
\put(100,60){\line(1,0){10}}
\put(100,50){\line(1,0){10}}
\put(100,50){\line(0,1){10}}
\put(110,50){\line(0,1){10}}}}

\put(93,40){$\overset{\phi^{-1}}{\longmapsto}
\quad 11213122$.}

\end{picture}
\end{equation}
This is the bottom right case of Example \ref{t:ex:sept2_1}.
The numbers on the left of Young diagrams are vacancy numbers,
which are subsidiary data but convenient to keep track of.
$$
\begin{picture}(455,130)(0,-37)
\setlength{\unitlength}{0.5mm}
\multiput(0,-10)(0,0){1}{
\put(35,70){$L=8$}
\multiput(0,0)(0,-10){2}{

\put(10,60){\line(1,0){20}}}

\multiput(0,0)(0,10){2}{
\put(10,30){\line(1,0){10}}}

\put(10,30){\line(0,1){30}}
\put(20,30){\line(0,1){30}}
\put(30,50){\line(0,1){10}}

\put(2,53){1}\put(2,38){3}
\put(35,52){1}
\put(23,42){2}
\put(23,32){1}

\multiput(-65,0)(0,0){1}{
\put(117,52){1}\put(140,52){0}

\multiput(25,0)(0,0){1}{
\put(100,60){\line(1,0){10}}
\put(100,50){\line(1,0){10}}
\put(100,50){\line(0,1){10}}
\put(110,50){\line(0,1){10}}}}
}

\put(90,35){$\overset{2}{\longmapsto}$}

\multiput(110,-10)(0,0){1}{
\put(30,70){$L=7$}
\multiput(0,0)(0,-10){2}{
\put(10,60){\line(1,0){10}}}

\multiput(0,0)(0,10){2}{
\put(10,30){\line(1,0){10}}}

\put(10,30){\line(0,1){30}}
\put(20,30){\line(0,1){30}}

\put(2,42){2}

\put(23,52){2}
\put(23,42){2}
\put(23,32){1}

\multiput(-70,0)(0,0){1}{
\put(117,52){1}\put(140,52){0}

\multiput(25,0)(0,0){1}{
\put(100,60){\line(1,0){10}}
\put(100,50){\line(1,0){10}}
\put(100,50){\line(0,1){10}}
\put(110,50){\line(0,1){10}}}}
}

\put(200,35){$\overset{2}{\longmapsto}$}

\multiput(220,-10)(0,0){1}{
\put(30,70){$L=6$}
\multiput(0,0)(0,-10){3}{
\put(10,60){\line(1,0){10}}}

\put(10,40){\line(0,1){20}}
\put(20,40){\line(0,1){20}}

\put(2,47){3}

\put(23,52){2}
\put(23,42){1}

\multiput(-70,0)(0,0){1}{
\put(117,52){0}\put(140,52){0}

\multiput(25,0)(0,0){1}{
\put(100,60){\line(1,0){10}}
\put(100,50){\line(1,0){10}}
\put(100,50){\line(0,1){10}}
\put(110,50){\line(0,1){10}}}}
}

\put(0,-17){$\overset{1}{\longmapsto}$}

\multiput(20,-70)(0,0){1}{
\put(30,70){$L=5$}
\multiput(5,0)(0,-10){3}{
\put(10,60){\line(1,0){10}}}

\put(15,40){\line(0,1){20}}
\put(25,40){\line(0,1){20}}

\put(7,47){2}

\put(28,52){2}
\put(28,42){1}

\multiput(-75,0)(0,0){1}{
\put(117,52){0}\put(140,52){0}

\multiput(25,0)(0,0){1}{
\put(100,60){\line(1,0){10}}
\put(100,50){\line(1,0){10}}
\put(100,50){\line(0,1){10}}
\put(110,50){\line(0,1){10}}}}
}

\put(100,-17){$\overset{3}{\longmapsto}$}

\multiput(115,-70)(0,0){1}{
\put(20,70){$L=4$}
\multiput(5,0)(0,-10){2}{
\put(10,60){\line(1,0){10}}}

\put(15,50){\line(0,1){10}}
\put(25,50){\line(0,1){10}}

\put(7,52){2} \put(28,52){1}

\put(40,52){$\varnothing$}
}

\put(170,-17){$\overset{1}{\longmapsto}$}

\multiput(185,-70)(0,0){1}{
\put(20,70){$L=3$}
\multiput(5,0)(0,-10){2}{
\put(10,60){\line(1,0){10}}}

\put(15,50){\line(0,1){10}}
\put(25,50){\line(0,1){10}}

\put(7,52){1} \put(28,52){1}

\put(40,52){$\varnothing$}
}

\put(240,-17){$\overset{2}{\longmapsto}$}

\multiput(245,-70)(0,0){1}{
\put(20,70){$L=2$}

\put(20,52){$\varnothing$}

\put(40,52){$\varnothing$}
}
\end{picture}
$$
Here we have exhibited the number $d$ in the algorithm (i) over the
arrows.
After this, we only get $1$ twice.
Thus we get the sequence
$\overset{2}{\mapsto}\,\overset{2}{\mapsto}\,
\overset{1}{\mapsto}\,\overset{3}{\mapsto}\,
\overset{1}{\mapsto}\,\overset{2}{\mapsto}\,
\overset{1}{\mapsto}\,\overset{1}{\mapsto}$
to reach the trivial rigged conf\/iguration.
Reading these numbers backwards we obtain (\ref{k:samp}).
\end{example}

\section{Proof of Theorem \ref{k:th:rc}}\label{k:app:proof1}

For $p=b_1 \otimes \cdots \otimes b_L \in B$,
def\/ine another version of energy $\hat{E}^{(r)}_l(p)=e_1+\cdots + e_L$
with the diagram (\ref{k:hone}) by taking $v = u^{r,l}$, the
highest element of $B^{r,l}$.
$u^{r,l}$ is the semistandard tableau
whose entries in the $j$th row are all $j$.
For example, $u^{2,3}={111 \atop 222}$.
Set
\begin{gather}
\hat{{\mathcal P}}(\mu)  =\Bigg\{p \in {\mathcal P}\mid
\hat{E}^{(a)}_l(p) = \sum_{i=1}^{g_a}\min\big(l,l^{(a)}_i\big)m^{(a)}_i\Bigg\},
\label{k:pmuh}\\
\hat{{\mathcal P}}_+(\mu)  =
\big\{p \in \hat{{\mathcal P}}(\mu)\mid
p: \text{highest} \big\}.
\label{k:hiph}
\end{gather}
These def\/initions resemble (\ref{k:pmu}) and (\ref{k:hip}).
The dif\/ferences are that here enters $\hat{E}^{(r)}_l$ instead of $E^{(r)}_l$,
and there is no requirement that $p$ be evolvable in (\ref{k:pmuh}).

\begin{theorem}[\cite{Sa}]\label{k:th:sa}
The KKR bijection \eqref{k:kkr} can be restricted into a finer bijection:
\begin{equation*}
\phi : \ \hat{{\mathcal P}}_+(\mu) \rightarrow \text{\rm RC}(\mu).
\end{equation*}
\end{theorem}

\begin{lemma}\label{k:le:ee}
The equality
$E^{(r)}_l(p) = \hat{E}^{(r)}_l(p)$ is valid for any
$p \in {\mathcal P}_+(\mu)$.
\end{lemma}

\begin{proof}
Take $p \in {\mathcal P}_+(\mu)$.
Then by Lemma \ref{k:lem:ext}, the rigged conf\/iguration for
$p\otimes p$ has the conf\/iguration obtained from $\mu$ by
replacing $m^{(a)}_i$ with $2m^{(a)}_i$.
Thus from (\ref{k:pmuh}) and Theorem \ref{k:th:sa},
one has $\hat{E}^{(r)}_l(p\otimes p) = 2\hat{E}^{(r)}_l(p)$.
On the other hand, def\/ine
$v^{r,l} \in B^{r,l}$ (and $\hat{p}$) by
$u^{r,l}\otimes p \simeq \hat{p}\otimes v^{r,l}$.
Then, using Lemma \ref{k:lem:ext} in this paper
and Lemma 8.5 in \cite{KSS},
one can show $v^{r,l} \otimes p \simeq p' \otimes v^{r,l}$
for some $p' \in B$.
Since $p$ is evolvable, we can claim that
$v^{r,l}$ works as a carrier to def\/ine $E^{(r)}_l(p)$.
Then from the relation
$u^{r,l}\otimes p \otimes p \simeq
\hat{p} \otimes v^{r,l} \otimes p \simeq
\hat{p} \otimes p' \otimes v^{r,l}$,
the energy $\hat{E}^{(r)}_l$ associated to
$p \otimes p$ is given by
$\hat{E}^{(r)}_l(p\otimes p) = \hat{E}^{(r)}_l(p) + E^{(r)}_l(p)$.
Comparing this with
$\hat{E}^{(r)}_l(p\otimes p) = 2\hat{E}^{(r)}_l(p)$,
we get $E^{(r)}_l(p) = \hat{E}^{(r)}_l(p)$.
\end{proof}

\begin{lemma}\label{k:le:pp}
There is a natural inclusion
${\mathcal P}_+(\mu) \hookrightarrow \hat{\mathcal P}_+(\mu)$.
\end{lemma}

\begin{proof}
Due to Lemma \ref{k:le:ee}, the dif\/ference of
$\hat{\mathcal P}_+(\mu)$ (\ref{k:hiph})
and ${\mathcal P}_+(\mu)$  (\ref{k:hip})
is only the extra condition that $p$ be evolvable in (\ref{k:pmu}).
\end{proof}

We expect that actually
${\mathcal P}_+(\mu) = \hat{\mathcal P}_+(\mu)$ holds,
namely, all the highest paths are evolvable.
However we do not need this fact in this paper.

\begin{proof}[Proof of Theorem \ref{k:th:rc}.]
Combine Theorem \ref{k:th:sa} and Lemma \ref{k:le:pp}.
\end{proof}

\subsection*{Acknowledgements}

A.K. thanks Rei Inoue,
Masato Okado, Reiho Sakamoto, Mark Shimozono,
Alexander Veselov, Yasuhiko Yamada for discussion, and
Claire Gilson, Christian Korf\/f and Jon Nimmo for
a warm hospitality during the conference,
{\it Geometric Aspects of Discrete and Ultra-Discrete Integrable Systems},
March~30~-- April~3, 2009, Glasgow, UK.
This work is partially supported by Grand-in-Aid for Scientif\/ic
Research JSPS No. 21540209.

\addcontentsline{toc}{section}{References}
\LastPageEnding

\end{document}